\documentclass[12pt,  reqno, oneside]{amsart}      
\usepackage{geometry}                       
\usepackage{amsthm}
\usepackage{bbm}
\usepackage{dsfont}
\usepackage{graphicx}               
\usepackage{amssymb}
\usepackage{cite}
\newcommand{\N}{\mathbb{N}}
\newcommand{\R}{\mathbb{R}}

\newcommand{\Heis}{\mathbb{H}}
\newcommand{\eps}{\varepsilon}
\newcommand{\Lmeas}{\mathcal{L}}
\newcommand{\Jac}{{\rm Jac}}
\newcommand{\Ent}{{\rm Ent}}
\newcommand{\dd}{{\rm d}}

\usepackage{color}
\usepackage[unicode]{hyperref}
\hypersetup{
	colorlinks = true,%
	citecolor = [rgb]{0.0,0.0,0.9},
	filecolor=black,%
	linkcolor = [rgb]{0.65,0.0,0.0},%
	anchorcolor = red,
	pagecolor = red,
	urlcolor= [rgb]{0.65,0.0,0.0}
}

\usepackage{amssymb,amsmath}
\usepackage{graphicx}
\usepackage{mathrsfs}
\usepackage{float}
\usepackage[active]{srcltx}

\numberwithin{equation}{section}

\newtheorem{proposition}{Proposition}[section]
\newtheorem{theorem}{Theorem}[section]
\newtheorem{lemma}{Lemma}[section]
\newtheorem{corollary}{Corollary}[section]

\newtheorem{remark}{Remark}[section]

\makeatletter
\def\l@section{\@tocline{1}{0pt}{1pc}{}{}}
\def\l@subsection{\@tocline{2}{0pt}{1pc}{4.6em}{}}
\def\l@subsubsection{\@tocline{3}{0pt}{1pc}{7.6em}{}}
\renewcommand{\tocsection}[3]{%
  \indentlabel{\@ifnotempty{#2}{\makebox[2.3em][l]{%
    \ignorespaces#1 #2.\hfill}}}#3}
\renewcommand{\tocsubsection}[3]{%
  \indentlabel{\@ifnotempty{#2}{\hspace*{2.3em}\makebox[2.3em][l]{%
    \ignorespaces#1 #2.\hfill}}}#3}
\renewcommand{\tocsubsubsection}[3]{%
  \indentlabel{\@ifnotempty{#2}{\hspace*{4.6em}\makebox[3em][l]{%
    \ignorespaces#1 #2.\hfill}}}#3}
\makeatother

\title[Geometric inequalities on Heisenberg groups]
{Geometric inequalities on Heisenberg groups}

\author{Zolt\'an M. Balogh, Alexandru Krist\'aly, and Kinga Sipos}	

\thanks{Z. M. Balogh was
	supported by the Swiss National Science Foundation, Grant Nr. {200020\_146477}.  A. Krist\'aly  was supported by  the STAR-UBB Institute.
	K. Sipos was supported by ERC Marie-Curie Research and Training Network MANET}

\makeatletter
\renewcommand{\@biblabel}[1]{[#1]\hfill}
\makeatother

\begin{document}

	\begin{abstract}
		We establish geometric inequalities in the sub-Riemannian setting of the Heisenberg group $\mathbb H^n$. Our results include a natural sub-Riemannian version of the celebrated curvature-dimension condition of Lott-Villani and Sturm and also a geodesic version of the  Borell-Brascamp-Lieb inequality akin to the one obtained by Cordero-Erausquin, McCann and Schmuc\-kenschl\"ager. The latter statement implies  sub-Riemannian versions of the geodesic Pr\'ekopa-Leindler and Brunn-Minkowski inequalities. The proofs are based on optimal mass transportation and Riemannian approximation of $\mathbb H^n$ developed by Ambrosio and Rigot. These results refute a gene\-ral point of view, according to which no  geometric inequalities can be derived by optimal mass transportation on singular spaces. 
	\end{abstract}
	
	\vspace*{-1.8cm}

	\maketitle
	
	\begin{center}\it Dedicated to Hans Martin Reimann on the occasion of his 75th birthday
		\\and to Cristian Guti\'errez on the occasion of  his 65th birthday.
	\end{center}

	\tableofcontents

	\vspace*{0.5cm}
	
	
	\medskip
	\medskip
	

	\vspace*{1.6cm}
	\section{Introduction and main results}
	\vspace*{0.7cm}
	\subsection{General background and motivation}\label{subsection-1.1}
	Due to the seminal papers by Lott and Villani \cite{LV} and Sturm \cite{Sturm1, Sturm2}, metric measure spaces with generalized lower Ricci curvature bounds support various geometric and functional inequalities including Borell-Brascamp-Lieb, Brunn-Minkowski, Bishop-Gromov  inequalities. A basic assumption for these results is the famous curvature-dimension condition $\textsf{CD}(K,N)$ which --in the case of a Riemannian manifold $M$--,  represents  the lower bound $K\in \mathbb R$ for the Ricci curvature on $M$ and the upper bound $N\in \mathbb R$ for the dimension of $M$, respectively. It is a fundamental  question whether the method used in  \cite{LV}, \cite{Sturm1, Sturm2}, based on optimal mass transportation works in the setting of {\it singular spaces} with no apriori lower curvature bounds. A large class of such spaces are the sub-Riemannian geometric structures or Carnot-Carath\'eodory geometries, see Gromov \cite{Gromov}.
	
	During the last decade considerable effort has been made to  establish geometric and functional inequalities on sub-Riemannian spaces. The quest for Borell-Brascamp-Lieb and Brunn-Minkowski type inequalities became a hard nut to crack even on simplest sub-Riemannian setting such  as the Heisenberg group $\mathbb H^n$ endowed  with the usual Carnot-Carath\'eo\-do\-ry metric $d_{CC}$ and $\mathcal L^{2n+1}$-measure. 
	One of the reasons for  this is that although there is a good first order Riemannian approximation (in the pointed Gromov-Hausdorff sense) of the sub-Riemannian metric structure of the Heisenberg group $\mathbb H^n$, there is no uniform lower bound  on the Ricci curvature in these approximations (see e.g.  Capogna, Danielli, Pauls and Tyson \cite[Section 2.4.2]{Capogna}); indeed, at 
	every point of $\mathbb H^n$ there is a Ricci curvature whose limit is $-\infty$  in the Riemannian approximation.  The lack of uniform lower Ricci bounds prevents a straightforward extension of the Riemannian Borell-Brascamp-Lieb and Brunn-Minkowski inequalities of Cordero-Erausquin, McCann and Schmuckenschl\"ager \cite{McCann} to the setting 
	of the Heisenberg group. Another serious warning is attributed to Juillet
	\cite{Juillet-IMNR} who proved that both the Brunn-Minkowski inequality and the curvature-dimension condition  $\textsf{CD}(K,N)$  fail on $(\mathbb H^n,d_{CC},\mathcal
	L^{2n+1})$ for every choice of $K$ and $N$.

	These facts tacitly established the view according to which there are no entropy-convexity and Borell-Brascamp-Lieb type inequalities on singular spaces such as the Heisenberg groups. The purpose of this paper is to deny this paradigm.	 Indeed, we show that the method of  optimal mass transportation is powerful enough to yield good results even in the absence of lower curvature bounds. By using convergence results for optimal transport maps in the  Riemannian approximation of $\mathbb H^n$ due to Ambrosio and Rigot \cite{AR} we are able to introduce the correct sub-Riemannian geometric quantities which can replace the lower curvature bounds and can be successfully used to establish geodesic Borell-Brascamp-Lieb, Pr\'ekopa-Leindler, Brunn-Minkowski and  entropy inequalities on the Heisenberg group $\mathbb H^n$. The main statements from the papers of Figalli and Juillet \cite{FJ} and Juillet \cite{Juillet-IMNR} will appear as special cases of our results.  
	
	Before stating our results we shortly recall  the aforementioned geometric inequalities of Borell-Brascamp-Lieb and the curvature dimension condition  of Lott-Sturm-Villani and indicate their behavior in the sub-Riemannian setting of  Heisenberg groups. 
	
	\newpage
	
	\subsection{An overview of geometric inequalities} The classical Borell-Brascamp-Lieb inequality in $\R^n$ states that for any
	fixed $s \in (0,1)$, $p \geq -\frac{1}{n}$ and integrable functions  $f, g, h: \R^n \to
	[0, \infty)$ which satisfy 
	\begin{eqnarray}\label{elso-BBL-feltetel}
	h((1-s)x + sy)
	\geq M_s^p \left( f(x), g(y) \right) \quad \mbox{ for all } x, y \in \R^n,
	\end{eqnarray}
	one has
	\begin{eqnarray*}
		\int_{\R^n} h \geq M_s^{\frac{p}{1 + np}} \left( \int_{\R^n} f,
		\int_{\R^n} g \right).
	\end{eqnarray*}
	Here and in the sequel, for every $s\in (0,1)$, $p\in \mathbb R\cup
	\{\pm\infty\}$ and $a,b\geq 0$, we consider  the $p$-mean
	$$M_s^p(a,b)=\left\{\begin{array}{lll}
	\left( (1-s)a^p + s b^p \right)^{1/p} &\mbox{if} &  ab\neq 0; \\
	0 &\mbox{if} &  ab=0,
	\end{array}\right.$$
	with the  conventions 
	$M_s^{-\infty}(a,b)=\min\{a,b\}$, and $M_s^0(a,b)=a^{1-s}b^s,$ and $M_s^{+\infty}(a,b)=\max\{a,b\}$ if $ab\neq 0$ and  $M_s^{+\infty}(a,b)=0$ if $ab=0.$
	The Borell-Brascamp-Lieb inequality reduces to the
	Pr\'ekopa-Leindler inequality for $p=0$, which in turn implies the Brunn-Minkowski
	inequality
	$$    \mathcal L^{n}((1-s)A+sB)^\frac{1}{n}\geq (1-s)\mathcal
	L^{n}(A)^\frac{1}{n}+s\mathcal L^{n}(B)^\frac{1}{n},
	$$
	where $ A$ and $B$ are positive and finite measure subsets of $\R^{n}$, and $\mathcal
	L^n$ denotes the $n$-dimensional Lebesgue measure. For a  comprehensive survey on geometric inequalities in  $\R^n$ and their applications to isoperimetric problems, sharp Sobolev inequalities and convex geometry, we refer to Gardner \cite{Gardner}.

	In his Ph.D. Thesis, McCann \cite[Appendix D]{McCann-PhD} (see also \cite{McCann_Adv_Math}) presented an optimal mass transportation approach to Pr\'ekopa-Leindler, Brunn-Minkowski and Brascamp-Lieb inequalities in the Euclidean setting. This pioneering idea led to the extension of    
	a geodesic version of the Borell-Brascamp-Lieb inequality on complete Riemannian manifolds via optimal mass transportation, established by Cordero-Erausquin, McCann and Schmuckenschl\"ager \cite{McCann}.  Closely related to the Borell-Brascamp-Lieb inequalities on Riemannian manifolds is the convexity of the entropy functional \cite{McCann}. The latter fact served as the starting point of the work of Lott and Villani \cite{LV} and Sturm \cite{Sturm1, Sturm2} who initiated independently the synthetic study of Ricci curvature on  metric measure spaces by introducing the  curvature-dimension condition  $\textsf{CD}(K,N)$ for $K\in \mathbb R$ and $N\geq 1$. Their approach is based on the effect of the curvature of the space encoded in the reference distortion coefficients
	{\small
		$$\tau_s^{K,N}(\theta)=\left\{
		\begin{array}{lll}
		+\infty, & {\rm if} & K\theta^2\geq (N-1)\pi^2;
		\\ s^\frac{1}{N}\left(\sin\left(\sqrt{\frac{K}{N-1}}s\theta\right)\big/\sin\left(\sqrt{\frac{K}{N-1}}\theta\right)\right)^{1-\frac{1}{N}},& {\rm if} &
		0<K\theta^2<(N-1)\pi^2;\\
		s, & {\rm if} & K\theta^2=0;\\
		s^\frac{1}{N}\left(\sinh\left(\sqrt{-\frac{K}{N-1}}s\theta\right)\big/\sinh\left(\sqrt{-\frac{K}{N-1}}\theta\right)\right)^{1-\frac{1}{N}},&
		{\rm if} & K\theta^2<0,
		\end{array}
		\right.$$}
	where $s\in (0,1)$, see e.g. Sturm \cite{Sturm2} and Villani \cite{Villani1}.  	To be more precise,   
	let $(M,d,\textsf{m})$ be a metric measure space,  
	$K\in \mathbb R$ and $N\geq 1$ be fixed, $\mathcal P_2(M,d)$ be the usual Wasserstein space, and ${\rm Ent}_{N'}(\cdot|\textsf{m}):\mathcal P_2(M,d)\to \mathbb R$ be the R\'enyi entropy functional given by  
	\begin{eqnarray} \label{EntDef}
	{\rm \Ent}_{N'}(\mu | \textsf{m}) = -\displaystyle\int_M \rho^{1-\frac{1}{N'}} \dd \textsf{m},
	\end{eqnarray} 
	where $\rho$ is the density function of $\mu$ w.r.t. $\textsf{m},$ and $N'\geq N.$ The metric  measure space $(M,d,\textsf{m})$ satisfies the 	curvature-dimension condition $\textsf{CD}(K,N)$ for $K\in \mathbb R$ and $N\geq 1$ if and only if for every $\mu_0,\mu_1\in \mathcal P_2(M,d)$ there exists an optimal coupling $q$ of $\mu_0=\rho_0 \textsf{m}$ and $\mu_1=\rho_1 \textsf{m}$ and a geodesic $\Gamma:[0,1]\to \mathcal P_2(M,d)$ joining $\mu_0$ and $\mu_1$ such that for all $s\in [0,1]$ and $N'\geq N$, 
	$${\rm Ent}_{N'}(\Gamma(s)|\textsf{m})\leq -\int_{M\times M}\left[\tau_{1-s}^{K,N'}(d(x_0,x_1))\rho_0(x_0)^{-\frac{1}{N'}}+\tau_s^{K,N'}(d(x_0,x_1))\rho_1(x_1)^{-\frac{1}{N'}}\right]{\rm d}q(x_0,x_1).$$
	It turns out that a Riemannian (resp. Finsler) manifold $(M,d,\textsf{m})$ satisfies the condition $\textsf{CD} ( K, N )$ if and only if the Ricci curvature on $M$ is not smaller than  $K$ and the dimension of $M$ is not greater than $N$, where $d$ is the natural metric on $M$ and $\textsf{m}$ is the canonical Riemannian (resp. Busemann-Hausdorff) measure on $M,$ see Sturm \cite{Sturm2} and Ohta \cite{Ohta}.

	Coming back to the Borell-Brascamp-Lieb inequality in curved spaces, e.g., when $(M,d,\textsf{m})$ is a complete $N$-dimensional Riemannian manifold, we have to replace the convex combination $(1-s)x + sy$ in (\ref{elso-BBL-feltetel})  by the set of $s$-intermediate points $Z_s(x,y)$ between $x$ and $y$ w.r.t. the Riemannian metric $d$ on $M$ defined by
	$$ Z_s(x,y)=  \{ z \in M : d(x,z) = s d(x,y),\
	d(z,y) = (1-s) d(x,y)\}.$$
	With this notation, we can state the result of Cordero-Erausquin, McCann and Schmuckenschl\"ager \cite{McCann} (see also Bacher \cite{Bacher}), as the  Borell-Brascamp-Lieb inequality $\textsf{BBL}(K,N)$ on $(M,d,\textsf{m})$ which holds if and only if for all $s \in (0,1)$, $p\geq -\frac{1}{N}$ and integ\-rable functions  $f, g, h: M \to
	[0, \infty)$ satisfying 
	\begin{eqnarray}\label{masodik-BBL-feltetel}
	h(z)
	\geq M_s^p \left( \frac{f(x)}{\left(\tilde\tau_{1-s}^{K,N}(d(x,y))\right)^N},\frac{g(y)}{\left(\tilde\tau_{s}^{K,N}(d(x,y))\right)^N} \right) \quad \mbox{ for all } x, y \in
	M,\ z\in Z_s(x,y),
	\end{eqnarray}
	one has
	\begin{eqnarray*}
		\int_{M} h {\rm d}\textsf{m} \geq M_s^{\frac{p}{1 + Np}} \left( \int_{M} f{\rm d}\textsf{m},
		\int_{M} g{\rm d}\textsf{m} \right),
	\end{eqnarray*}
	where $\tilde\tau_{s}^{K,N}=s^{-1}\tau_{s}^{K,N}.$ We would like to emphasize the fact that in \cite{McCann} the main ingredient is provided by a weighted Jacobian determinant inequality satisfied by the optimal transport interpolant map.  
	
	It turns out, even in the more general setting of non-branching geodesic metric spaces, that both $\textsf{CD}(K,N)$ and $\textsf{BBL}(K,N)$ imply the geodesic Brunn-Minkowski inequality $\textsf{BM}(K,N)$, see Bacher \cite{Bacher}, i.e., if $(M,d,\textsf{m})$ is such a space, for Borel sets $A,B\subset M$ with $\textsf{m}(A)\neq 0\neq \textsf{m}(B)$ and $s\in(0,1)$,
	\begin{equation}\label{BM-CD-version}
	\textsf{m}(Z_s(A,B))^\frac{1}{N}\geq \tau_{1-s}^{K,N}(\theta_{A,B})\textsf{m}(A)^\frac{1}{N}+\tau_{s}^{K,N}(\theta_{A,B})\textsf{m}(B)^\frac{1}{N}.
	\end{equation}
	Here $Z_s(A,B)$ is the set of $s$-intermediate points between the
	elements of  the sets $A$ and $B$ w.r.t. the metric $d$, defined by
	$Z_s(A,B) = \bigcup_{(x,y) \in A \times B} Z_s(x,y),$
	and
	\begin{eqnarray*}
		\theta_{A,B}= \left\{
		\begin{array}{lll}
			\inf_{(x,y)\in A\times B}d(x,y)
			\  &\mbox{if} &  K\geq 0; \\
			\sup_{(x,y)\in A\times B}d(x,y)
			\  &\mbox{if} &  K< 0.
		\end{array}\right.
	\end{eqnarray*}
	
	As we already pointed out, Juillet
	\cite{Juillet-IMNR} proved that the Brunn-Minkowski inequality $\textsf{BM}(K,N)$ fails on $(\mathbb H^n,d_{CC},\mathcal
	L^{2n+1})$ for every choice of $K$ and $N$; therefore,  both $\textsf{CD}(K,N)$ and $\textsf{BBL}(K,N)$ fail too.  In fact, a closer investigation shows that  the failure of these  inequalities  on $\mathbb H^n$ is
	not surprising: indeed,  the distortion coefficient 
	$\tau_s^{K,N}$ is a 'pure Riemannian' object coming from the behavior of Jacobi fields along geodesics in Riemannian space forms. More quantitatively, since certain Ricci curvatures tend to $-\infty$ in the Riemannian approximation of the first Heisenberg group $\mathbb H^1$  (see   Capogna, Danielli, Pauls and Tyson \cite[Section 2.4.2]{Capogna}) and $\lim_{K\to-\infty}\tau_s^{K,N}(\theta)=0$ for every $s\in (0,1)$ and $\theta>0$,  some Riemannian quantities blow up  and they fail to capture the subtle sub-Riemannian metric structure of the Heisenberg group. In particular, assumption (\ref{masodik-BBL-feltetel}) in $\textsf{BBL}(K,N)$ degenerates to an impossible condition. 
	
	On the other hand, there is a positive effect in the Riemannian approximation (see \cite[Section 2.4.2]{Capogna}) that would be unfair to conceal. It turns out namely, that the two remaining Ricci curvatures in $\mathbb H^1$ will blow up to $+\infty$ in the Riemannian approximation scheme. This can be interpreted as a sign of hope for a certain cancellation that could save the day at the end. This will be indeed the case: appropriate geodesic versions of Borell-Brascamp-Lieb and Brunn-Minkowski inequalities still hold on the Heisenberg group as we show in the sequel.

	\subsection{Statement of main results} According to Gromov \cite{Gromov}, the Heisenberg group $\Heis^n$ with its sub-Riemannian, or Carnot-Carath\'eodory metric, can be seen as the simplest prototype of a singular space. In this paper we shall use a model of  $\Heis^n$ that is identified with its Lie
	algebra $\mathbb R^{2n+1} \simeq \mathbb C^n \times \mathbb R$ via canonical
	exponential coordinates. At this point we just recall the bare minimum that is needed of the metric structure of $\Heis^n$ in order to state our results. In the next section we present a more detailed exposition of the Heisenberg geometry, its Riemannian approximation and the connection between their optimal mass transportation maps. We denote a
	point in $\Heis^n$ by $x = (\xi, \eta, t) = (\zeta, t)$, where $\xi
	= (\xi_1, \ldots, \xi_n) \in \R^n$, $\eta = (\eta_1, \ldots, \eta_n)
	\in \R^n$, $t \in \R,$ and we identify the pair $(\xi, \eta)$ with
	$\zeta \in \mathbb C^n$ having coordinates $\zeta_j = \xi_j + i \eta_j$ for all $j = 1, \ldots, n$.
	The correspondence with its Lie algebra through the exponential
	coordinates induces the  group law
	$$(\zeta, t) \cdot (\zeta', t') = \left( \zeta + \zeta', t + t' + 2{\rm Im} \langle\zeta , {\zeta'}\rangle \right), \quad \forall (\zeta, t),\ (\zeta', t') \in \mathbb C^n \times \mathbb R,$$
	where ${\rm Im}$ denotes the imaginary part of a complex number and
	$\langle\zeta , {\zeta'}\rangle =\sum\limits_{j = 1}^n \zeta_j
	\overline{\zeta_j'}$ is the Hermitian inner product. In these
	coordinates the neutral element of $\Heis^n$ is $0_{\mathbb
		H^n}=(0_{\mathbb C^n},0)$ and the inverse element of $(\zeta, t)$ is
	$(-\zeta, -t)$. Note that $x = (\xi, \eta, t) = (\zeta, t)$ form a
	real coordinate system for $\mathbb H^n$ and the system of vector
	fields given as differential operators
	$$X_j = \partial_{\xi_j} + 2 \eta_j \partial_{t}, \quad Y_j = \partial_{\eta_j} - 2 \xi_j \partial_{t}, \quad  j \in \{ 1, \ldots n\},\ \ T = \partial_t,$$
	forms a basis for the left invariant vector fields of $\mathbb H^n.$
	The vectors $X_j,Y_j, j \in \{1,...,n\}$ form the basis of the
	horizontal bundle and we denote by $d_{CC}$ the associated Carnot-Carath\'eodory
	metric.

	Following the notations of Ambrosio and Rigot \cite{AR} and Juillet \cite{Juillet-IMNR}, 
	we parametrize the sub-Riemannian geodesics starting from the origin as follows. For  every $(\chi,\theta)\in \mathbb C^n  \times \mathbb R$ we consider the curve $\gamma_{\chi,\theta}:[0,1]\to \mathbb
	H^n$  defined by
	\begin{eqnarray}\label{explicit-geodet}
	\gamma_{\chi,\theta}(s) = \left\{ \begin{array}{lll}
	\left( i \frac{e^{-i\theta s }- 1}{\theta} \chi,  2 |\chi|^2
	\frac{\theta s - \sin(\theta s)}{\theta^2}\right) \  &\mbox{if} &  \theta\neq 0; \\
	(s\chi,0) &\mbox{if} &  \theta=0.
	\end{array}\right.
	\end{eqnarray}
	For the parameters $(\chi,\theta) \in \left( \mathbb C^n \setminus \{ 0_{\mathbb C^n} \}\right)\times
		[-2\pi,2\pi]$, the paths $\gamma_{\chi,\theta}$ are length-minimizing non-constant geodesics in $\mathbb H^n$ joining $0_{\mathbb H^n}$ and
	$\gamma_{\chi,\theta}(1)$. If $\theta \in (-2\pi,2\pi)$ then it follows that the geodesics connecting $0_{\mathbb H^n}$ and
	$\gamma_{\chi,\theta}(1) \neq 0_{\Heis^n}$ are unique, while for $\theta \in \{-2\pi, 2\pi\}$ the uniqueness fails. Let 
	{$$\Gamma_1(\chi,\theta)=\gamma_{\chi,\theta}(1)={\rm the\ endpoint\
		of}\ \gamma_{\chi,\theta},$$} 
and $L=\{(0_{\mathbb C^n},t):t\in \mathbb R\}$ be the center of the group $\mathbb H^n$. {The cut-locus of $0_{\mathbb H^n}$ is $L^*=L\setminus \{0_{\mathbb H^n}\}$.  If $\gamma_{\chi,\theta}(1) \notin L$ then  $\Gamma_1^{-1}(\gamma_{\chi,\theta}(1)) = (\chi, \theta) \in
	{\left( \mathbb C^n \setminus \{ 0_{\mathbb C^n}\} \right)\times (-2\pi,2\pi)}$ is well defined. {Otherwise, $\Gamma_1^{-1}(\gamma_{\chi,\theta}(1)) {\subseteq}
	\mathbb C^n\times \{-2\pi,2\pi\}$ (if $\gamma_{\chi,\theta}(1) \in L^*$) or $\Gamma_1^{-1}(\gamma_{\chi,\theta}(1)) =
		\{0_{\mathbb C^n}\}\times [-2\pi,2\pi]$ (if $\gamma_{\chi,\theta}(1) = 0_{\Heis^n}$).
	
	In analogy to $\tau_ s^{K,N}$ we introduce for  $s \in (0,1)$ the {\it Heisenberg distortion coefficients}
	$\tau_s^n:[0,2\pi]\to [0,\infty]$ defined by
	\begin{eqnarray}\label{concentration}
	\tau_s^n(\theta) = \left\{
	\begin{array}{lll}
	+\infty &\mbox{if} &  \theta=2\pi;\\
	{s^\frac{1}{2n+1}}
	\left(\frac{\sin\frac{\theta s}{2}}{\sin\frac{\theta }{2}}\right)^\frac{2n-1}{2n+1}\left(\frac{\sin\frac{\theta s}{2}-\frac{\theta s}{2}\cos\frac{\theta s}{2}}{\sin\frac{\theta }{2}-\frac{\theta }{2}\cos\frac{\theta
		}{2}}\right)^\frac{1}{2n+1}
	\  &\mbox{if} &  \theta\in( 0,2\pi); \\
	s^\frac{2n+3}{2n+1} &\mbox{if} &  \theta=0.\\
	\end{array}\right.
	\end{eqnarray}
	The function $\theta\mapsto
	\tau_s^n(\theta)$ is increasing on $[0,2\pi]$ (cf. Lemma \ref{lemma-novekvo}), in particular $\tau_s^n(\theta)\to
	+\infty$ as $\theta\to 2\pi$; and also: 
	\begin{equation}\label{tau-becsles}
	\tau_s^n(\theta)\geq \tau_s^n(0)=
	s^\frac{2n+3}{2n+1}\ \ {\rm for\
		every}\ \theta\in [0,2\pi], s\in (0,1).
	\end{equation}
	For $s\in (0,1)$, we introduce the notation 
	\begin{eqnarray}\label{tauk-kozotti-osszefugges}
	\tilde\tau_s^n= s^{-1} \tau_s^n.
	\end{eqnarray}
	If  $x,y\in \mathbb H^n$, $x\neq y$ we let $\theta(x,y) = |\theta|$ with the property that $(\chi, \theta) \in \Gamma_1^{-1}(x^{-1}\cdot y)$. Observe, that $\theta(x,y)$ is well defined and $\theta(x,y)= \theta (y,x)$. If 
	$x=y$ we set $\theta(x,y)= 0$.  
	
	A rough comparison of the Riemannian and Heisenberg distortion coefficients is in order. First of all, 
	both quantities $\tau_s^{K,N}$ and $\tau_s^n$ encode the effect of the curvature in geometric inequalities. Moreover, both of them depend on the dimension of the space, {as indicated by the parameter $N$ in the Riemannian case and $n$ in the Heisenberg case}. However, by $\tau_s^{K,N}$ there is an explicit dependence of the lower bound of the Ricci curvature $K$, while in the expression of $\tau_s^n$ no such dependence shows up.
	
	Let us recall that in case of $\mathbb{R}^n$ the  elegant proof of the Borell-Brascamp-Lieb inequality by the method of optimal mass transportation, see e.g. Villani \cite{Villani1, Villani}  is based on  the concavity of $\det(\cdot)^{\frac{1}{n}}$ defined on the set of $n\times n$-dimensional real symmetric positive semidefinite matrices. In a similar fashion, Cordero-Erausquin, McCann and Schmuckenschl\"ager derive the Borell-Brascamp-Lieb inequality on Riemannian manifolds by the optimal  mass transportation approach from a concavity-type property of $\det(\cdot)^{\frac{1}{n}}$ as well, which holds for the $n \times n$-dimensional matrices, obtained as Jacobians of the map $x \mapsto \exp^{M}_x (-s\nabla^{M} \varphi_{M}(x))$. Here $\varphi_M$ is a $c = \frac{d^2}{2}$-concave map defined on the complete  Riemannian manifold $(M,g)$, $d$ is the Riemannian metric, and $\exp^M$ and  $\nabla^{M}$ denote the exponential map and Riemannian gradient on $(M,g)$. {Here, the concavity is for the Jacobian matrices} $s \mapsto \mathrm{Jac}(\psi^M_s)(x)$, where $\psi^M_s$ is the interpolant map defined for $\mu_0$-a.e. $x \in M$ as 
	$$\psi^M_s(x) = Z_s^M(x,\psi^M(x)).$$
	Here $Z_s^M(A,B)$ is the set of $s$-inter\-mediate points 
	between $A,B\subset M$ w.r.t. to the Riemannian metric $d$, and  
	$\psi^M:M\to M$ is the optimal transport map between the absolutely continuous probability measures $\mu_0$ and $\mu_1$ defined on $M$ minimizing the transportation cost w.r.t. the quadratic cost function $\frac{d^2}{2}$.
	
	Our first result is an appropriate version of the  Jacobian determinant inequality on the Heisenberg group. In order to formulate the precise statement we need to introduce some more notations. 
	
	Let $s\in (0,1)$. Hereafter, $Z_s(A,B)$ denotes the $s$-inter\-mediate set  associated to the nonempty sets $A,B
	\subset \Heis^n$ w.r.t. the Carnot-Carath\'eodory metric $d_{CC}$. Note that  $(\Heis^n,d_{CC})$
	is a geodesic metric space, thus $Z_s(x,y) \neq \emptyset$ for every
	$x,y\in \Heis^n$.
	
	Let $\mu_0$ and $\mu_1$ be two compactly supported probability measures on $\mathbb H^n$ that are  absolutely continuous  w.r.t. $\mathcal L^{2n+1}$. 
	According to Ambrosio and Rigot \cite{AR}, there exists a unique optimal transport map 
	$\psi:\Heis^n\to \Heis^n$ transporting $\mu_0$ to $\mu_1$ associated to the cost function $\frac{d^2_{CC}}{2}$.  If  $\psi_s$ denotes the interpolant optimal transport map associated to $\psi$, defined as
	$$\psi_s(x) = Z_s(x,\psi(x)) \mbox{ for } \mu_0\mbox{-a.e. } x \in \mathbb H^n,$$
	the push-forward measure $\mu_s=(\psi_s)_{\#}\mu_0$ is also absolutely continuous w.r.t. $\mathcal L^{2n+1}$, see Figalli and Juillet \cite{FJ}. 
	Note that the maps $\psi$ and $\psi_s$ are essentially injective thus their inverse functions $\psi^{-1}$ and $\psi_s^{-1}$ are well defined  $\mu_1$-a.e. and $\mu_s$-a.e., respectively, see Figalli and Rifford \cite[Theorem 3.7]{FR} and Figalli and Juillet \cite[p. 136]{FJ}. If $\psi(x)$ is not in the Heisenberg cut-locus of  $x \in \Heis^n$ (i.e., $x^{-1}\cdot \psi(x)\notin L^*,$ which happens $\mu_0$-a.e.) {and  $\psi(x)\neq x$}, there exists a unique 'angle' $\theta_x \in (0, 2\pi)$ defined by $\theta_x = |\theta(x)|$, where $(\chi(x), \theta(x)) \in \left(\mathbb C^n \setminus \{ 0_{\mathbb C^n} \}\right) \times (-2\pi, 2\pi)$ is the unique pair such that $x^{-1}\cdot\psi(x) = \Gamma_1(\chi(x), \theta(x))$. If $\psi(x) = x$, we set $\theta_x = 0$. 
	Observe that the map $x\mapsto \tau_{s}^n(\theta_x)$ is Borel measurable on $\mathbb H^n$. 
	
	Our main result  can now be stated as follows.
	\begin{theorem}\label{TJacobianDetIneq}{\bf (Jacobian determinant inequality on $\mathbb H^n$)} Let $s \in (0,1)$ and assume that $\mu_0$ and $\mu_1$ are two compactly supported, Borel probability measures, both absolutely continuous w.r.t. $\mathcal L^{2n+1}$ on  $\Heis^n$.  Let $\psi:\Heis^n\to \Heis^n$ be the unique optimal transport map transporting $\mu_0$ to $\mu_1$ associated to the cost function $\frac{d^2_{CC}}{2}$ and $\psi_s$ its interpolant map. Then the following Jacobian determinant inequality holds:
		\begin{equation}\label{Jacobi-inequality-elso}
		\left({\rm Jac}(\psi_s)(x)\right)^\frac{1}{2n+1}\geq \tau_{1-s}^n(\theta_x)+\tau_{s}^n(\theta_x)\left({\rm Jac}(\psi)(x)\right)^\frac{1}{2n+1}\mbox{ for } \mu_0 \mbox{-a.e. } x \in \mathbb H^n.
		\end{equation}
	\end{theorem}
	\medskip 
	
	\noindent If $\rho_0$, $\rho_1$ and $\rho_s$ are the density functions of the measures $\mu_0$, $\mu_1$ and  $\mu_s=(\psi_s)_{\#}\mu_0$ w.r.t. to $\mathcal L^{2n+1}$, respectively,  the Monge-Amp\`ere equations 
\begin{equation}\label{MA-bevezeto}
\rho_{0}(x) = \rho_{s}(\psi_{s}(x)) {\rm Jac}(\psi_s)(x), \ \ \rho_{0}(x) = \rho_{1}(\psi(x)) {\rm Jac}(\psi)(x)\mbox{ for } \mu_0 \mbox{-a.e. } x \in \mathbb H^n,
\end{equation}
	show the equivalence of (\ref{Jacobi-inequality-elso}) to 
	\begin{equation}\label{Jacobi-inequality-elso-ekvivalens}
	\rho_s(\psi_s(x))^{-\frac{1}{2n+1}}\geq  
	\tau_{1-s}^n(\theta_x) (\rho_0(x))^{-\frac{1}{2n+1}}   + \tau_{s}^n(\theta_x) (\rho_1(\psi(x)))^{-\frac{1}{2n+1}}\mbox{ for } \mu_0 \mbox{-a.e. } x \in \mathbb H^n.
	\end{equation}

	\noindent It turns out that a version of Theorem \ref{TJacobianDetIneq} holds even in the case when only $\mu_{0}$ is required to be absolutely continuous. In this case we consider only the first term on the right hand side of (\ref{Jacobi-inequality-elso-ekvivalens}). Inequality  (\ref{tau-becsles}) shows that 
	$$\rho_s(y)\leq \frac{1}{(1-s)^{2n+3}}\rho_0(\psi_s^{-1}(y))\mbox{ for } \mu_s \mbox{-a.e. } y \in \mathbb H^n,$$
	which is the main estimate of Figalli and Juillet \cite[Theorem 1.2]{FJ}; for further details see Remark \ref{remark-mikor-nem-feltetlenul-absz-folyt} and Corollary  \ref{Corollary-Figalli-Juillet}.

	\medskip
	

	
	The first application of Theorem \ref{TJacobianDetIneq} is an entropy inequality. In order to formulate the result, we recall that for a function $U:[0, \infty) \to \mathbb R$ one defines the $U$-entropy of an absolutely continuous measure $\mu$ w.r.t. $\mathcal L^{2n+1}$ on $\mathbb H^n$ as
	$${\rm Ent}_U(\mu | \mathcal L^{2n+1}) = \int_{\mathbb H^n} U\left( \rho(x) \right) \dd \mathcal L^{2n+1}(x),$$
	where $\rho=\frac{\dd \mu}{\dd \mathcal L^{2n+1}}$ is the density of $\mu.$
	
	\medskip 
	Our entropy inequality is stated as follows:
	
	\begin{theorem}\label{TEntIneqHeisGen} {\bf (General entropy inequality on $\Heis^n$)}
		Let $s \in (0,1)$ and assume that $\mu_0$ and $\mu_1$ are two compactly supported, Borel probability measures, both absolutely continuous w.r.t. $\mathcal L^{2n+1}$ on  $\Heis^n$ with densities $\rho_{0}$ and $\rho_{1},$ respectively. Let $\psi:\Heis^n\to \Heis^n$ be the unique optimal transport map transporting $\mu_0$ to $\mu_1$ associated to the cost function $\frac{d^2_{CC}}{2}$ and $\psi_s$ its interpolant map.  If $\mu_s=(\psi_s)_{\#}\mu_0$ is
		the interpolant measure between $\mu_0$ and $\mu_1,$ and  $U: [0, \infty) \to \mathbb R$ is a function such that $U(0)=0$ and $t \mapsto t^{2n+1} U\left(\frac{1}{t^{2n+1}}\right)$ is non-increasing and convex, the following entropy inequality holds:
		\begin{eqnarray*}
			{\rm Ent}_{U}(\mu_s | \mathcal L^{2n+1})   &\leq&  (1-s) \int_{\mathbb H^n} \left(\tilde{\tau}_{1-s}^n(\theta_x)\right)^{2n+1} U\left(\frac{\rho_0(x)}{\left(\tilde{\tau}_{1-s}^n(\theta_x)\right)^{2n+1}}\right) \dd \mathcal L^{2n+1}(x) \\
			&&+ s \int_{\mathbb H^n} \left(\tilde{\tau}_{s}^n(\theta_{\psi^{-1}(y)})\right)^{2n+1} U\left(\frac{\rho_1(y)}{\left(\tilde{\tau}_{s}^n(\theta_{\psi^{-1}(y)})\right)^{2n+1}}\right) \dd \mathcal L^{2n+1}(y).
		\end{eqnarray*}
	\end{theorem}
	\medskip	
	\noindent	 Inequality (\ref{tau-becsles}),  Theorem \ref{TEntIneqHeisGen} and the assumptions made for $U$ give the uniform entropy estimate (see also Corollary \ref{Corollary-TEntIneqHeisGen}): 
	\begin{eqnarray*}
		{\rm Ent}_{U}(\mu_s | \mathcal L^{2n+1})   &\leq&  (1-s)^3 \int_{\mathbb H^n}  U\left(\frac{\rho_0(x)}{(1-s)^2}\right) \dd \mathcal L^{2n+1}(x) + s^3 \int_{\mathbb H^n}  U\left(\frac{\rho_1(y)}{s^2}\right) \dd \mathcal L^{2n+1}(y).
	\end{eqnarray*}
	
	\noindent 	
	
\noindent	Various relevant choices of admissible functions $U: [0,\infty) \to \R$ will be presented in the sequel. In particular,  Theorem \ref{TEntIneqHeisGen} provides an 
	{\it curvature-dimension condition} on the  metric measure space $(\mathbb H^n,d_{CC},\mathcal
	L^{2n+1})$ for the choice of  $$U_R(t)=-t^{{1 - \frac{1}{2n+1}}},$$ see  Corollary \ref{TEntIneqHeis-corollary}. Further  consequences of Theorem \ref{TEntIneqHeisGen} are also presented for the Shannon entropy in Corollary \ref{Shannon-TEntIneqHeis-corollary}.
	
	\medskip
	
	Another consequence of Theorem \ref{TJacobianDetIneq} is the following Borell-Brascamp-Lieb inequality:
	
	\begin{theorem}\label{TRescaledBBLWithWeights} {\bf (Weighted Borell-Brascamp-Lieb inequality on $\mathbb H^n$)}
		Fix $s\in (0,1)$ and $p \geq -\frac{1}{2n+1}$. 
		Let $f,g,h:\mathbb H^n\to [0,\infty)$ be Lebesgue integrable
		functions with the property that for all $(x,y)\in \Heis^n\times \Heis^n, z\in Z_s(x,y),$
		\begin{eqnarray}\label{ConditionRescaledBBLWithWeights}
		h(z) \geq M^{p}_s
		\left(\frac{f(x)}{\left(\tilde{\tau}_{1-s}^n(\theta(y,x))\right)^{2n+1}},\frac{g(y)}{\left(\tilde{\tau}_s^n(\theta(x,y))\right)^{2n+1}} \right).
		\end{eqnarray}
		Then the following inequality holds:
		\begin{eqnarray*}
			\int_{\Heis^n} h \geq M^\frac{p}{1+(2n+1)p}_s \left(\int_{\Heis^n}
			f, \int_{\Heis^n} g \right).
		\end{eqnarray*}
	\end{theorem}
	
	
	Consequences of Theorem \ref{TRescaledBBLWithWeights} are uniformly weighted and non-weighted Borell-Brascamp-Lieb inequalities on $\mathbb H^n$ which are stated in Corollaries \ref{CLighterWeightedBBL} and \ref{CRescaledBBLWithoutWeights}, respectively. As particular cases we obtain Pr\'ekopa-Leindler-type inequalities on $\mathbb H^n$,  stated in Corollaries \ref{C-weighted-prekopa-leindler}-\ref{C-nonweighted-prekopa-leindler}.

	Let us emphasize the  difference between the Riemannian and sub-Riemannian versions of the entropy and Borell-Brascamp-Lieb inequalites. In the Riemannian case, we notice the appearance of  
	the distance function in the expression of  $\tau_s^{K,N}(d(x,y))$.  
	The 
	explanation of this phenomenon is that in the Riemannian case the effect of the curvature 
	accumulates in dependence of the distance between $x$ and $y$ in a controlled way, estimated by the lower bound $K$ of the Ricci curvature.	In contrast to this fact, in the sub-Riemannian framework the argument $\theta(x,y)$ appearing in the weight $\tau_s^n(\theta(x,y))$ is {\it not} a distance but a quantity measuring the deviation from the horizontality of the points $x$ and $y$, respectively.  Thus,  in the Heisenberg case the effect of  positive curvature occurs along geodesics between points that are situated in a more vertical position with respect to each other. On the other hand an effect of negative curvature is manifested between points that are in a relative `horizontal position' to each other. The size of the angle $\theta(x,y)$ measures the 'degree of verticality' of the relative positions of $x$ and $y$ which contributes to the curvature. 
	

	\medskip
	
	The geodesic Brunn-Minkowski inequality on the Heisenberg group $\Heis^n$ will be a
	consequence of Theorem \ref{TRescaledBBLWithWeights}. 
	For two nonempty measurable sets $A,B \subset \mathbb H^n$ we introduce the quantity
	$$\Theta_{A,B}=\sup_{A_0,B_0} \inf_{(x,y) \in A_0\times B_0} \left\{|\theta|\in [0,2\pi]:(\chi,\theta)\in \Gamma_1^{-1}(x^{-1}\cdot
	y)\right\},$$
	where the sets $A_0$ and $B_0$ are nonempty, full measure subsets of $A$ and $B$, respectively.

	\begin{theorem}\label{Brunn-Minkowski}  {\bf (Weighted Brunn-Minkowski inequality on $\mathbb H^n$)}
		Let $s\in (0,1)$ and $A$ and $B$ be two nonempty measurable sets of $\mathbb
		H^n$. 
		Then the following geodesic Brunn-Minkowski inequality holds:
		\begin{eqnarray}\label{BM-1}
		\mathcal L^{2n+1}(Z_s(A,B))^\frac{1}{2n+1} \geq
		\tau_{1-s}^n(\Theta_{A,B})\mathcal
		L^{2n+1}(A)^\frac{1}{2n+1}+\tau_{s}^n(\Theta_{A,B})\mathcal
		L^{2n+1}(B)^\frac{1}{2n+1}.
		\end{eqnarray}
		
	\end{theorem}
	Here we consider the outer Lebesgue measure whenever $Z_s(A,B)$ is not measurable, and the convention $+\infty\cdot 0=0$ for the right hand side of {\rm(\ref{BM-1})}. The latter case may happen e.g. when $A^{-1}\cdot B\subset L=\{0_{\mathbb C^n}\}\times \mathbb R;$ indeed, in this case $\Theta_{A,B}=2\pi$ and $\mathcal
	L^{2n+1}(A)=\mathcal
	L^{2n+1}(B)=0.$
	
	The value $\Theta_{A,B}$ represents a typical Heisenberg
	quantity indicating a lower bound of the deviation of  an {\it essentially horizontal} position of the sets $A$ and $B$. An intuitive description of the role of weights $\tau_{1-s}^n(\Theta_{A,B})$ and $\tau_s^n(\Theta_{A,B})$ 
	in 
	(\ref{BM-1}) will be given in Section \ref{section-BM}.
	
	By  Theorem \ref{Brunn-Minkowski} we deduce several forms of the Brunn-Minkowski inequality, see Corollary \ref{CBrunn-Minkowski-2}. Moreover, the weighted Brunn-Minkowski inequality implies  the measure contraction property {\rm{\textsf{ MCP}}}$(0,2n+3)$ 
	on $\mathbb H^n$ proved
	by Juillet \cite[Theorem 2.3]{Juillet-IMNR}, see also  Corollary \ref{MCP-1}, namely, 
	for every  $s\in (0,1)$, $x\in \mathbb H^n$
	and nonempty measurable set $E\subset \mathbb H^n$,
	$$\displaystyle \mathcal L^{2n+1}(Z_s(x,E)) \geq s^{2n+3}\mathcal
	L^{2n+1}(E).$$


	Our proofs are based on techniques of optimal mass transportation and Riemannian approximation of the sub-Riemannian structure. We use extensively the machinery developed by Cordero-Erausquin, McCann and Schmuckenschl\"ager \cite{McCann} on Riemannian manifolds and the results of Ambrosio and Rigot \cite{AR} and Juillet \cite{Juillet-IMNR} on $\mathbb H^n$. 
	In our approach we can avoid the blow-up of the Ricci curvature to $-\infty$  by {\it not considering} limits of the expressions of $\tau_s^{K,N}$. Instead of this, 
	we apply the limiting procedure to the coefficients expressed in terms of volume distortions. 
	It turns out that one can directly calculate these 
	volume distortion coefficients in terms of Jacobians of exponential maps in the Riemannian 
	approximation. These quantities behave in a much better way under the limit, avoiding blow-up phenomena. The calculations are based on an explicit parametrization of the Heisenberg group and the approximating Riemannian manifolds by an appropriate set of spherical coordinates that are based on a fibration of the space by geodesics.

	The paper is organized as follows. In the second section we present a series of preparatory lemmata obtaining the Jacobian representations of the volume distortion coefficients in the Riemannian approximation of the Heisenberg group and we discuss their limiting behaviour. In the third section we present the proof of our main results, i.e., the Jacobian determinant inequality, various entropy inequalities and   Borell-Brascamp-Lieb inequalities.  
	The forth section is devoted to geometric aspects of the Brunn-Minkowski inequality. In the last section we indicate further perspectives related to this research. The results of this paper have been announced in \cite{BKS-CR}.

	\medskip
	
	\noindent 	{\bf Acknowledgements.} The authors wish to express their gratitude to Luigi Ambrosio,  Nicolas Juillet, Pierre Pansu, Ludovic Rifford, S\'everine Rigot and Jeremy Tyson for helpful conversations on various topics related to this paper.

	\section{Preliminary  results}\label{Section-2}
	

	\subsection{Volume distortion coefficients in $\mathbb H^n$}
	The left translation $l_x: \Heis^n \to \Heis^n$ by the element $x
	\in \Heis^n$ is given by
	$l_x(y) = x \cdot y$ for all $y \in \Heis^n.$
	One can observe that $l_x$ is affine, associated to a matrix with
	determinant 1. Therefore the Lebesgue measure of $\R^{2n+1}$ will be
	the Haar measure on $\Heis^n$ (uniquely defined up to a positive
	multiplicative constant).
	
	For $\lambda > 0$ define the nonisotropic dilation $\rho_{\lambda}: \Heis^n \to
	\Heis^n$ as
	$\rho_{\lambda} \left(\zeta,t\right) = (\lambda \zeta, \lambda^2 t), \quad \forall (\zeta, t) \in \Heis^n.$
	Observe that for any  measurable set $A\subset \mathbb H^n$, $$\Lmeas^{2n+1}(\rho_{\lambda}(A)) = \lambda^{2n+2}
	\Lmeas^{2n+1}(A),$$  thus the homogeneity
	dimension of the Lebesgue measure $\Lmeas^{2n+1}$ is $2n+2$ on $\mathbb H^n$.
	
	In order to equip the Heisenberg group with the
	Carnot-Carath\'eodory metric we consider the basis of the space of the horizontal
	left invariant vector fields $\{X_1, \dots, X_n,$ $ Y_1,\ldots Y_n\}$. A horizontal curve is an absolutely continuous curve $\gamma:
	[0, r] \to \Heis^n$ for which there exist measurable functions $h_j:
	[0,r] \to \R$ ($j = 1, \ldots, 2n$) such that
	$$\dot{\gamma}(s) = \sum\limits_{j=1}^n \left[h_j(s) X_j(\gamma(s)) + h_{n+j}(s) Y_j(\gamma(s))\right] \quad \mbox{a.e. } s \in [0,r].$$
	The length of this curve is
	$$l(\gamma) = \int\limits_{0}^{r} ||\dot{\gamma}(s)|| \dd s= \int\limits_{0}^{r} \sqrt{\sum\limits_{j = 1}^{n} [h^2_j(s) + h_{n+j} ^2(s)]} \dd s = \int\limits_{0}^{r} \sqrt{\sum\limits_{j = 1}^{n} [\dot{\gamma}^2_j(s) + \dot{\gamma}_{n+j} ^2(s)]} \dd s.$$
	The classical Chow-Rashewsky theorem assures that any two points
	from the Heisenberg group can be joined by a horizontal curve, thus
	it makes sense to define the distance of two points as the infimum
	of lengths of all horizontal curves connecting the points, i.e.,
	$$d_{CC}(x,y) = \inf \{l(\gamma): \gamma \mbox{ is a horizontal curve joining } x \mbox{ and } y\};$$
	$d_{CC}$  is called the Carnot-Carath\'eodory metric.
	The left invariance and homogeneity of the vector fields $X_1,
	\dots, X_n, Y_1,\ldots Y_n$ are inherited by the distance $d_{CC}$,
	thus
	$$d_{CC}(x, y) = d_{CC}(0_{\mathbb H^n}, x^{-1} \cdot y) \quad {\rm for\ every}\ x,y \in \Heis^n,$$
	and
	$$d_{CC}(\rho_\lambda(x), \rho_\lambda(y)) = \lambda d_{CC}(x, y)\quad {\rm for\ every}\  x,y \in \Heis^n \mbox{ and } \lambda > 0.$$
	
	We recall the curve $\gamma_{\chi,\theta}$ introduced in
	(\ref{explicit-geodet}). One can observe that for every $x\in
	\mathbb H^n\setminus L$, there exists a unique minimal geodesic
	$\gamma_{\chi,\theta}$ joining $0_{\mathbb H^n}$ and $x$, where
	$L=\{0_{\mathbb C^n}\}\times \mathbb R$ is the center of $\mathbb H^n.$
	 In the sequel, following Juillet
	\cite{Juillet-IMNR}, we consider the diffeomorphism $\Gamma_s : \left(\mathbb C^n
	\setminus \{0_{\mathbb C^n}\}\right) \times (-2\pi, 2\pi) \to \Heis^n \setminus L
	$ defined by
	\begin{equation}\label{Gamma-s-definicio}
	\Gamma_s(\chi,\theta)=\gamma_{\chi,\theta}(s).
	\end{equation}
	By
	\cite[Corollary 1.3]{Juillet-IMNR}, the Jacobian of $\Gamma_s$ for
	$s\in (0,1]$ and  $(\chi,\theta)\in (\mathbb C^n\setminus \{0_{\mathbb
		C^n}\})\times (-2\pi,2\pi)$ is
	\begin{eqnarray}\label{Jacobian-Juillet}
	{\rm Jac}(\Gamma_s)(\chi,\theta) = \left\{
	\begin{array}{lll}
	2^{2n+2} s |\chi|^2
	\left(\frac{\sin\frac{\theta s}{2}}{\theta}\right)^{2n-1}\frac{\sin\frac{\theta s}{2}-\frac{\theta s}{2}\cos\frac{\theta s}{2}}{\theta^3}
	\  &\mbox{if} &  \theta\neq 0; \\
	\frac{s^{2n+3}|\chi|^2}{3} &\mbox{if} &  \theta=0.
	\end{array}\right.
	\end{eqnarray}
	{In particular,  ${\rm Jac}(\Gamma_s)(\chi,\theta)\neq 0$ for every
		$s\in (0,1]$ and  $(\chi,\theta)\in (\mathbb C^n\setminus \{0_{\mathbb
			C^n}\})\times (-2\pi,2\pi)$. Moreover, by (\ref{concentration}) and (\ref{Jacobian-Juillet}), we have for every $\theta\in [0,2\pi)$ (and $\chi\neq 0_{\mathbb C^n}$) that 
		$$\tau_s^n(\theta)=\left(\frac{{\rm Jac}(\Gamma_s)(\chi,\theta)}{{\rm Jac}(\Gamma_1)(\chi,\theta)}\right)^\frac{1}{2n+1}.$$}
	\begin{lemma}\label{lemma-novekvo}
		{Let $s\in (0,1).$ The function $\tau_s^n$ is increasing on $[0,2\pi]$.}
	\end{lemma}
	
	{ {\it Proof.} 
		Let $s\in (0,1)$ be fixed and consider the functions $f_{i,s}:(0,\pi)\to \mathbb R$, $i\in \{1,2\},$ given by $$f_{1,s}(t)=\frac{\sin{(t s)}}{\sin t}\ \ {\rm and}\ \ f_{2,s}(t)=\frac{\sin(ts)-ts\cos(ts)}{\sin t-t\cos t }.$$ 
		Note that both functions $f_{i,s}$ are positive on $(0,\pi)$, $i\in \{1,2\}$. 
		First, for every $ t\in (0,\pi)$ one has 
		$$\frac{f_{1,s}'(t)}{f_{1,s}(t)}=s\cot{(t s)}-\cot t=2\pi^2t(1-s^2)\sum_{k=1}^\infty\frac{k^2}{(\pi^2k^2-t^2)(\pi^2k^2-(ts)^2)}> 0,$$
		where we use the Mittag-Leffler expansion of the cotangent function  
		$$\cot t=\frac{1}{t}+2t\sum_{k=1}^\infty\frac{1}{t^2-\pi^2k^2}.$$ Therefore, $f_{1,s}$ is increasing on  $(0,\pi)$.  In a similar way, we have that 
		\begin{eqnarray*}
			f_{2,s}'(t)&=&t\frac{\sin t \sin(ts)}{(\sin t-t\cos t )^2}(s^2(1-t\cot t)-(1-ts\cot(ts)))\\
			&=&2t^5s^2(1-s^2)\frac{\sin t \sin(ts)}{(\sin t-t\cos t )^2}\sum_{k=1}^\infty\frac{1}{(\pi^2k^2-t^2)(\pi^2k^2-(ts)^2)}> 0.
		\end{eqnarray*}
		Thus, $f_{2,s}$ is also increasing on  $(0,\pi)$. Since  
		$$\tau_s^n(\theta)=s^\frac{1}{2n+1}f_{1,s}\left({\theta}/{2}\right)^\frac{2n-1}{2n+1}f_{2,s}\left({\theta}/{2}\right)^\frac{1}{2n+1},$$ the claim follows. 
		\hfill $\square$}
	\medskip
	
	Let $s\in (0,1)$ and $x,y\in \mathbb H^n$ be such that $x\neq y$.  If $B(y,r)=\{w\in \mathbb H^n:d_{CC}(y,w)<r\}$ is the open CC-ball with center $y\in \mathbb H^n$ and radius $r>0$, we introduce the
	{\it Heisenberg volume distortion coefficient}
	$$
	v_s(x,y) =  \limsup\limits_{r \to 0}\frac{\Lmeas^{2n+1} \left(
		Z_s(x, B(y, r))\right)}{ \Lmeas^{2n+1} \left( B(y, sr)\right)}.$$ 

	The following property gives a formula for the Heisenberg volume distortion
	coefficient in terms of the Jacobian $\Jac (\Gamma_s)$.
	\begin{proposition}\label{proposition-repres-jacobi} Let $s \in (0,1)$ be fixed.
		For  $x, y \in \Heis^n$ such that $x^{-1}\cdot y \notin L$  let
		$(\chi,\theta)=\Gamma_1^{-1}(x^{-1}\cdot y)$. Then
		\begin{itemize}
			\item[{\rm (i)}] $\displaystyle v_s(x,y) = \frac{1}{s^{2n+2}} \frac{\Jac
				(\Gamma_s)(\chi,\theta)}{\Jac (\Gamma_1)(\chi,\theta)};$
			\item[{\rm (ii)}] $\displaystyle v_{1-s}(y,x) = \frac{1}{(1-s)^{2n+2}} \frac{\Jac
				(\Gamma_{1-s})(\chi,\theta)}{\Jac (\Gamma_1)(\chi,\theta)}. $
		\end{itemize}
	\end{proposition}
	{\it Proof.} (i) By left translation, we have that $Z_s(x, B(y,
	r))=x\cdot  Z_s(0_{\mathbb H^n}, B(x^{-1}\cdot y, r))$. Thus, on one
	hand, we have
	\begin{eqnarray*}
		v_s(x,y)  &=& \lim\limits_{r \to 0}\frac{\Lmeas^{2n+1} \left( x\cdot  Z_s(0_{\mathbb H^n},
			B(x^{-1}\cdot y, r))\right)}{ \Lmeas^{2n+1} \left(x\cdot
			B(x^{-1}\cdot y, sr)\right)} = \lim\limits_{r \to
			0}\frac{\Lmeas^{2n+1} \left( Z_s(0_{\mathbb H^n}, B(x^{-1}\cdot y,
			r))\right)}{ \Lmeas^{2n+1} \left( B(x^{-1}\cdot y, sr)\right)} \\
		&=& \lim\limits_{r \to 0}\frac{\Lmeas^{2n+1} \left(  B(x^{-1}\cdot y, r)\right)}{
			\Lmeas^{2n+1} \left( B(x^{-1}\cdot y, sr)\right)}\frac{\Lmeas^{2n+1} \left( Z_s(0_{\mathbb H^n},
			B(x^{-1}\cdot y, r))\right)}{ \Lmeas^{2n+1} \left( B(x^{-1}\cdot y,
			r)\right)}.
	\end{eqnarray*}
	Because of the homogeneities of $d_{CC}$ and $\mathcal L^{2n+1}$, we
	have
	\begin{eqnarray*}
		v_s(x,y) = \frac{1}{s^{2n+2}}\lim\limits_{r \to
			0}\frac{\Lmeas^{2n+1} \left( Z_s(0_{\mathbb H^n}, B(x^{-1}\cdot y,
			r))\right)}{ \Lmeas^{2n+1} \left( B(x^{-1}\cdot y, r)\right)}.
	\end{eqnarray*}
	Since $x^{-1}\cdot y \notin L$,  we have that $ B(x^{-1}\cdot y,r)
	\cap L = \emptyset$ for  $r$ small enough, thus the map $\Gamma_s
	\circ \Gamma^{-1}_1$ realizes a diffeomorphism between the sets
	$B(x^{-1}\cdot y,r)$ and $Z_s(0_{\mathbb H^n},B(x^{-1}\cdot y,r))$.
	This constitutes the basis for the following change of variable
	\begin{eqnarray*}
		\Lmeas^{2n+1} \left( Z_s(0_{\mathbb H^n}, B(x^{-1} \cdot y,
		r))\right) = \int\limits_{Z_s(0_{\mathbb H^n}, B(x^{-1} \cdot y,
			r))}  \dd\Lmeas^{2n+1}  = \int\limits_{B(x^{-1} \cdot y,r)} \frac{\Jac
			(\Gamma_s)(\Gamma_1^{-1}(w))}{\Jac (\Gamma_1)(\Gamma_1^{-1}(w))} \dd
		\Lmeas^{2n+1}(w).
	\end{eqnarray*}
	By the continuity of the integrand in the latter expression, the
	volume derivative of $ Z_s(0_{\mathbb H^n}, \cdot)$ at the point
	$x^{-1} \cdot y$ is
	\begin{eqnarray*}
		\frac{\Jac (\Gamma_s)(\Gamma_1^{-1}(x^{-1} \cdot y))}{\Jac
			(\Gamma_1)(\Gamma_1^{-1}(x^{-1} \cdot y))},
	\end{eqnarray*}
	which gives precisely the claim.

	
	
	(ii) At first glance, this property  seems to be just the symmetric
	version of (i). Note however that
	$$v_{1-s}(y,x)=v_{1-s}(0_{\mathbb H^n},y^{-1}\cdot
	x)=v_{1-s}(0_{\mathbb H^n},-x^{-1}\cdot y),$$ thus we need the
	explicit form of the geodesic from $0_{\mathbb H^n}$ to
	$-x^{-1}\cdot y$ in terms of $(\chi,\theta)$. A direct computation
	based on (\ref{Gamma-s-definicio}) shows that
	$$\Gamma_1\left(-\chi e^{-i \theta},-\theta\right)=-\Gamma_1(\chi,\theta)=-x^{-1}\cdot
	y.$$ Therefore, the minimal geodesic joining $0_{\mathbb H^n}$ and
	$-x^{-1}\cdot y$ is given by the curve  $s\mapsto
	\Gamma_s\left(-\chi e^{-i \theta},-\theta\right),$ $s\in [0,1].$
	Now, it remains to apply (i) with the corresponding modifications,
	obtaining
	\begin{eqnarray*}
		v_{1-s}(y,x) &=& \frac{1}{(1-s)^{2n+2}}\frac{{\rm Jac}(\Gamma_{1-s})\left(-\chi e^{-i \theta},-\theta\right)}{{\rm
				Jac}(\Gamma_1)\left(-\chi e^{-i \theta},-\theta\right)} =
		\frac{1}{(1-s)^{2n+2}} \frac{\Jac (\Gamma_{1-s})(\chi,\theta)}{\Jac
			(\Gamma_1)(\chi,\theta)},
	\end{eqnarray*}
	which concludes the
	proof. \hfill $\square$\\
	
	{For further use (see Proposition \ref{bridge-proposition}), we consider 
		\begin{eqnarray}\label{Heisenberg-distortion}
		v_s^0(x,y)
		=
		\left\{
		\begin{array}{lll}
		sv_s(x,y)
		\  &\mbox{if} &  x\neq y; \\
		s^2 &\mbox{if} &  x=y.
		\end{array}\right.
		\end{eqnarray}}
	{
		\begin{corollary}\label{corollary-1} Let $s\in (0,1)$ and  $x, y \in \Heis^n$  such that $x\neq y.$ The following properties hold:
			\begin{enumerate}
				\item[{\rm (i)}] if $x^{-1}\cdot y \notin
				L$, then $(\chi,\theta)=\Gamma_1^{-1}(x^{-1}\cdot y)\in (\mathbb C^n\setminus
				\{0_{\mathbb C^n}\})\times (-2\pi,2\pi)$ and \begin{eqnarray*}
					v_s(x,y)=  \left\{
					\begin{array}{lll}
						s^{-(2n+1)}
						\left(\frac{\sin\frac{\theta s}{2}}{\sin\frac{\theta }{2}}\right)^{2n-1}\frac{\sin\frac{\theta s}{2}-\frac{\theta s}{2}\cos\frac{\theta s}{2}}{\sin\frac{\theta }{2}-\frac{\theta }{2}\cos\frac{\theta }{2}}
						\  &\mbox{if} &  \theta\neq 0; \\
						s &\mbox{if} &  \theta=0;
					\end{array}\right.
				\end{eqnarray*}
				\item[{\rm (ii)}]  if $x^{-1}\cdot y \in
				L$, then $v_s(x,y)=+\infty.$
			\end{enumerate}
			Moreover, we have for every $s\in (0,1)$ and $x,y\in \mathbb H^n$ that
			\begin{equation}\label{kapcsolat-volume-es-CD}
			\left(\tilde\tau_s^n(\theta(x,y))\right)^{2n+1}=v_s^0(x,y)\geq s^2.
			\end{equation}
			Similar relations  hold
			for $v_{1-s}(y,x)$ by replacing  $s$ by $(1-s)$. 
		\end{corollary}
	}

	{\it Proof.} (i) Directly follows by Proposition  \ref{proposition-repres-jacobi} and relation (\ref{Jacobian-Juillet}).

	(ii) Let $t\in \mathbb R\setminus \{0\}$ be such that $x^{-1}\cdot y=(0_{\mathbb C^n},t)\in L^*=L\setminus \{0_{\mathbb H^n}\}$; for simplicity, we assume that $t>0.$ Let us choose $r < \frac{\sqrt {t}}{2}$. Then for every $w\in B(x^{-1}\cdot y,r)\setminus L$, there exists a unique $(\chi_w,\theta_w)\in \left(\mathbb C^n
	\setminus \{0_{\mathbb C^n}\}\right) \times (0, 2\pi)$ such that $(\chi_w,\theta_w)=\Gamma_1^{-1}(w)$. Moreover, by (\ref{explicit-geodet}), it follows that  
	\begin{equation}\label{kesobb-kell-peldaba}
	0<\sin\left(\frac{\theta_w}{2}\right)\leq c_1r\ \ {\rm for\ every} \ w\in  B(x^{-1}\cdot y,r)\setminus L,
	\end{equation}
	where $c_1>0$ is a constant which depends on $t>0$ (but not on $r>0$). To check inequality (\ref{kesobb-kell-peldaba}) we may replace the ball $B(x^{-1} \cdot y, r)$ in the Carnot-Carath\'eodory metric $d_{CC}$ by the ball in the Kor\'anyi metric $d_K$ (introduced as the gauge metric in \cite{AR}). Since the two metrics are bi-Lipschitz equivalent, it is enough to check (\ref{kesobb-kell-peldaba}) for the Kor\'anyi ball; for simplicity, we keep the same notation. 
	
	Let $w=\Gamma_1(\chi_w, \theta_w) \in B(x^{-1} \cdot y, r)\setminus L$.   Since $r< \frac{\sqrt{t}}{2} $, it is clear that
	$B(x^{-1} \cdot y, r) \cap \left(\mathbb C^{n} \times \{0\}\right) = \emptyset$; therefore,  $\theta_w \neq 0$. 
	Using the notation $(\zeta_w,t_w) \in \mathbb{C}^{n} \times \mathbb{R}$ for the point $w \in \mathbb H^n$, 
	due to the properties of the Kor\'anyi metric, from $d_K(\Gamma_1(\chi_w,\theta_w), (0_{\mathbb C^n},t)) \leq r$  it follows that $|\zeta_w| \leq r$ and $\sqrt{|t_w-t|} \leq r$. By (\ref{explicit-geodet}) we obtain the estimates
	\begin{eqnarray}\label{sin-also-korlathoz}
	|\zeta_w| = 2 \sin\left(\frac{\theta_w}{2}\right) \frac{|\chi_w|}{|\theta_w|} \leq r, 
	\end{eqnarray}
	\begin{eqnarray}\label{chi-also-korlathoz}
	\sqrt{|t_w-t|} = \sqrt{\left|\frac{2|\chi_w|^2}{\theta_w^2}(\theta_w-\sin(\theta_w))-t\right|} \leq r. 
	\end{eqnarray}
	Recalling that $r < \frac{\sqrt{t}}{2}$, by inequality (\ref{chi-also-korlathoz}) we obtain that
	\begin{eqnarray}\label{chi-also-korlathoz2}
	|\chi_w|^2\frac{\theta_w-\sin(\theta_w)}{\theta_w^2} \geq \frac{3t}{8}.
	\end{eqnarray}
	Since $\frac{\theta-\sin(\theta)}{\theta^2}\in \left(0,\frac{1}{\pi}\right]$ for every $\theta\in (0,2\pi)$, 
	by inequality (\ref{chi-also-korlathoz2}) we conclude that 
	$\frac{1}{|\chi_w|} \leq \sqrt{\frac{8}{3t\pi}}.$
	Combining this estimate with inequality (\ref{sin-also-korlathoz}), it yields that
	$\sin\left(\frac{\theta_w}{2}\right) \leq r \frac{|\theta_w|}{2|\chi_w|} \leq r \frac{\pi}{|\chi_w|} \leq r\sqrt{\frac{8\pi}{3t}} ,$
	proving inequality (\ref{kesobb-kell-peldaba}).
	
Note that $\theta_w$ is close to $2\pi$ whenever $r$ is very small. Therefore, by continuity reasons, since $r< \frac{\sqrt{t}}{2} $, one has that $\sin\frac{\theta_w s}{2}-\frac{\theta_w s}{2}\cos\frac{\theta_w s}{2}\geq c_2^1$,  $\sin\frac{\theta_w }{2}-\frac{\theta_w }{2}\cos\frac{\theta_w }{2}\leq c_2^2$ and $\sin\frac{\theta_w s}{2}\geq c_2^3$ for every $w\in B(x^{-1} \cdot y, r)\setminus L$, where the numbers $c_2^1,c_2^2,c_2^3>0$ depend only on $s\in (0,1)$, $t>0$ and $n\in \mathbb N$. Consequently, by relation  (\ref{Jacobian-Juillet}) one has for every $w\in B(x^{-1} \cdot y, r)\setminus L$ that 
$$\frac{\Jac
	(\Gamma_s)(\chi_w,\theta_w)}{\Jac (\Gamma_1)(\chi_w,\theta_w)}\geq \frac{c_2}{\left(\sin\left(\frac{\theta_w}{2}\right)\right)^{2n-1}},$$
where $c_2>0$ depends on $c_2^i>0$, $i\in \{1,2,3\}. $
	
	Since the map $\Gamma_s
	\circ \Gamma^{-1}_1$ is a diffeomorphism between 
	the sets $B(x^{-1}\cdot y,r)\setminus L$ and $Z_s(0_{\mathbb H^n},B(x^{-1}\cdot y,r)\setminus L)$, a similar argument as in the proof of Proposition \ref{proposition-repres-jacobi}  gives 
	\begin{eqnarray*}
		\Lmeas^{2n+1} \left( Z_s(0_{\mathbb H^n}, B(x^{-1} \cdot y,
		r))\right) &\geq &\Lmeas^{2n+1} \left( Z_s(0_{\mathbb H^n}, B(x^{-1} \cdot y,
		r)\setminus L)\right) \\ &= & \int\limits_{Z_s(0_{\mathbb H^n}, B(x^{-1} \cdot y,
			r)\setminus L)}  \dd\Lmeas^{2n+1} \\& =& \int\limits_{B(x^{-1} \cdot y,r)\setminus L} \frac{\Jac
			(\Gamma_s)(\chi_w,\theta_w)}{\Jac (\Gamma_1)(\chi_w,\theta_w)} \dd
		\Lmeas^{2n+1}(w)\\&\geq &c_2\int\limits_{B(x^{-1} \cdot y,r)\setminus L} \frac{1}{\left(\sin\left(\frac{\theta_w}{2}\right)\right)^{2n-1}}\dd
		\Lmeas^{2n+1}(w).
	\end{eqnarray*}
	By the latter estimate and  (\ref{kesobb-kell-peldaba}) we have
	\begin{equation}\label{kesobb-jo-lesz}
	\Lmeas^{2n+1} \left( Z_s(0_{\mathbb H^n}, B(x^{-1} \cdot y,
	r))\right) \geq\frac{c_2}{c_1^{2n-1}}\frac{\Lmeas^{2n+1}(B(x^{-1} \cdot y,r)\setminus L)}{r^{2n-1}}.
	\end{equation}
	Consequently, since $\mathcal L^{2n+1}(L)=0$, 
	we have 
	\begin{eqnarray*}
		v_s(x,y) &=& v_s(0_{\mathbb H^n},x^{-1}\cdot y) =\limsup\limits_{r \to 0}\frac{\Lmeas^{2n+1} \left(
			Z_s(0_{\mathbb H^n}, B(x^{-1}\cdot y, r))\right)}{ \Lmeas^{2n+1} \left( B(x^{-1}\cdot y, sr)\right)}\\ &\geq& \frac{c_2}{c_1^{2n-1}s^{2n+2}}\limsup\limits_{r \to 0}\frac{1}{r^{2n-1}}\\&=&+\infty.
	\end{eqnarray*}

	
	
	{The first part of relation (\ref{kapcsolat-volume-es-CD}) follows by (\ref{concentration}), (\ref{Heisenberg-distortion}) and (i)\&(ii), while the inequality $v_s^0(x,y)\geq s^2$ is a consequence of Lemma \ref{lemma-novekvo}. }
	\hfill $\square$
	
	\begin{remark}\rm 
		The fact that $v_s(x,y)=+\infty$ for $x^{-1}\cdot y\in L^*$ encompasses another typical sub-Riemannian feature of the Heisenberg group $\mathbb H^n$ showing that on  'vertical directions' the curvature blows up even in small scales (i.e., when $x$ and $y$ are arbitrary close to each other), described by the behavior of the Heisenberg volume distortion coefficient. This phenomenon shows another aspect of the singular space structure of the Heisenberg group $\mathbb H^n$. 
	\end{remark}
	
	\medskip
	\subsection{Volume distortion coefficients in the Riemannian approximation $M^\eps$ of $\mathbb H^n$}\label{section-aprox-Riemann} We introduce specific Riemannian manifolds in order to approximate the Heisenberg group $\mathbb H^n$, following Ambrosio and Rigot \cite{AR} and Juillet \cite{Juillet-IMNR, Juillet-calculus}.
	
	For every $\eps > 0$,  let $M^\eps=\R^{2n+1}$ be equipped with the
	usual Euclidean topology and with the Riemannian structure where the
	orthonormal basis of the metric tensor $g^\eps$ at the point
	$x=(\xi, \eta, t)$ is given by the vectors (written as differential operators): 
	$$X_j = \partial_{\xi_j} + 2 \eta_j \partial_{t}, \quad Y_j = \partial_{\eta_j} - 2 \xi_j \partial_{t} \quad {\rm for\ every}\ j = 1, \ldots n,$$
	and
	$$T^\eps = \eps \partial_t = \eps T.$$
	On $(M^\eps,g^\eps)$ we consider the  measure $\textsf{m}_\eps$ with the canonical volume element 
	$${\rm d}\textsf{m}_{\eps}=\sqrt{{\rm det}g^\eps}{\rm d}\mathcal
	L^{2n+1}=\frac{1}{\eps}{\rm d}\mathcal L^{2n+1},$$ and ${\rm
		Vol}^\eps(A)=\displaystyle\int_A {\rm d}\textsf{m}_{\eps}$ for every measurable set
	$A\subset M^\eps.$
	The length of a piecewise $C^1$ curve $\gamma: [0,1] \to M^\eps$ is
	defined as
	$$l_{\eps}(\gamma) = \int\limits_{0}^{1} ||\dot{\gamma}(s)||_\eps \dd s = \int\limits_{0}^{1} \sqrt{\sum\limits_{j=1}^{n} [\dot{\gamma}_j(s)^2 + \dot{\gamma}_{n+j}(s)^2] + \dot{\gamma}_\eps(s)^2} \dd s,$$
	where $(\gamma_1(s), \ldots, \gamma_{2n+1}(s))$ are the cartesian
	coordinates of $\gamma(s)$ expressed in the canonical basis of
	$M^\eps = \R^{2n+1}$ and $(\dot{\gamma}_1(s), \ldots,
	\dot{\gamma}_{2n}(s), \dot{\gamma}_\eps(s))$ are the coordinates of
	$\dot{\gamma}(s) \in T_{\gamma(s)}M^\eps$ in the basis
	$X_1(\gamma(s)), \ldots, X_n(\gamma(s)), Y_1(\gamma(s)), \ldots,
	Y_n(\gamma(s)), T^\eps(\gamma(s))$. One can check that
	$\dot{\gamma}_{j}(s)$ is equal with the $j$-th cartesian coordinate
	of $\dot{\gamma}(s)$, $j=1,...,2n,$ and
	$$\dot{\gamma}_\eps(s) = \frac{1}{\eps} \left( \dot{\gamma}_{2n+1}(s) - 2 \sum\limits_{j=1}^n \left[\gamma_{n+j}(s) \dot{\gamma}_j(s) - \gamma_j(s) \dot{\gamma}_{n+j}(s)\right] \right).$$
	
	
	The induced Riemannian distance on $(M^\eps,g^\eps)$ is
	$$d^\eps(x,y) = \inf \{l_\eps(\gamma): \gamma \mbox{ is a piecewise}\ C^1\ \mbox{curve in } M^\eps \mbox{ joining } x \mbox{ and } y\}.$$
	Note that $(M^\eps,d^\eps)$ is complete and the distance $d^\eps$
	inherits the left invariance of the vector
	fields $X_1, \ldots, X_n, Y_1, \ldots, Y_n, T^\eps,$ similarly as in
	the Heisenberg group $\mathbb H^n$. Moreover, one can observe
	that $d^\eps$ is {decreasing} w.r.t. $\eps>0$ and due to Juillet \cite{Juillet-calculus} for a fixed $c>0$ constant,
	\begin{equation}\label{dist-eps-Juillet}
	{d_{CC}(x,y) - c\pi \eps \leq d^{\eps}(x,y)\leq d_{CC}(x,y) \quad {\rm for\ every}\ x,y \in \Heis^n,  \eps>0 ;}
	\end{equation}
	thus 
	$d_{CC}(x,y) = \sup\limits_{\eps > 0} d^{\eps}(x,y) = \lim\limits_{\eps \searrow 0} d^{\eps}(x,y)$ for every  $x,y \in \Heis^n.$
	
	For $\varepsilon>0$ fixed, we recall from Ambrosio and  Rigot
	\cite{AR} that the $\eps$-{\it geodesic} $\gamma^\eps:[0,1]\to
	M^\eps$ with starting point $0_{\mathbb H^n}$ and initial vector
	$$w^\eps=\sum_{j=1}^nw_j^\eps X_j(0_{\mathbb H^n})+\sum_{j=1}^nw_{j+n}^\eps
	Y_j(0_{\mathbb H^n})+w_{2n+1}^\eps T^\eps(0_{\mathbb H^n})\in
	T_{0_{\mathbb H^n}} M^\eps$$ is
	\begin{equation}\label{explicit-geodet-epsz}
	\gamma^\eps(s)=\exp_{0_{\mathbb H^n}}^\eps(sw^\eps).
	\end{equation}
	Using the complex notation $\mathbb C^n \times \mathbb R$ for the  Heisenberg
	group $\Heis^n$, we can write the expression of the $\eps$-geodesics
	$\gamma^\eps$ explicitly as
	\begin{eqnarray}\label{eps-geodetikus}
	\gamma^\eps(s) = \left\{
	\begin{array}{lll}
	\left( i \frac{e^{-i \theta^\eps s} -
		1}{\theta^\eps} \chi^\eps, \frac{\eps^2}{4}(\theta^\eps s) \right.
	\left. +2 |\chi^\eps|^2 \frac{\theta^\eps s - \sin \theta^\eps s
	}{(\theta^\eps)^2}  \right) 
	\  &\mbox{if} &  \theta^\eps\neq 0; \\
	(s\chi^\eps,0) &\mbox{if} &  \theta^\eps=0, 
	\end{array}\right.
	\end{eqnarray}
	where
	\begin{equation}\label{jeloles-chi}
	\theta^\eps = \frac{4
		w_{2n+1}^\eps}{\eps}\ \ {\rm and} \ \ \chi^\eps=(w_1^\eps+iw_{n+1}^\eps,...,w_{n}^\eps+iw_{2n}^\eps)\in
	\mathbb C^{n}.
	\end{equation}

	With these  notations, let
	$$   \Gamma^\eps_s(\chi^\eps,\theta^\eps)=\gamma^\eps(s).$$
	For further use, let cut$^\eps(x)$ be the cut-locus of  $x\in M^\eps$ in the Riemannian manifold $(M^\eps,g^\eps).$
	
	\begin{lemma}\label{epsz-Jacobian}
		Let $s,\eps\in (0,1]$  be fixed and assume that
		$\gamma^\eps(1)\notin {\rm cut}^\eps(0_{\mathbb H^n}).$ Then the
		Jacobian of $\Gamma^\eps_s$ at $(\chi^\eps,\theta^\eps)$ is
		\begin{eqnarray*}
			\Jac (\Gamma^\eps_s)(\chi^\eps,\theta^\eps) = \left\{
			\begin{array}{lll}
				2^{2n+2} s |\chi^\eps|^2 \left(\frac{\sin \frac{\theta^\eps
						s}{2}}{\theta^\eps}\right)^{2n-1} \frac{\sin\frac{\theta^\eps s}{2}
					- \frac{\theta^\eps s}{2}\cos\frac{\theta^\eps s}{2}}{\left(\theta^\eps\right)^3} + 
				2^{2n}\frac{\eps^2}{4} s \left( \frac{ \sin \frac{\theta^\eps
						s}{2}}{{\theta^\eps}} \right)^{2n} & \mbox{if }&\theta^\eps \neq 0;\\
				\frac{s^{2n+3}|\chi^\eps|^2}{3}+ s^{2n+1}\frac{\eps^2}{4} &\mbox{if
				}& \theta^\eps = 0.
			\end{array}
			\right.
		\end{eqnarray*}
	\end{lemma}
	
	{\it Proof.} {We prove the relation in case of $\theta^\eps \neq 0$. When $\theta^\eps = 0$ the formula can be obtained as a continuous limit of the previous case.} 
	
	We may decompose the differential of $\Gamma_s^\eps$
	calculated at $(\chi^\eps,\theta^\eps)$ into blocks as $$\left(
	\begin{array}{ll}
	J^{s,\eps}_{(1...2n,1...2n)} & J^{s,\eps}_{(1...2n,2n+1)} \vspace{0.20cm}\\
	J^{s,\eps}_{(2n+1,1...2n)} & J^{s,\eps}_{(2n+1,2n+1)}
	\end{array}
	\right),$$ where the components are calculated as follows:
	\begin{itemize}
		\item[*] 
		{The complex representation of the $2n \times 2n$ dimensional real matrix $J_{(1...2n,1...2n)}^{s,\eps}$ is}
		$$ i \frac{e^{-i \theta^\eps s} - 1}{\theta^\eps} I_n,$$
		where $I_n$ is
		the identity matrix in $M_n(\mathbb C)$.
		\item[*] {The column block $J_{(1...2n,
				2n+1)}^{s,\eps}$ represented as a vector in $\mathbb C^n$ is}
		$$ \left( s\frac{e^{-i\theta^\eps s}}{\theta^\eps} -
		i \frac{e^{-i\theta^\eps s}-1}{(\theta^\eps)^2} \right) \chi^\eps.$$
		\item[*] {The row block $J_{(2n+1, 1...2n)}^{s,\eps}$ can be identified with the complex representation}
		$$ {4}
		\frac{\theta^\eps s - \sin(\theta^\eps s)}{(\theta^\eps)^2}
		\chi^\eps.$$ 
		\item[*] The single element of the matrix in the lower right corner is
		$$J_{(2n+1,2n+1)}^{s,\eps} =\frac{\eps^2}{4} s+ 2  |\chi^\eps|^2
		\left( \frac{2 \sin(\theta^\eps s)-\theta^\eps s(1+ \cos(\theta^\eps
			s))}{(\theta^\eps)^3}  \right).$$
	\end{itemize}
	One can observe as in  Juillet \cite{Juillet-IMNR} that
	$\Jac(\Gamma_s^\eps)(\chi^\eps, \theta^\eps) =
	\Jac(\Gamma_s^\eps)(\chi^\eps_{0}, \theta^\eps)$ for
	$|\chi^\eps| = |\chi^\eps_{0}|$.
	{Indeed, let $W\in {\rm U}(n)$ be a unitary matrix identified with a real $2n \times 2n$ matrix and consider the linear map  $\tilde{W}(\chi,\theta) = (W\chi, \theta)$.
		Notice that $\Gamma_s^{\eps}(W \chi, \theta) = \tilde W(\Gamma_s^\eps(\chi, \theta))$; thus the  chain rule implies that ${\rm Jac}(\Gamma_s^\eps)(W \chi, \theta) = {\rm Jac(\Gamma_s^\eps)(\chi, \theta)}$. Since $|\chi^\eps| = |\chi^\eps_{0}|$, we may choose the unitary matrix $W$ such that $W \chi^\eps = \chi_0^\eps$, which proves the  claim.
	}

	In particular, we can simplify the computations by setting $\chi^\eps_0 =
	(0,...,0, |\chi^\eps|)$; in this way the above matrix has several
	zeros, and its determinant is the product of $n-1$ identical determinants
	corresponding to the matrix
	\begin{eqnarray*}
		\left(
		\begin{array}{cc}
			\frac{\sin(\theta^\eps s)}{\theta^\eps} & \frac{1 - \cos(\theta^\eps s)}{\theta^\eps}
			\vspace{0.20cm}\\
			\frac{\cos(\theta^\eps s) - 1}{\theta^\eps} & \frac{\sin(\theta^\eps
				s)}{\theta^\eps}
		\end{array}
		\right)
	\end{eqnarray*}
	with the determinant of
	\begin{eqnarray*}
		\left(
		\begin{array}{ccc}
			\frac{\sin(\theta^\eps s)}{\theta^\eps} & \frac{1 - \cos(\theta^\eps s)}{\theta^\eps} & |\chi^\eps| \left( \frac{s\cos(\theta^\eps s)}{\theta^\eps} - \frac{\sin(\theta^\eps s)}{(\theta^\eps)^2}\right) \vspace{0.20cm} \\
			\frac{\cos(\theta^\eps s) - 1}{\theta^\eps} & \frac{\sin(\theta^\eps s)}{\theta^\eps} & -|\chi^\eps| \left( \frac{s\sin(\theta^\eps s)}{\theta^\eps} + \frac{\cos(\theta^\eps s) - 1}{(\theta^\eps)^2}\right) \vspace{0.20cm} \\
			|\chi^\eps| \frac{\theta^\eps s - \sin(\theta^\eps
				s)}{(\theta^\eps)^2} & 0 & \frac{\eps^2}{4} s+ 2  |\chi^\eps|^2
			\left( \frac{2 \sin(\theta^\eps s)-\theta^\eps s(1+ \cos(\theta^\eps
				s))}{(\theta^\eps)^3}  \right)
		\end{array}
		\right).
	\end{eqnarray*}
	The rest of the computation is straightforward. \hfill $\square$\\

	For a fixed
	$s \in [0,1]$ and $(x,y) \in M^\eps \times M^\eps$ let
	\begin{eqnarray}\label{eps-intermediate}
	Z^\eps_s(x,y) = \{ z \in M^\eps : d^\eps(x,z) = s d^\eps(x,y),\
	d^\eps(z,y) = (1-s) d^\eps(x,y)\}
	\end{eqnarray}
	and
	$$
	Z^\eps_s(A,B) = \bigcup\limits_{(x,y) \in A \times B} Z^\eps_s(x,y)
	$$
	for any two nonempty subsets $A,B \subset M^\eps$. Since
	$(M^\eps,d^\eps)$ is complete, $Z^\eps_s(x,y)\neq \emptyset$ for
	every $x,y\in M^\eps.$ Let $B^\eps(y,r)=\{w\in
	M^\eps:d^\eps(y,w)<r\}$ for every $r>0.$

	Following Cordero-Erausquin, McCann and Schmuckenschl\"ager
	\cite{McCann}, we consider the volume distortion coefficient in
	$(M^\eps, g^\eps)$ as
	\begin{eqnarray*}
		&&v^\eps_s(x,y) = \lim\limits_{r \to 0}\frac{{\rm Vol}^{\eps} \left( Z^\eps_s(x, B^\eps(y, r))\right)}{{\rm  Vol}^{\eps} \left( B^\eps(y, sr)\right)} \mbox{ when } s \in (0,1].
	\end{eqnarray*}
	{Note that $v^\eps_1(x,y)=1$ for every $x,y\in \mathbb H^n$. Moreover, the local behavior of geodesic balls shows that  $v^\eps_s(x,x)=1$ for every $s\in (0,1)$ and $x\in \mathbb H^n$.}
	
	
	The following statement provides an expression for the volume
	distortion coefficient in terms of the Jacobian $\Jac
	(\Gamma_s^\eps)$.
	
	\begin{proposition}\label{property-Jacobi-representation-1}
		Let $x,y\in M^\eps$, {$x\neq y,$} and assume that $y\notin {\rm cut}^\eps(x)$. Let
		$\gamma^\eps:[0,1]\to M^\eps$ be the unique minimal geodesic joining
		$0_{\mathbb H^n}$ to $x^{-1}\cdot y$ given by
		$\gamma^\eps(s)=\exp_{0_{\mathbb H^n}}^\eps(sw^\eps)$ for some
		$w^\eps=\sum_{j=1}^nw_j^\eps X_j(0_{\mathbb
			H^n})+\sum_{j=1}^nw_{j+n}^\eps Y_j(0_{\mathbb H^n})+w_{2n+1}^\eps
		T^\eps(0_{\mathbb H^n})\in T_{0_{\mathbb H^n}} M^\eps$. Then for
		every $s\in (0,1)$ we have
		\begin{itemize}
			\item[(i)] $\displaystyle v_s^\epsilon(x,y)=\frac{1}{s^{2n+1}}\frac{{\rm
					Jac}(\Gamma_{s}^\eps)(\chi^\eps,\theta^\eps)}{{\rm
					Jac}(\Gamma_1^\eps)(\chi^\eps, \theta^\eps)};$\\
			
			\item[(ii)] $\displaystyle v_{1-s}^\epsilon(y,x) =\frac{1}{(1-s)^{2n+1}}\frac{{\rm
					Jac}(\Gamma_{1-s}^\eps)(\chi^\eps ,\theta^\eps)}{{\rm
					Jac}(\Gamma_1^\eps)(\chi^\eps ,
				\theta^\eps)},$
		\end{itemize}
		where $\theta^\eps$ and $\chi^\eps$  come from {\rm
			(\ref{jeloles-chi})}.
	\end{proposition}
	
	{\it Proof.} Since $y\notin {\rm cut}^\eps(x)$ and ${\rm
		cut}^\eps(x)$ is closed, there exists $r>0$ small enough such that
	$B^\eps(y,r)\cap {\rm cut}^\eps(x)=\emptyset.$ In particular, the
	point $x$ and every element from $B^\eps(y,r)$ can be joined  by a
	unique minimal $\eps$-geodesic and $Z_s^\eps(x, z)$ is a
	singleton for every $z\in B^\eps(y,r).$ By the left-translation
	(valid also on the $(2n+1)-$dimensional Riemannian manifold
	$(M^\eps,g^\eps)$), we observe that  $Z_s^\eps(x,z) = x \cdot
	Z_s^\eps(0_{\mathbb H^n}, x^{-1} \cdot z)$ for all $z \in
	B^\eps(y,r).$ Thus, 
	\begin{eqnarray*}
		v_s^\eps(x,y)  &=& \lim\limits_{r \to 0}\frac{{\rm Vol}^\eps \left( x\cdot  Z_s^\eps(0_{\mathbb H^n},
			B^\eps(x^{-1}\cdot y, r))\right)}{{\rm Vol}^\eps \left(x\cdot
			B^\eps(x^{-1}\cdot y, sr)\right)} = \lim\limits_{r \to 0}\frac{{\rm
				Vol}^\eps \left( Z_s^\eps(0_{\mathbb H^n}, B^\eps(x^{-1}\cdot y,
			r))\right)}{ {\rm Vol}^\eps \left( B^\eps(x^{-1}\cdot y, sr)\right)} \\
		&=& \lim\limits_{r \to 0}\frac{{\rm Vol}^\eps \left(
			B^\eps(x^{-1}\cdot y, r)\right)}{ {\rm Vol}^\eps \left(
			B^\eps(x^{-1}\cdot y, sr)\right)}\frac{{\rm Vol}^\eps \left( Z_s^\eps(0_{\mathbb H^n},
			B^\eps(x^{-1}\cdot y, r))\right)}{ {\rm Vol}^\eps \left(
			B^\eps(x^{-1}\cdot y, r)\right)}.
	\end{eqnarray*}
	Because of the asymptotic behaviour of the volume of small balls in the
	Riemannian geometry (see Gallot,  Hulin and  Lafontaine \cite{GHL}), we have
	\begin{eqnarray*}
		v_s^\eps(x,y) &=& \frac{1}{s^{2n+1}}\lim\limits_{r \to 0}\frac{{\rm
				Vol}^\eps \left( Z_s^\eps(0_{\mathbb H^n}, B^\eps(x^{-1}\cdot y,
			r))\right)}{ {\rm Vol}^\eps \left( B^\eps(x^{-1}\cdot y, r)\right)}
		\\&=& \frac{1}{s^{2n+1}}\lim\limits_{r \to 0}\frac{\Lmeas^{2n+1} \left(
			Z_s^\eps(0_{\mathbb H^n}, B^\eps(x^{-1}\cdot y, r))\right)}{
			\Lmeas^{2n+1} \left( B^\eps(x^{-1}\cdot y, r)\right)}.
	\end{eqnarray*}
	In the last step we used ${\rm d}\textsf{m}_{\eps}=\frac{1}{\eps}{\rm
		d}\mathcal L^{2n+1}$. The rest of the proof goes in the same way as
	in case of Proposition \ref{proposition-repres-jacobi} (i); see also Cordero-Erausquin, McCann and Schmuckenschl\"ager
	\cite{McCann}.
	


	(ii) Taking into account that
	$$v_{1-s}^\eps(y,x)=v_{1-s}^\eps(0_{\mathbb H^n},y^{-1}\cdot
	x)=v_{1-s}^\eps(0_{\mathbb H^n},-x^{-1}\cdot y)$$ and the
	$\eps$-geodesic joining $0_{\mathbb H^n}$  and $-x^{-1}\cdot y$ is
	given by the curve  $s\mapsto \Gamma^\eps_s\left(-\chi^\eps e^{-i
		\theta^\eps},-\theta^\eps\right),\ \ s\in [0,1],$ a similar
	argument works as in Proposition \ref{proposition-repres-jacobi}
	(ii). \hfill
	$\square$

	\subsection{Optimal mass transportation on $\mathbb
		H^n$ and $M^\eps$}\label{subsection-optimal-mass} Let us fix two
	functions $f, g:\mathbb H^n\to [0,\infty)$ and assume  that
	$$\int_{\mathbb H^n}f=\int_{\mathbb H^n}g=1.$$ Let $\mu_0=f \mathcal
	L^{2n+1}$ and $\mu_1=g \mathcal L^{2n+1}$. By the theory of
	optimal mass transportation on $\mathbb H^n$ for $c=d_{CC}^2/2$, see Ambrosio and Rigot 
	\cite[Theorem 5.1]{AR}, there exists a unique optimal transport map
	from $\mu_0$ to $\mu_1$ which is induced by the map
	\begin{equation}\label{optimal-map}
	\psi(x)=x\cdot \Gamma_1(-X \varphi(x)-iY \varphi(x),-4T\varphi(x)) \
	{\rm  a.e.}\ x\in {\rm supp}f,
	\end{equation}
	for some $c$-concave and locally
	Lipschitz map $\varphi,$ where $\Gamma_1$ comes from
	(\ref{Gamma-s-definicio}). In fact, according to Figalli and Rifford
	\cite{FR}, there exists a Borel set $C_0 \subset {\rm
		supp}f$ of null $\mathcal L^{2n+1}$-measure  such that for every $x\in {\rm
		supp}f\setminus C_0$, there exists a unique minimizing geodesic from
	$x$ to $\psi(x)$; this geodesic is represented by
	{\begin{eqnarray} \label{DefIntMapH}
		s\mapsto \psi_s(x):=
		\left\{
		\begin{array}{lll}
		x\cdot \Gamma_s(-X \varphi(x)-iY \varphi(x),-4T\varphi(x))
		\  &\mbox{if} &  \psi(x)\neq x; \\
		x &\mbox{if} &  \psi(x)=x.
		\end{array}\right.
		\end{eqnarray}
		The sets $\mathcal M_\psi=\{x\in \mathbb H^n:\psi(x)\neq x\}\ {\rm  and }\ \mathcal S_\psi=\{x\in \mathbb H^n:\psi(x)= x\}$ correspond to the moving and static sets of the transport map $\psi$, respectively.} 

	On the Riemannian manifold $(M^{\eps},g^{\eps})$, we may consider
	the unique optimal transport map $\psi^{\eps}$ from $\mu_0^{\eps}=(\eps f) \textsf{m}_{\eps}$ to $\mu_1^{\eps}=(\eps
	g) \textsf{m}_{\eps}$. The existence and uniqueness of such a map is
	provided by McCann \cite[Theorem 3.2]{McCann-GAFA}. This 
	map is defined by a $c_{\eps}=(d^\eps)^2/2$-concave function
	$\varphi_{\eps}$ via
	$$\psi^{\eps}(x) = \exp^{\eps}_x (-\nabla^{\eps}
	\varphi_{\eps}(x)) = x \cdot \Gamma_1^\eps(-X \varphi_\eps(x) - i Y
	\varphi_\eps(x), -4T \varphi_\eps(x)) \ {\rm a.e.}\ x \in {\rm
		supp}f,$$ where $$
	\nabla^{\eps} \varphi_{\eps}(x)=\sum_{j=1}^n
	X_j\varphi_{\eps}(x)X_j(x)+Y_j\varphi_{\eps}(x)Y_j(x)+T^{\eps}\varphi_{\eps}(x)T^{\eps}(x)\in
	T_xM^{\eps}, $$ see Ambrosio and Rigot \cite[p. 292]{AR}. Note that we may always assume
	that $\varphi_{\eps}(0_{\mathbb H^n})=0.$ Due to Cordero-Erausquin, McCann and Schmuckenschl\"ager \cite[Theorem
	4.2]{McCann}, there exists a Borel set $C_{\eps} \subset {\rm
		supp}f$   of null $\textsf{m}_{\eps}$-measure such that $\psi^{\eps}(x) \notin {\rm
		cut}^{\eps}(x)$ for every $x \in {\rm supp}f \backslash C_{\eps}$.
	Now we consider the interpolant map 
	
	\begin{equation}\label{imapeps} 
	\psi_s^{\eps}(x) =\exp^{\eps}_x (-s\nabla^{\eps} \varphi_{\eps}(x)),\ \ x\in{\rm
		supp}f \backslash C_{\eps}.
	\end{equation}Ê
	
	Using again a left-translation, we
	equivalently have $$\psi_s^{\eps}(x) =x\cdot  \exp^{\eps}_0
	(-sw^{\eps}(x))= x \cdot \Gamma_s^\eps(-X \varphi_\eps(x) - i Y
	\varphi_\eps(x), -4T \varphi_\eps(x)),\ \ x\in{\rm supp}f \backslash
	C_{\eps},$$ where
	\begin{equation}\label{epsz-gradient-2}
	w^{\eps}(x)=\sum_{j=1}^n X_j\varphi_{\eps}(x)X_j(0_{\mathbb
		H^n})+Y_j\varphi_{\eps}(x)Y_j(0_{\mathbb
		H^n})+T^{\eps}\varphi_{\eps}(x)T^{\eps}(0_{\mathbb H^n})\in
	T_{0_{\mathbb H^n}}M^{\eps}.
	\end{equation}
	
	With the above notations we summarize the results in this section,
	establishing a bridge between notions in $\mathbb H^n$ and $M^\eps$ which will be crucial in the proof of our main
	theorems:
	
	\begin{proposition}\label{bridge-proposition} There exists a sequence $\{ \eps_k\}_{k \in \N}\subset (0,1)$ 
		converging to $0$ and a full $\mu_{0}$-measure set $D\subset \mathbb H^n$ such that $f $ is positive on $D$ and
		for every $x\in D$ we have:
		\begin{itemize}
			\item[{\rm (i)}] $\lim\limits_{k
				\to \infty} \psi^{\eps_k}_s(x) = \psi_s(x)$ for every $s\in (0,1];$
			{\item[{\rm (ii)}] $\liminf\limits_{k \to
					\infty} v_s^{\eps_k}(x,\psi^{\eps_k}(x))\geq v^0_s(x,\psi(x))$ for every $s\in (0,1);$ 
				\item[{\rm (iii)}] $  \liminf\limits_{k \to \infty}
				v_{1-s}^{\eps_k}(\psi^{\eps_k}(x),x)\geq v^0_{1-s}(\psi(x),x)$ for every $s\in (0,1)$.
			}
		\end{itemize}
	\end{proposition}
	
	\begin{remark}\rm 
	Note that the limiting value of the distortion coefficients in the Riemannian approximation (i.e., (ii) and (iii)) are {\it not} the Heisenberg volume distortion coefficients $v_s(x,y)$. The appropriate limits are given by $v_s^0(x,y)$, see (\ref{Heisenberg-distortion}).  
	\end{remark}
	
	{\it Proof of Proposition \ref{bridge-proposition}.}  Let us start with an arbitrary sequence $\{ \eps_k\}_{k \in \N}$ of positive numbers
	such that  $\lim_{k\to \infty}\eps_k=0$ and $C = C_0 \cup \left(\cup_{k \in \N}
	C_{\eps_k}\right)$, where $C_0$ and $C_{\eps_k}$ are the sets with
	null $\mathcal L^{2n+1}$-measure coming from  the previous
	construction, i.e., there is a unique minimizing geodesic from $x$
	to $\psi(x)$ and $\psi^{\eps_k}(x) \notin {\rm cut}^{\eps_k}(x)$ for
	every $x\in {\rm supp}f\setminus C$. We define $D= \{x\in \mathbb H^n : f(x) > 0\} \setminus C.$ 
	Notice that $D$ has full $\mu_0$-measure by its definition.
	It is clear that every volume distortion
	coefficient appearing in (ii) and (iii) is well-defined for every
	$x\in D$. The set $D$ from the claim will be obtained in the course of the proof by subsequently discarding
	null measure sets several times from $D$. In order to simplify the notation we shall keep the notation $D$ for the 
	sets that are obtained in this way. Similarly, we shall keep the notation for  $\{\eps_k\}_{k\in \mathbb N}$ when we pass to a subsequence.

	
	Accordingly,  by Ambrosio and Rigot \cite[Theorem 6.2]{AR} we have
	that 
	$$\lim\limits_{k
		\to \infty} \psi^{\eps_k}(x) = \psi(x)\  \text{for a.e.} \ x \in D.$$
	
	In the proof of (i) we shall distinguish two cases. Let $s\in (0,1]$ be fixed.  
	
	{\bf Case 1:} {\it  the moving set $ \mathcal M_\psi.$} By using \cite[Theorem
	6.11]{AR} of Ambrosio and Rigot,  up to the removal of a null measure set and up to passing to a subsequence 
	we have
	\begin{equation}\label{limits-1}
	\lim_{k\to \infty} X_j\varphi_{\eps_k}(x)=X_j\varphi(x), \ \
	\lim_{k\to \infty} Y_j\varphi_{\eps_k}(x)=Y_j\varphi(x), \ \
	\lim_{k\to \infty} T\varphi_{\eps_k}(x)=T\varphi(x),
	\end{equation}
	where $\varphi_{\eps_k}$ and $\varphi$
	are the $c_{\eps_k}$-concave and $c$-concave functions appearing in
	(\ref{optimal-map}) and (\ref{epsz-gradient-2}). 
	
	Due to the
	form of $w^{\eps_k}(x) \in T_{0_{\mathbb H^n}}M^{\eps_k}$ from
	(\ref{epsz-gradient-2}), we introduce the complex vector-field $\chi^{\eps_k}= (\chi^{\eps_k}_{1}, \ldots , \chi^{\eps_k}_{n})$ by 
	$
	\chi^{\eps_k}_{j}(x)=X_{j}\varphi_{\eps_k}(x)+iY_{j}\varphi_{\eps_k}(x).
	$
	Let also $ w_{2n+1}^{\eps_k}(x)=T^{\eps_k}\varphi_{\eps_k}(x)$.

	The limits
	in (\ref{limits-1}) imply that for a.e. $x \in D \cap \mathcal M_\psi$ we have 
	\begin{equation}\label{limit-alap-1-1}
	\lim_{k\to \infty}\chi^{\eps_k}= X\varphi(x)+iY\varphi(x),
	\end{equation}
	\begin{eqnarray}\label{limit-alap-2-1}
	\lim_{k\to \infty} \theta^{\eps_k} &=&\lim_{k\to \infty}\frac{4w_{2n+1}^{\eps_k}}{\eps_k}=  \lim_{k\to \infty}\frac{4T^{\eps_k}\varphi_{\eps_k}(x)}{\eps_k}=
	\lim_{k\to
		\infty}\frac{4{\eps_k}T\varphi_{\eps_k}(x)}{\eps_k}=4\lim_{k\to
		\infty} T\varphi_{\eps_k}(x)\nonumber  \\
	&=&4T\varphi(x).
	\end{eqnarray}
	From the representations (\ref{eps-geodetikus}) and
	(\ref{explicit-geodet}) of the $\eps_k$-geodesics and Heisenberg
	geodesics, respectively, relations (\ref{limit-alap-1-1}) and
	(\ref{limit-alap-2-1}) imply that
	\begin{eqnarray*}
		\lim\limits_{k
			\to \infty} \psi^{\eps_k}_s(x)&=&\lim\limits_{k \to \infty}x\cdot
		\exp^{\eps_k}_0 (-sw^{\eps_k}(x))=
		x\cdot \Gamma_s(-X \varphi(x)-iY \varphi(x),-4T\varphi(x))= \psi_s(x).
	\end{eqnarray*}
	
	{\bf Case 2:} {\it  the static set $ \mathcal S_\psi.$}  From the representation (\ref{DefIntMapH}) we have that $\psi_{s}(x)  = x$ for any $x \in \mathcal S_\psi$.
	Clearly, we only need to consider values of $\eps_{k}$ for which $\psi^{\eps_{k}}(x) \neq x$.  Again, by \cite[Theorem 6.2]{AR} of Ambrosio and Rigot,
	$\lim\limits_{k
		\to \infty} \psi^{\eps_k}(x) = \psi(x)=x\  \text{for a.e.} \ x \in D \cap \mathcal S_\psi.$ According to (\ref{imapeps}) the point $\psi^{\eps_k}_{s}(x)$ lies on the $\eps_{k}$-geodesic 
	connecting $x$ and $\psi^{\eps_k}(x)$. The latter limit and the estimate (\ref{dist-eps-Juillet}) imply 
	the following chain of inequalities
	\begin{eqnarray*}
		{
			s \left( d_{CC}\left(x, \psi^{\eps_k}(x)\right) - c \pi \eps_k\right) \leq s d^{\eps_k}\left(x, \psi^{\eps_k}(x)\right) = d^{\eps_k} \left(x, \psi^{\eps_k}_s(x)\right) \leq d_{CC}\left(x, \psi^{\eps_k}_s(x)\right) \leq} \\
		{\leq d^{\eps_k}\left(x, \psi^{\eps_k}_s(x)\right) + c \pi \eps_k = s d^{\eps_k} \left(x, \psi^{\eps_k}(x)\right) + c\pi\eps_k \leq s d_{CC}\left(x, \psi^{\eps_k}(x)\right) + c\pi \eps_k,
		}
	\end{eqnarray*}
	{so}
	$\lim\limits_{k
		\to \infty}\psi^{\eps_k}_{s}(x) = x$, which ends the 
	proof of (i). 
	
	\medskip

	To prove inequality (ii) we distinguish again two cases.  
	
	{\bf Case 1:} {\it the moving set $ \mathcal M_\psi.$} Let $x\in D\cap \mathcal M_\psi.$ Since $\lim\limits_{k
		\to \infty} \psi^{\eps_k}(x)=\psi(x)\neq x$, there exists $k_0\in \mathbb N$ such that  $\psi^{\eps_k}(x)\neq x$ for every $k\geq k_0$.  Thus,  we have 
	\begin{eqnarray*}
		\lim_{k\to \infty}v^{\eps_k}_s(x,\psi^{\eps_k}(x)) &=& \frac{1}{s^{2n+1}}\lim_{k\to \infty}\frac{{\rm
				Jac}(\Gamma_{s}^{\eps_k})(-\chi^{\eps_k},-\theta^{\eps_k})}{{\rm
				Jac}(\Gamma_1^{\eps_k})(-\chi^{\eps_k},-\theta^{\eps_k})}
		\ \ \ \ \ \ \ \ \ \ \ \ \  \    ({\rm cf.\ Proposition\
			\ref{property-Jacobi-representation-1}})\\
		&=&\frac{1}{s^{2n+1}} \frac{{\rm
				Jac}(\Gamma_{s})(-X\varphi(x)-iY\varphi(x),-4T\varphi(x))}{{\rm
				Jac}(\Gamma_1)(-X\varphi(x)-iY\varphi(x),-4T\varphi(x))}\\
		&& \ \ \ \ \ \ \ \ \ \ \ \ \ \ \ \ \ \  \ \ \ \ \ \ \ \ \ \ \ \ \ \
		\ \ \ \ \ \ \ \ \ \ \  \ \ \ \ ({\rm cf.\ Lemma\
			\ref{epsz-Jacobian}}\ \&\ (\ref{limit-alap-1-1}),
		(\ref{limit-alap-2-1})) \\
		&=&sv_s(x,\psi(x))\ \ \ \ \ \ \ \ \ \ \ \ \ \ \ \ \ \ \ \ \ \  \ \ \ \ \ \ \ \ \    (x^{-1}\cdot \psi(x)\notin L\ \& {\rm\ Proposition\ \ref{proposition-repres-jacobi}})\\
		&=&
		v_s^0(x,\psi(x)). \ \ \ \ \ \ \ \ \ \ \ \ \ \ \ \ \ \ \ \ \ \  \ \ \ \ \ \ \ \ \ \ \ \ \ \ \ \ \ \ \ \ \ \ \ \ \  \ \  ({\rm cf.\ (\ref{Heisenberg-distortion})})
	\end{eqnarray*}
	
	{\bf Case 2:} {\it the static set $ \mathcal S_\psi.$} Let $x\in D\cap \mathcal S_\psi.$ If $\psi^{\eps_k}(x)= x$ then by (\ref{Heisenberg-distortion}) we have  $v_s^{\eps_k}(x,\psi^{\eps_k}(x))=1\geq s^2=v_s^0(x,\psi(x))$. If $\psi^{\eps_k}(x)\neq  x$, by Proposition  \ref{property-Jacobi-representation-1} and Lemma \ref{epsz-Jacobian}, we have 
	$$v^{\eps_k}_s(x,\psi^{\eps_k}(x))=\frac{1}{s^{2n+1}}\frac{A_s^{\eps_k}+B_s^{\eps_k}}{A_1^{\eps_k}+B_1^{\eps_k}},$$
	where 
	\begin{eqnarray*}
		A_s^\eps = \left\{
		\begin{array}{lll}
			2^{2n+2} s |\chi^\eps|^2 \left(\frac{\sin \frac{\theta^\eps
					s}{2}}{\theta^\eps}\right)^{2n-1} \frac{\sin\frac{\theta^\eps s}{2}
				- \frac{\theta^\eps s}{2}\cos\frac{\theta^\eps s}{2}}{\left(\theta^\eps\right)^3} & \mbox{if }&\theta^\eps \neq 0;\\
			\frac{s^{2n+3}|\chi^\eps|^2}{3} &\mbox{if
			}& \theta^\eps = 0,
		\end{array}
		\right.
	\end{eqnarray*}
	and 
	\begin{eqnarray*}
		B_s^\eps = \left\{
		\begin{array}{lll}
			2^{2n}\frac{\eps^2}{4} s \left( \frac{ \sin \frac{\theta^\eps
					s}{2}}{{\theta^\eps}} \right)^{2n} & \mbox{if }&\theta^\eps \neq 0;\\
			s^{2n+1}\frac{\eps^2}{4} &\mbox{if
			}& \theta^\eps = 0.
		\end{array}
		\right.
	\end{eqnarray*}
	The elementary inequality $\sin(\alpha s) \geq s\sin(\alpha)$ for $\alpha \in [0, \pi]$ and $s \in [0,1]$ shows that $B_s^\eps \geq s^{2n+1}B_1^\eps$. By Lemma \ref{lemma-novekvo}, 
	$A_s^\eps\geq s^{2n+3}A_1^\eps.$
	Therefore, the above inequalities imply that $$v^{\eps_k}_s(x,\psi^{\eps_k}(x))\geq s^2=v_s^0(x,\psi(x)),$$ which concludes the proof.

	Claim (iii) for  $v^0_{1-s}(\psi(x),x)$ is proven similarly as claim (ii) for $v^0_{s}(\psi(x),x)$.
	{ \hfill
		$\square$ }
	
	{\begin{remark}\rm 
			In the second case of the above proof (i.e., $x\in \mathcal S_\psi$) we could expect a better lower bound than $s^2$ for $v^{\eps_k}_s(x,\psi^{\eps_k}(x))$  as $k\to \infty$ since no explicit presence of Heisenberg volume distortion coefficient is expected. However, in the general case  $s^2$ is the optimal bound. Indeed, since $x\in \mathcal S_\psi$ we first notice that  $|\chi^{\eps_k}|\to 0$ as $k\to \infty$. Thus, if $\theta^{\eps_k}\to 0$ and we assume that  $|\chi^{\eps_k}|=\mathcal O(\eps_k^\alpha)$ as $k\to \infty$ for some $0<\alpha<1$, we have   $$\liminf_{k\to \infty}v^{\eps_k}_s(x,\psi^{\eps_k}(x))=s^2.$$
		\end{remark}
	}
	\medskip\medskip
	
	\section{Proof of main results}
	
	
	\subsection{Jacobian determinant inequality on $\mathbb H^n$} In this subsection we shall prove our Jacobian determinant inequality on $\mathbb H^n$ as the key result of the paper.  \medskip
	
	\noindent {\bf Proof of Theorem \ref{TJacobianDetIneq}.}
	We shall consider the sequence $\{ \eps_k\}_{k \in \N}\subset (0,1]$ such that $\lim_{k\to \infty}\eps_k=0$ and the statement of Proposition \ref{bridge-proposition} holds. Let $(M^{\eps_k},g^{\eps_k})$ be the Riemannian 
	manifolds  approximating $\mathbb H^n$, $k\in \mathbb N.$  
	
	Let us consider the measures
	$\mu_0=\rho_0 \mathcal L^{2n+1}$, $\mu_1=\rho_1 \mathcal
	L^{2n+1}$, $\mu_0^{\eps_k}=(\eps_k \rho_0) \textsf{m}_{\eps_k}=\mu_0$,
	$\mu_1^{\eps_k}=(\eps_k \rho_1) \textsf{m}_{\eps_k}=\mu_1$, the associated optimal transport maps $\psi, \psi^{\eps_{k}}$ and their interpolants $\psi_s$ and $\psi_s^{\eps_k}$, respectively.
	Let us keep the notations $C$ and $D$ from \S \ref{subsection-optimal-mass} and  Proposition \ref{bridge-proposition}. According to Figalli and Rifford \cite[Theorem 3.7]{FR}, Figalli and Juillet \cite[p. 136]{FJ} and Cordero-Erausquin, McCann and Schmuckenschl\"ager \cite[Lemma 5.3]{McCann}, the maps $\psi$, $\psi_s$  and $\psi_s^{\eps_k}$  are essentially injective on $D$, respectively. Consequently, there is a set $D_0\subset D$ of null $\mathcal L^{2n+1}$-measure such that the maps  $\psi$, $\psi_s$ and $\psi_s^{\eps_k}$ $(k\in \mathbb N)$ 
	are injective on $D\setminus D_0$; for simplicity, we keep the notation $D$ for $D\setminus D_0$. 
	Let $\mu_s=(\psi_s)_\#\mu_0$ and $\mu_s^{\eps_k} =(\psi_s^{\eps_k})_{\#} \mu_0$ be the push-forward measures on $\mathbb H^n$ and $M^{\eps_k}$, and $\rho_s$ and $\eps_k\rho_s^{\eps_k}$ be their density functions w.r.t. to the measures $\mathcal L^{2n+1}$ and $\textsf{m}_{\eps_k}$,  respectively.

	Let $A_i\subset \mathbb H^n$ be the support of the measures $\mu_i$, $i\in \{0,1\}.$   On account of (\ref{dist-eps-Juillet}), definition (\ref{eps-intermediate}) and the compactness of the sets $A_0$ and $A_1$, one has for every $x\in D$ that 
	\begin{eqnarray}\label{d-cc-d-epsz-estimate}
	{d^{\eps_k}(x,\psi_s^{\eps_k}(x)) = sd^{\eps_k}(x,\psi^{\eps_k}(x)) \leq sd_{CC}(x,\psi^{\eps_k}(x)) \leq s\max_{(x,y)\in A_0\times A_1} d_{CC}(x,y).}
	\end{eqnarray}
	Since by (\ref{dist-eps-Juillet}) we have that 
	$$d_{CC}(x,\psi_s^{\eps_k}(x)) \leq  d^{\eps_k}(x,\psi_s^{\eps_k}(x)) + \eps_{k}c \pi, $$ 
	the estimate   (\ref{d-cc-d-epsz-estimate})
	assures the existence of $R>0$ such that the ball $B(0,R)$  contains the supports of the measures $\mu_s=(\psi_s)_\#\mu_0$ and $\mu_s^{\eps_k} = (\psi_s^{\eps_k})_{\#} \mu_0^{\eps_k}=(\psi_s^{\eps_k})_{\#} \mu_0$,  $k \in \N$; in fact, we can choose 
	\begin{equation}\label{R-sugar}
	{R=\max_{x\in A_0} d_{CC}(0_{\mathbb H^n},x)+\max_{(x,y)\in A_0\times A_1} d_{CC}(x,y)+1.}
	\end{equation}
	Clearly, $A_0,A_1\subset B(0,R)$. 
	Thus, it is enough to take $\textsf{m}=\mathcal L^{2n+1}|_{B(0,R)}$ as the reference measure. 
	
	The proof is based on the Jacobian determinant inequality from
	\cite[Lemma 6.1]{McCann} on $M^{\eps_k}$, i.e.,  for every $x\in D$,
	\begin{eqnarray*}
		\left({\rm
			Jac}(\psi^{\eps_k}_s)(x)\right)^\frac{1}{2n+1} &\geq&
		(1-s)\left(v^{\eps_k}_{1-s}(\psi^{\eps_k}(x),x)\right)^\frac{1}{2n+1}
		+s\left(v^{\eps_k}_{s}(x,\psi^{\eps_k}(x))\right)^\frac{1}{2n+1}\left({\rm
			Jac}(\psi^{\eps_k})(x)\right)^\frac{1}{2n+1}.
	\end{eqnarray*}
	
	The technical difficulty is that we cannot simply pass to a point-wise limit in the latter inequality because we do not have 
	an almost everywhere convergence of Jacobians. To overcome this issue we aim to prove a weak version of the inequality 
	by multiplying by a continuous test function and integrating. As we shall see in the sequel, this trick allows the process of 
	passing to the limit and we can obtain an integral version of the Jacobian inequality which in turn gives us the desired point-wise 
	inequality almost everywhere.

	To carry out the aforementioned program, we combine the above Jacobian determinant inequality with the Monge-Amp\`ere equations on $M^{\eps_k}$, namely,
	\begin{equation}\label{MAEq-1}
	\eps_k \rho_0(x)=\eps_k \rho_1(\psi^{\eps_k}(x)){\rm Jac}(\psi^{\eps_k})(x),\ \  \eps_k \rho_0(x)=\eps_k \rho^{\eps_k}_s(\psi_s^{\eps_k}(x)){\rm Jac}(\psi_s^{\eps_k})(x),\ \ x\in D.
	\end{equation}
	Thus, we obtain for every $x \in D$ that 
	\begin{eqnarray}\label{jacobi-inequality-eps}
	\left(\rho_s^{\eps_k} (\psi_s^{\eps_k}(x))\right)^{-\frac{1}{2n+1}} \geq  (1-s) (v_{1-s}^{\eps_k}(\psi^{\eps_k}(x),x))^{\frac{1}{2n+1}} (\rho_0(x))^{-\frac{1}{2n+1}} \nonumber \\
	+ s (v_s^{\eps_k}(x, \psi^{\eps_k}(x))^{\frac{1}{2n+1}} (\rho_1(\psi^{\eps_k}(x)))^{-\frac{1}{2n+1}}.
	\end{eqnarray}
	{Let us fix an arbitrary non-negative test function $h\in C_c(\mathbb H^n)$ with support in $B(0,R)$;  for simplicity of notation, let $S={\rm supp}(h)$.} 
	Multiplying (\ref{jacobi-inequality-eps}) by $h(\psi_s^{\eps_k}(x))
	\geq 0$, an integration on $D$  w.r.t. the measure $\mu_0=\rho_0 \textsf{m}$ gives
\begin{equation}\label{left-right}
L_s^k\geq R_{s,1}^k+R_{s,2}^k,
\end{equation}
	where 
	$$L_s^k:= \int_{D} h(\psi_s^{\eps_k}(x)) \left(\rho_s^{\eps_k} (\psi_s^{\eps_k}(x))\right)^{-\frac{1}{2n+1}} {\rho_0(x)} \dd \textsf{m}(x),$$
	$$R_{s,1}^k:=\int_D h(\psi_s^{\eps_k}(x))(1-s) (v_{1-s}^{\eps_k}(\psi^{\eps_k}(x),x))^{\frac{1}{2n+1}} (\rho_0(x))^{1-\frac{1}{2n+1}} \dd \textsf{m}(x),$$
	and 
	$$R_{s,2}^k:=\int_D h(\psi_s^{\eps_k}(x)) s(v_s^{\eps_k}(x, \psi^{\eps_k}(x)))^{\frac{1}{2n+1}} (\rho_1(\psi^{\eps_k}(x)))^{-\frac{1}{2n+1}} {\rho_0(x)}\dd \textsf{m}(x).$$
Note that by Fatou's lemma, the continuity of $h$ and Proposition \ref{bridge-proposition}, we have
\begin{eqnarray}\label{right-1}
\liminf\limits_{k \to \infty}R_{s,1}^k &=& \liminf\limits_{k \to \infty}\int_D h(\psi_s^{\eps_k}(x))(1-s) (v_{1-s}^{\eps_k}(\psi^{\eps_k}(x),x))^{\frac{1}{2n+1}} (\rho_0(x))^{1-\frac{1}{2n+1}} \dd \textsf{m}(x)\nonumber \\&\geq & \int_D h(\psi_s(x))(1-s) (v_{1-s}^0(\psi(x),x))^{\frac{1}{2n+1}} (\rho_0(x))^{1-\frac{1}{2n+1}} \dd \textsf{m}(x).
\end{eqnarray} 
By the Monge-Amp\`ere equations (\ref{MA-bevezeto}) and (\ref{MAEq-1}), it turns out that for every $k\in \mathbb N$ we have $ \psi^{\eps_k}(D)= \psi(D)={\rm supp}(\mu_1)$ (up to a null measure set). Therefore, by performing a change of variables $y=\psi^{\varepsilon_k}(x)$
in the integrand $R_{s,2}^k$, we obtain by (\ref{MAEq-1}) that
\begin{eqnarray*}
	R_{s,2}^k &= &	\int_D h(\psi_s^{\eps_k}(x)) s(v_s^{\eps_k}(x, \psi^{\eps_k}(x)))^{\frac{1}{2n+1}} (\rho_1(\psi^{\eps_k}(x)))^{-\frac{1}{2n+1}} {\rho_0(x)}\dd \textsf{m}(x)\nonumber \\&=& 	\int_{\psi^{\eps_k}(D)} h(\psi_s^{\eps_k}\circ (\psi^{\eps_k})^{-1}(y)) s(v_s^{\eps_k}((\psi^{\eps_k})^{-1}(y), y))^{\frac{1}{2n+1}} (\rho_1(y))^{1-\frac{1}{2n+1}} \dd \textsf{m}(y)\nonumber \\&=&\int_{\psi(D)} h(\psi_s^{\eps_k}\circ (\psi^{\eps_k})^{-1}(y)) s(v_s^{\eps_k}((\psi^{\eps_k})^{-1}(y), y))^{\frac{1}{2n+1}} (\rho_1(y))^{1-\frac{1}{2n+1}} \dd \textsf{m}(y). 
\end{eqnarray*}
Taking the lower limit as $k\to \infty,$  Fatou's lemma, the continuity of $h$ and  Proposition \ref{bridge-proposition}  imply that 
$$	\liminf\limits_{k \to \infty}R_{s,2}^k\geq \int_{\psi(D)} h(\psi_s\circ \psi^{-1}(y)) s(v_s^0(\psi^{-1}(y), y))^{\frac{1}{2n+1}} (\rho_1(y))^{1-\frac{1}{2n+1}} \dd \textsf{m}(y).$$
Changing back the variable $y=\psi(x)$, it follows by (\ref{MA-bevezeto}) that
\begin{equation}\label{right-2}
\liminf\limits_{k \to \infty}R_{s,2}^k\geq \int_Dh(\psi_s(x)) s (v_s^0(x, \psi(x)))^{\frac{1}{2n+1}} (\rho_1(\psi(x)))^{-\frac{1}{2n+1}} {\rho_0(x)} \dd \textsf{m}(x).
\end{equation}
	By Corollary \ref{corollary-1}, relations (\ref{concentration}), (\ref{kapcsolat-volume-es-CD}) and (\ref{tauk-kozotti-osszefugges}), we observe that for every $x\in D$,
	\begin{eqnarray*}
		(1-s)[v_{1-s}^0(\psi(x),x)]^{\frac{1}{2n+1}} = \tau_{1-s}^n(\theta_x) \ \ {\rm and}\ \
		s[ v_s^0(x, \psi(x)))]^{\frac{1}{2n+1}} = \tau_{s}^n(\theta_x).
	\end{eqnarray*}
	Therefore, by the estimates (\ref{right-1}) and (\ref{right-2}) we obtain
	\begin{eqnarray}\label{elo-entropia}
	\liminf\limits_{k \to \infty} L_s^k &\geq & 
	\int_D h(\psi_s(x))\tau_{1-s}^n(\theta_x) \rho_0(x)^{1-\frac{1}{2n+1}} \dd \textsf{m}(x) \nonumber \\
	&& + \int_{D}h(\psi_s(x))\tau_{s}^n(\theta_x) (\rho_1(\psi(x)))^{-\frac{1}{2n+1}} {\rho_0(x)} \dd \textsf{m}(x).
	\end{eqnarray}

	In the sequel, we shall prove that 
	\begin{equation}\label{entropiak-weak-1}
	\int_D h(\psi_s(x))\rho_s(\psi_s(x))^{-\frac{1}{2n+1}}\rho_0(x)\dd \textsf{m}(x)\geq \liminf\limits_{k \to \infty} L_s^k.
	\end{equation}
	
	Let us notice first  that $\mu_s^{\eps_k} \rightharpoonup \mu_s$ as  $k \to \infty$. Indeed, let $\varphi:\mathbb H^{n} \to \R$ be a continuous test function with support in $B(0,R)$. By the definition of interpolant measures    $\mu_s^{\eps_k} = (\psi_s^{\eps_k})_{\#} \mu_0$ and $\mu_s = (\psi_s)_{\#} \mu_0$ it follows 
	\begin{eqnarray}
	\int \varphi(y) \dd \mu_s^{\eps_k}(y) = \int \varphi(\psi_s^{\eps_k}(x)) \dd \mu_0(x) \ {\rm and}\ 
	\int \varphi(y) \dd \mu_s(y) = \int \varphi(\psi_s(x)) \dd \mu_0(x), \label{push-forwardHeis}
	\end{eqnarray} 
	where all integrals are over $B(0,R)$. 
	Since $\mu_0$ is compactly supported and $\lim_{k\to \infty}\psi_s^{\eps_k}(x) = \psi_s(x)$ for $\mu_0$-a.e. $x$ (cf. Proposition \ref{bridge-proposition}),  the Lebesgue dominated convergence theorem implies
	$\displaystyle
	\int \varphi(\psi_s^{\eps_k}(x)) \dd \mu_0(x) \to \int \varphi(\psi_s(x)) \dd \mu_0(x)
	$ as  $k\to \infty.$
	Combined the latter limit with    (\ref{push-forwardHeis}) the claim follows, i.e.,  $\displaystyle
	\int \varphi(y) \dd \mu_s^{\eps_k}(y) \to \int \varphi(y) \dd \mu_s(y) \mbox{ as } k \to \infty.
	$ In particular, since $\frac{\dd \mu_s^{\eps_k}}{\dd \textsf{m}} = \frac{\dd \mu_s^{\eps_k}}{\dd \textsf{m}_{\eps_k}} \frac{\dd \textsf{m}_{\eps_k}}{\dd \textsf{m}}=\eps_k\rho_s^{\eps_k} \frac{1}{\eps_k}=\rho_s^{\eps_k}$ and $\frac{\dd \mu_s}{\dd \textsf{m}} = \rho_s, $ the latter limit implies  that 
	$
	\displaystyle \int \varphi\left(\rho_s^{\eps_k}-\rho_s\right) \dd \textsf{m}\to 0\ 
	$  as  $ k \to \infty.$
	
	In what follows we need an inequality version of this weak convergence result valid for upper semicontinuous functions. We shall formulate the result as 
	the following:
	
	{\bf Claim.} {\it Let $\varphi: \mathbb H^{n} \to [0,\infty)$ be a bounded, upper semicontinous function. Then the following inequality holds:}
	\begin{equation} \label{limsup}
	\limsup_{{k}\to \infty}\int \varphi(y) \rho^{\eps_{k}}_{s}(y) \dd \textsf{m}(y) \leq \int \varphi(y) \rho_{s}(y) \dd \textsf{m}(y).
	\end{equation} 
	
	To prove the claim let us notice that by the definition of densities and push-forwards of measures the inequality (\ref{limsup}) is equivalent to
	\begin{equation} \label{limsup-1}
	\limsup_{{k}\to \infty}\int \varphi(\psi^{\eps_{k}} _{s}(x)) \dd \mu_{0}(x) \leq  \int \varphi(\psi _{s}(x))   \dd \mu_{0}(x).
	\end{equation} 
	By the upper semicontinuity of $\varphi$ and from the fact that  $\lim_{{k}\to \infty}\psi^{\eps_{k}} _{s}(x) = \psi _{s}(x)$ for every  $x\in D$
	(c.f. Proposition \ref{bridge-proposition}), we obtain
	\begin{equation}\label{limsup-2}
	\limsup_{k\to \infty} \varphi(\psi^{\eps_{k}} _{s}(x)) \leq \varphi (\psi _{s}(x)).
	\end{equation}


	Let $M>0$ be an upper bound of $\varphi$, i.e.,  $0\leq \varphi \leq M$. For an arbitrarily fixed $\delta >0$ we shall prove  that there exists $k_\delta \in \N$ such that for $k\geq k_\delta$,  
	\begin{equation} \label{limsup-3} 
	\int \varphi(\psi^{\eps_{k}} _{s}(x) \dd \mu_{0}(x)) \leq \int \varphi(\psi _{s}(x))  \dd \mu_{0}(x) + (M+1) \delta.
	\end{equation} 
	Since $\delta >0$ is arbitrarily small, the claim (\ref{limsup-1}) would follow from (\ref{limsup-3}). In order to show (\ref{limsup-3}) let us introduce for all $l\in \N$ the set
	$$ S^{l}_{\delta}: = \{ x \in D: \varphi(\psi^{\eps_{k}} _{s}(x)) \leq \varphi (\psi _{s}(x)) + \delta \ \ \text{for all} \ k\geq l \}. $$
	Note that $S^{l}_{\delta} \subseteq S^{l+1}_{\delta}$ for all $l\in \N$ and $\cup_{l} S^{l}_{\delta}= D$; the latter property follows by (\ref{limsup-2}). Since $D$ is a full $\mu_{0}$-measure set it follows that for $\delta >0$ there exists $k_\delta\in \mathbb N$ such that for $k\geq k_\delta$ we have $\mu_{0}(S^{k}_{\delta}) \geq 1 -\delta$. This implies that for every $k\geq k_\delta$ we have the
	estimates
	\begin{eqnarray*}
		\int \varphi(\psi^{\eps_{k}} _{s}(x)) \dd \mu_{0}(x) & \leq & 	\int_{S^{k}_{\delta}} \varphi(\psi^{\eps_{k}} _{s}(x)) \dd \mu_{0}(x) +
		M \mu_{0}(\mathbb H^{n} \setminus S^{k}_{\delta})\\  & \leq & \int_{S^{k}_{\delta}}\varphi(\psi _{s}(x)) \dd \mu_{0}(x)+ 
		\delta \mu_{0}(S^{k}_{\delta}) + M \mu_{0}(\mathbb H^{n} \setminus S^{k}_{\delta}) \\ & \leq & \int \varphi(\psi _{s}(x)) \dd \mu_{0}(x) + (M+1) \delta,
	\end{eqnarray*} 
	concluding the proof of the claim. 
	
	\medskip
	
	We  resume now the proof of the theorem. Since   $\rho_s\in L^1(\dd \textsf{m})$,  there exists a  decreasing sequence of non-negative lower semicontinuous  functions $\{\rho_s^i\}_{i \in \N}$ approximating $\rho_s$ from 
	above. More precisely, we have that $\rho_s^{i}\geq  \rho_s$ and $\rho_s^{i} \to \rho_s$ in $L^1(\dd \textsf{m})$ as $i \to \infty$. Replacing $\rho_s^i$ by $\rho_s^i + \frac{1}{i}$ if necessary, we can even assume that  $\rho_s^i > \rho_s$. In particular, $\rho_s^i$ is  strictly positive   and lower semicontinuous. This implies that  $(\rho_s^i)^{-\frac{1}{2n+1}}$  is positive,  bounded from above and  upper semicontinuous for every $i \in \N$. We introduce the sequence of functions defined by  
	\begin{eqnarray*}
		\rho_s^{\eps_k,i} = \rho^{\eps_k}_s + \rho_s^i - \rho_s,\ i \in \N.
	\end{eqnarray*} 
	
	Note that $\rho_s^{\eps_k,i} >0$ on $D$.  
	To continue the proof of the theorem we notice that the injectivity of the function $\psi_s^{\eps_k}$ on $D$, relation (\ref{MAEq-1}) and  a change of variable  $y = \psi_s^{\eps_k}(x)$ give that
	$$L_s^k=\int_{D} h(\psi_s^{\eps_k}(x)) \left(\rho_s^{\eps_k} (\psi_s^{\eps_k}(x))\right)^{-\frac{1}{2n+1}} {\rho_0(x)} \dd \textsf{m}(x)=\int_{\psi_s^{\eps_k}(D)} h(y) \left(\rho_s^{\eps_k} (y)\right)^{1-\frac{1}{2n+1}}  \dd \textsf{m}(y).$$
	The sub-unitary triangle inequality (i.e., $|a+b|^\alpha \leq |a|^\alpha + |b|^\alpha$ for $a,b \in \R$ and $\alpha\leq 1$), and the convexity of the function $t \mapsto -t^{1-\frac{1}{2n+1}},$ $t>0$ imply the following chain of inequalities 
	{\begin{eqnarray*}
			L_s^k  &=&	 \int_{\psi_s^{\eps_k}(D)} h(y) \left(\rho_s^{\eps_k} (y)\right)^{1-\frac{1}{2n+1}}  \dd \textsf{m}(y)\\
			&\leq& 
			\int_{S} h(y) \left(\rho_s^{\eps_k} (y)\right)^{1-\frac{1}{2n+1}}  \dd \textsf{m}(y)
			\ \ \ \ \ \ \ \ \ \ \ \ \  ({\rm  supp}(h) {=} S)\\
			&\leq&
			\int_{S} h(y)(\rho_s^{\eps_k,i}(y))^{1-\frac{1}{2n+1}} \dd \textsf{m}(y)+ \int_{S} h(y){\left( \rho_s^{\eps_k,i}(y) - \rho^{\eps_k}_s(y) \right)}^{1-\frac{1}{2n+1}} \dd \textsf{m} \\
			&\leq&  \int_{S} h(y)(\rho_s^i(y))^{1-\frac{1}{2n+1}} \dd \textsf{m}(y)+ \frac{2n}{2n+1} \int_{S} h(y)(\rho_s^i(y))^{-\frac{1}{2n+1}}(\rho^{\eps_k}_s(y) - \rho_s(y)) \dd \textsf{m}(y)\\&&  + \int_{S} h(y){\left( \rho_s^i(y)-\rho_s(y) \right)}^{1-\frac{1}{2n+1}} \dd \textsf{m}\\
			&\leq&\int_{S} h(y)(\rho_s(y))^{1-\frac{1}{2n+1}} \dd \textsf{m}(y)  + \frac{2n}{2n+1} \int_{S} h(y) (\rho_s^i(y))^{-\frac{1}{2n+1}}(\rho^{\eps_k}_s(y) - \rho_s(y)) \dd \textsf{m}(y)\\&&+ 2\int_{S} h(y){\left( \rho_s^i(y)-\rho_s(y) \right)}^{1-\frac{1}{2n+1}} \dd \textsf{m}.
		\end{eqnarray*}}
		%
		
		Let $\delta>0$ be arbitrarily fixed. On one hand, by H\"older's inequality and the fact that $\rho_s^{i} \to \rho_s$ in $L^1(\dd \textsf{m})$ as $i \to \infty$, it follows the existence of $i_\delta\in \mathbb N$ such that for every $i\geq i_\delta$, 
		\begin{eqnarray*} 
			\int_{S} h|\rho^{i}_s-\rho_s|^{1-\frac{1}{2n+1}}  \dd \textsf{m} \leq \|h\|_{L^\infty(S)}\left(\int_{S} |\rho^i_s-\rho_s| \dd \textsf{m} \right)^{1-\frac{1}{2n+1}} \left(\textsf{m}(\R^{2n+1})\right)^{\frac{1}{2n+1}}<\frac{\delta}{4}. 
		\end{eqnarray*}
		On the other hand, since $y\mapsto \varphi(y)=h(y)(\rho_s^{i_\delta}(y))^{-\frac{1}{2n+1}}$ is positive, bounded from above and upper semicontinuous, by (\ref{limsup})  we find $k_\delta \in \N$ such that 
		\begin{eqnarray*}\label{Conv2}
			\frac{2n}{2n+1} \int_{S} h(y)(\rho_s^{i_\delta}(y))^{-\frac{1}{2n+1}}\left(\rho^{\eps_{k}}_s(y) - \rho_s(y)\right) \dd \textsf{m}(y) < \frac{\delta}{2} \quad \mbox{ for all } k \geq k_\delta.
		\end{eqnarray*}
		{Summing up the above estimates, for every $k \geq k_\delta$ we have 
			$$L_s^k \leq \int_{S} h(y)\left(\rho_s(y)\right)^{1-\frac{1}{2n+1}} \dd \textsf{m}(y) + \delta.$$
			Thus, the arbitrariness of $\delta>0$ implies that 
			\begin{equation}\label{ujabb-limit}
			\liminf\limits_{k \to \infty} L_s^k\leq \int_{S} h(y)\left(\rho_s(y)\right)^{1-\frac{1}{2n+1}} \dd \textsf{m}(y)=\int_{S\cap {\rm supp}(\rho_s)} h(y)\left(\rho_s(y)\right)^{1-\frac{1}{2n+1}} \dd \textsf{m}(y).
			\end{equation}
			Since ${\rm supp}(\rho_s)\subseteq \psi_s(D)$, by (\ref{ujabb-limit}) we have that 
			$$\liminf\limits_{k \to \infty} L_s^k\leq \int_{\psi_s(D)} h(y)\left(\rho_s(y)\right)^{1-\frac{1}{2n+1}} \dd \textsf{m}(y).$$
			Now, the injectivity of the map $D\ni x\mapsto \psi_s(x)$, a  change of variable  $y = \psi_s(x)$ in the right hand side of the latter estimate,  and the 
			Monge-Amp\`ere equation  
			\begin{equation}\label{Monge-Ampere}
			\rho_0(x)=\rho_s(\psi_s(x)){\rm Jac}(\psi_s)(x),\ \
			x\in D,
			\end{equation}
			give   the inequality in 
			(\ref{entropiak-weak-1}). }
		
	\noindent 	Combining the estimates (\ref{elo-entropia}) and (\ref{entropiak-weak-1}), we obtain 
		\begin{eqnarray*}
			\int_D h(\psi_s(x))\left(\rho_s(\psi_s(x))\right)^{-\frac{1}{2n+1}}\rho_0(x)\dd \textsf{m}(x)&\geq & 
			\int_D h(\psi_s(x))\tau_{1-s}^n(\theta_x) \left(\rho_0(x)\right)^{1-\frac{1}{2n+1}} \dd \textsf{m}(x) \nonumber \\
			&& + \int_{D}h(\psi_s(x))\tau_{s}^n(\theta_x) \left(\rho_1(\psi(x))\right)^{-\frac{1}{2n+1}} {\rho_0(x)} \dd \textsf{m}(x).
		\end{eqnarray*}
		Applying the change of variables $y= \psi_{s}(x)$ and (\ref{Monge-Ampere}) we obtain 
		\begin{eqnarray*}
			\int_{\psi_{s}(D)} h(y)\left(\rho_s(y)\right)^{1-\frac{1}{2n+1}} \dd \textsf{m}(y)&\geq & 
			\int_{\psi_{s}(D)} h(y)\tau_{1-s}^n(\theta_{\psi_{s}^{-1}(y)}) \left(\rho_0(\psi_{s}^{-1}(y))\right)^{-\frac{1}{2n+1}} \rho_{s}(y)\dd \textsf{m}(y) \nonumber \\
			&& + \int_{\psi_{s}(D)}h(y)\tau_{s}^n(\theta_{\psi_{s}^{-1}(y)}) \left((\rho_1\circ \psi)({\psi_{s}^{-1}(y)})\right)^{-\frac{1}{2n+1}} \rho_{s}(y) \dd \textsf{m}(y).
		\end{eqnarray*}
		
		Observe that the function on the left side of the above estimate that multiplies  $h$ is $\rho_{s }^{1-\frac{1}{2n+1}}$ which is in $L^{1}(\dd \textsf{m})$. Since we are considering only 
		positive functions it follows that the function on the right side multiplying $h$ is also in $L^{1}(\dd \textsf{m})$. We shall use the well-known fact that convolutions with mollifiers 
		converge point-wise almost everywhere to the function values for functions in $L^{1}(\dd \textsf{m})$. 
		
		Since the test function  $h\geq 0$ is arbitrary,  it can play the role of convolution 
		kernels. From here we can conclude that the latter integral inequality  implies the point-wise inequality: 
		\begin{eqnarray*}
			\left(\rho_s(y)\right)^{-\frac{1}{2n+1}} 
			\geq \tau_{1-s}^n(\theta_{\psi_{s}^{-1}(y)}) (\rho_0(\psi_{s}^{-1}(y)))^{-\frac{1}{2n+1}}   + \tau_{s}^n(\theta_{\psi_{s}^{-1}(y)}) \left((\rho_1\circ \psi)(\psi_{s}^{-1}(y))\right)^{-\frac{1}{2n+1}}
		\end{eqnarray*}
		for a.e. $y\in \psi_{s}(D)$.
		Composing with $\psi_{s}$ the above estimate, it yields 
		\begin{eqnarray}\label{surusegek}
		\left(\rho_s(\psi_s(x))\right)^{-\frac{1}{2n+1}}&\geq & 
		\tau_{1-s}^n(\theta_x) (\rho_0(x))^{-\frac{1}{2n+1}}   + \tau_{s}^n(\theta_x) \left(\rho_1(\psi(x))\right)^{-\frac{1}{2n+1}} \ \ {\rm  a.e.} \ x\in D.
		\end{eqnarray}
		By  the 
		Monge-Amp\`ere equations  (\ref{Monge-Ampere}) and $\rho_0(x)=\rho_1(\psi(x)){\rm Jac}(\psi)(x),$ $x\in D$, we obtain the inequality 
		$$\left({\rm Jac}(\psi_s)(x)\right)^\frac{1}{2n+1}\geq \tau_{1-s}^n(\theta_x)+\tau_{s}^n(\theta_x)\left({\rm Jac}(\psi)(x)\right)^\frac{1}{2n+1} \ \ {\rm  a.e.} \ x\in D,$$
		which concludes the proof.  \hfill $\square$
		
		\begin{remark} \rm \label{remark-mikor-nem-feltetlenul-absz-folyt} Observe that the Jacobian identity on the Riemannian manifolds $M^{\eps_k}$ (cf. \cite[Lemma 6.1]{McCann}) that constitutes the starting point of the proof of our determinant inequality holds also in the case when $\mu_1$ is not necessarily absolutely continuous w.r.t. the $\mathcal L^{2n+1}$-measure. In this case our arguments are based on the inequality that we obtain by canceling the second term of the right side, namely
			\begin{eqnarray*}
				\left({\rm
					Jac}(\psi^{\eps_k}_s)(x)\right)^\frac{1}{2n+1} &\geq&
				(1-s)\left(v^{\eps_k}_{1-s}(\psi^{\eps_k}(x),x)\right)^\frac{1}{2n+1}.
			\end{eqnarray*}
			Now we can perform the same steps as in the proof of Theorem \ref{TJacobianDetIneq}
			by obtaining
			\begin{eqnarray}\label{surusegek-nemabszfolyt-esetben}
			\left(\rho_s(\psi_s(x))\right)^{-\frac{1}{2n+1}}&\geq & 
			\tau_{1-s}^n(\theta_x) (\rho_0(x))^{-\frac{1}{2n+1}} \ \ {\rm  a.e.} \ x\in D,
			\end{eqnarray}
			or equivalently 
			\begin{equation}\label{Jacobi-degenerate}
			\left({\rm Jac}(\psi_s)(x)\right)\geq \tau_{1-s}^n(\theta_x)^{2n+1}\ \ {\rm  a.e.} \ x\in D.
			\end{equation}
		\end{remark}
		
		\noindent A direct consequence of (\ref{surusegek-nemabszfolyt-esetben}) and (\ref{tau-becsles}) is the main estimate from the paper of Figalli and Juillet \cite{FJ} formulated and refined in the following statement:

		\begin{corollary}{\bf (Interpolant density estimate  on $\mathbb H^n$)}\label{Corollary-Figalli-Juillet} Under the same assumptions as in Theorem \ref{TJacobianDetIneq} $($except  the absolutely continuous property of $\mu_1$$)$,  we have
			$$\rho_s(y) \leq  
			\left(\tau_{1-s}^n\left(\theta_{\psi_s^{-1}(y)}\right)\right)^{-(2n+1)} \rho_0(\psi_s^{-1}(y)) \leq \frac{1}{(1-s)^{2n+3}}\rho_0(\psi_s^{-1}(y))\mbox{ for } \mu_s \mbox{-a.e. } y \in \mathbb H^n.$$
		\end{corollary}
		
		\medskip \medskip
		
		\begin{remark} \rm \label{tau-improvement} 
			A closer inspection of inequality (\ref{Jacobi-inequality-elso}) from Theorem \ref{TJacobianDetIneq} shows that it can be improved in the presence of a positive measure set of static points. Indeed, if  $x$ is a static point of  $\psi$ than it follows that it will be a static point for $\psi_{s}(x)$ as well. Considering density points of the static set, (i.e. discarding a null set if necessary) we obtain that both Jacobians  ${\rm Jac}(\psi_s)(x)= {\rm Jac}(\psi)(x)= 1$ on a full measure of stationary points. 
			This implies that  relation  (\ref{Jacobi-inequality-elso}) holds with $\tau_s^n(\theta_x)=s$ and $\tau_{1-s}^n(\theta_x)=1-s$.
			
			Based on this observation it is natural to  define a new, optimal transport based Heisenberg distortion coefficient $\hat{\tau}_{s,\psi}^n$  which  depends directly on $x \in \Heis^n$  rather than on the angle $\theta_x$. If $s \in (0,1)$, we consider  
			\begin{eqnarray}\label{better-concentration}
			\hat{\tau}_{s,\psi}^n(x) = \left\{
			\begin{array}{lll}
			\tau_s^n(\theta_x) & \mbox{if} & x\in \mathcal M_\psi; \\
			s &\mbox{if} &  x\in \mathcal S_\psi.
			\end{array}\right.
			\end{eqnarray}
			With this notation,  under the assumptions of Theorem \ref{TJacobianDetIneq}, the following improved version of (\ref{Jacobi-inequality-elso}) holds: 
			\begin{equation}\label{Jacobi-inequality-improvement}
			\left({\rm Jac}(\psi_s)(x)\right)^\frac{1}{2n+1}\geq \hat{\tau}_{1-s, \psi}^n(\theta_x)+\hat{\tau}_{s, \psi}^n(\theta_x)\left({\rm Jac}(\psi)(x)\right)^\frac{1}{2n+1}\mbox{ for } \mu_0 \mbox{-a.e. } x \in \mathbb H^n.
			\end{equation}
		\end{remark}
		
		\medskip
		
		\subsection{Entropy inequalities on $\mathbb H^n$} As a first application of the Jacobian determinant inequality we prove several entropy inequalities on $\mathbb H^n$. 
		
		\medskip
		
		\noindent {\bf Proof of Theorem \ref{TEntIneqHeisGen}.}  
		We shall keep the notations from the proof of Theorem \ref{TJacobianDetIneq}. Since the function $t \mapsto t^{2n+1}U(t^{-(2n+1)})$ is non-increasing, relation   (\ref{surusegek}) implies that  for a.e. $x\in D$ we have 
		\begin{eqnarray*}
			\frac{U\left(\rho_s(\psi_s(x))\right)}{\rho_s(\psi_s(x))} \leq &&\left(\tau_{1-s}^n(\theta_x) (\rho_0(x))^{-\frac{1}{2n+1}}   + \tau_{s}^n(\theta_x) (\rho_1(\psi(x)))^{-\frac{1}{2n+1}}\right)^{(2n+1)} \times\\
			&& \times U\left(\left(\tau_{1-s}^n(\theta_x) (\rho_0(x))^{-\frac{1}{2n+1}}   + \tau_{s}^n(\theta_x) (\rho_1(\psi(x)))^{-\frac{1}{2n+1}}\right)^{-(2n+1)}\right).
		\end{eqnarray*}
		Recalling relation $s \tilde{\tau}_s^n = \tau_s^n$, the right hand side of the above inequality  can be written as  
		\begin{eqnarray*}
			\left((1-s) \left(\frac{\rho_0(x)}{\left(\tilde{\tau}_{1-s}^n(\theta_x)\right)^{2n+1}}\right)^{-\frac{1}{2n+1}}   + s\left(\frac{\rho_1(\psi(x))}{\left(\tilde{\tau}_{s}^n(\theta_x)\right)^{2n+1}}\right)^{-\frac{1}{2n+1}}\right)^{(2n+1)} \times\\
			\times U\left(\left((1-s) \left(\frac{\rho_0(x)}{\left(\tilde{\tau}_{1-s}^n(\theta_x)\right)^{2n+1}}\right)^{-\frac{1}{2n+1}}   + s\left(\frac{\rho_1(\psi(x))}{\left(\tilde{\tau}_{s}^n(\theta_x)\right)^{2n+1}}\right)^{-\frac{1}{2n+1}}\right)^{-(2n+1)}\right).
		\end{eqnarray*}
		By using the convexity of  $t \mapsto t^{2n+1}U(t^{-(2n+1)})$, the latter term can be estimated from above by
		$$(1-s)\frac{\left(\tilde{\tau}_{1-s}^n(\theta_x)\right)^{2n+1}}{\rho_0(x)}U\left(\frac{\rho_0(x)}{\left(\tilde{\tau}_{1-s}^n(\theta_x)\right)^{2n+1}}\right)
		+ s\frac{\left(\tilde{\tau}_s^n(\theta_x)\right)^{2n+1}}{\rho_1(\psi(x))}U\left(\frac{\rho_1(\psi(x))}{\left(\tilde{\tau}_s^n(\theta_x)\right)^{2n+1}}\right).$$
		Summing up, for a.e. $x \in D$ we have
		$$ \frac{U\left(\rho_s(\psi_s(x))\right)}{\rho_s(\psi_s(x))} \leq (1-s)\frac{\left(\tilde{\tau}_{1-s}^n(\theta_x)\right)^{2n+1}}{\rho_0(x)}U\left(\frac{\rho_0(x)}{\left(\tilde{\tau}_{1-s}^n(\theta_x)\right)^{2n+1}}\right)
		+ s\frac{\left(\tilde{\tau}_s^n(\theta_x)\right)^{2n+1}}{\rho_1(\psi(x))}U\left(\frac{\rho_1(\psi(x))}{\left(\tilde{\tau}_s^n(\theta_x)\right)^{2n+1}}\right).$$
		An integration of the above inequality on $D$  w.r.t. the measure $\mu_0=\rho_0 \textsf{m}$  gives
		\begin{eqnarray*} 
			\int_{D}U\left(\rho_s(\psi_s(x))\right)\frac{\rho_0(x)}{\rho_s(\psi_s(x))} \dd \textsf{m}(x) \leq (1-s)\int_{D}{\left(\tilde{\tau}_{1-s}^n(\theta_x)\right)^{2n+1}}U\left(\frac{\rho_0(x)}{\left(\tilde{\tau}_{1-s}^n(\theta_x)\right)^{2n+1}}\right) \dd \textsf{m}(x)\\
			\ \ \ \ \ \ \ + s\int_{D}\left(\tilde{\tau}_s^n(\theta_x)\right)^{2n+1}U\left(\frac{\rho_1(\psi(x))}{\left(\tilde{\tau}_s^n(\theta_x)\right)^{2n+1}}\right)\frac{\rho_0(x)}{\rho_1(\psi(x))} \dd \textsf{m}(x).
		\end{eqnarray*}
		Recall that $\psi_s$ and $\psi$ are injective on $D$; thus we can perform the changes of variables $z = \psi_s(x)$ and $y = \psi(x)$ and by the Monge-Amp\`ere equations $\rho_0(x)=\rho_s(\psi_s(x)){\rm Jac}(\psi_s)(x)$  and $\rho_0(x)=\rho_1(\psi(x)){\rm Jac}(\psi)(x),$ $x\in D$, we obtain the required entropy inequality. 
		\hfill $\square$
		
		\medskip
		
		By (\ref{kapcsolat-volume-es-CD}) and the monotonicity of $t \mapsto t^{2n+1}U(t^{-(2n+1)})$ we obtain a sub-Riemannian displacement convexity property of the entropy:

		\begin{corollary}\label{Corollary-TEntIneqHeisGen} {\bf (Uniform entropy inequality on $\Heis^n$)}
			Under the same assumptions as in Theorem \ref{TJacobianDetIneq},  the following entropy inequality holds:
			\begin{eqnarray*}
				{\rm Ent}_{U}(\mu_s | \mathcal L^{2n+1})   &\leq&  (1-s)^3 \int_{\mathbb H^n}  U\left(\frac{\rho_0(x)}{(1-s)^2}\right) \dd \mathcal L^{2n+1}(x) + s^3 \int_{\mathbb H^n}  U\left(\frac{\rho_1(y)}{s^2}\right) \dd \mathcal L^{2n+1}(y).
			\end{eqnarray*}
		\end{corollary}
		\medskip
		\noindent Some relevant admissible functions $U: [0, \infty) \to \mathbb R$ in Theorem \ref{TEntIneqHeisGen} are as follows: 
		\begin{enumerate}
			\item[$\bullet$]  {\it R\'enyi-type entropy}: $U_R(t) = -t^\gamma$ with $\gamma\in [{1 - \frac{1}{2n+1}},1]$; for $\gamma={1 - \frac{1}{2n+1}}$ we have precisely the R\'enyi entropy Ent$_{U_R}$=Ent$_{2n+1}$ from (\ref{EntDef}).  
			\item[$\bullet$] {\it Shannon entropy}: $U_S(t) = t \log t$ for $t>0$ and $U_S(0)=0$.
			\item[$\bullet$] {\it Kinetic-type entropy}: $U_K(t) = t^\gamma$ with $\gamma\geq 1.$ 
			\item[$\bullet$] {\it Tsallis entropy}: $U_T(t) = \frac{t^\gamma-t}{\gamma-1}$ with $\gamma\in [{1 - \frac{1}{2n+1}},\infty)\setminus \{1\};$  the limiting case $\gamma\to 1$  reduces to the Shannon entropy. 
		\end{enumerate}

		\noindent 	Applying Theorem \ref{TEntIneqHeisGen} together with (\ref{tau-becsles})    to  $U_R(t)=-t^{{1 - \frac{1}{2n+1}}},$ one has
		
		\begin{corollary}\label{TEntIneqHeis-corollary} {\bf (R\'enyi entropy inequality on $\Heis^n$)}
			Under the same assumptions as in Theorem \ref{TJacobianDetIneq},  the following entropy inequality holds:
			\begin{eqnarray*}
				{\rm Ent}_{2n+1}(\mu_s|\mathcal L^{2n+1}) & \leq & - \int_{\mathbb H^n} \tau_{1-s}^n(\theta_x) \rho_0(x)^{1 - \frac{1}{2n+1}} \dd \Lmeas^{2n+1}(x)  \\&&\ \ \ \ \ \  \ \  \ - \int_{\mathbb H^n} \tau_{s}^n(\theta_{\psi^{-1}(y)}) \rho_1(y)^{1 - \frac{1}{2n+1}} \dd \Lmeas^{2n+1}(y)\\
				&\leq &(1-s)^{\frac{2n+3}{2n+1}}\ {\rm Ent}_{2n+1}(\mu_0|\mathcal L^{2n+1}) + s^{\frac{2n+3}{2n+1}}\ {\rm Ent}_{2n+1}(\mu_1|\mathcal L^{2n+1}).
			\end{eqnarray*}
		\end{corollary}
		\medskip

		
		Let $U_S(t) = t \log t$ for $t>0$ and $U_S(0) =0$; Corollary \ref{Corollary-TEntIneqHeisGen} implies the following  convexity-type property of the Shannon entropy $s\mapsto {\rm Ent}_{U_S}(\mu_s|{\sf{m}})$ on $\mathbb H^n$:
		
		\begin{corollary}\label{Shannon-TEntIneqHeis-corollary} {\bf (Uniform Shannon entropy inequality on $\Heis^n$)}
			Under the same assumptions as in Theorem \ref{TJacobianDetIneq},  the following  entropy inequality holds:
			\begin{eqnarray*}
				{\rm Ent}_{U_S}(\mu_s|\mathcal L^{2n+1}) 
				&\leq &(1-s)\ {\rm Ent}_{U_S}(\mu_0|\mathcal L^{2n+1}) + s\ {\rm Ent}_{U_S}(\mu_1|\mathcal L^{2n+1})-2\log ((1-s)^{1-s}s^s).
			\end{eqnarray*}
		\end{corollary}

		\begin{remark}\rm 
			The positive concave function  $w(s)=-2\log ((1-s)^{1-s}s^s)$ compensates  the lack of convexity of $s\mapsto {\rm Ent}_{U_S}(\mu_s|{\sf{m}})$, $s\in (0,1)$. Notice also that we have  $0< w(s)\leq \log 4=w\left(\frac{1}{2}\right)$ for every $s\in (0,1)$, and $\lim_{s\to 0}w(s)=\lim_{s\to 1}w(s)=0.$  
		\end{remark}
		
		\begin{remark} \rm \label{tau-improvement2} 
			
			Based on  Remark \ref{tau-improvement} we can also define optimal transport based coefficients $\hat{\tilde{\tau}}_{s,\psi}^n$ as 
			\begin{eqnarray}\label{better-concentration2}
			\hat{\tilde{\tau}}_{s,\psi}^n(x) = \left\{
			\begin{array}{lll}
			\tilde{\tau}_s^n(\theta_x) & \mbox{if} & x\in \mathcal M_\psi; \\
			1 &\mbox{if} &  x\in \mathcal S_\psi,
			\end{array}\right.
			\end{eqnarray}
			and state a corresponding version of Theorem \ref{TEntIneqHeisGen} with respect to these coefficients.
		\end{remark}
		
		\medskip 
		\subsection{Borell-Brascamp-Lieb and Pr\'ekopa-Leindler inequalities on $\mathbb H^n$} In this subsection we prove various Borell-Brascamp-Lieb and Pr\'ekopa-Leindler inequalities on $\mathbb H^n$ by showing another powerful  application of the Jacobian determinant inequality. 
		
		\medskip 
		\noindent {\bf Proof of Theorem \ref{TRescaledBBLWithWeights}.} We assume that
		hypotheses of Theorem \ref{TRescaledBBLWithWeights} are fulfilled.
		Let $s\in (0,1)$ and $p\geq -\frac{1}{2n+1}$.  Note that if either $\displaystyle\int_{\Heis^n} f=0$
		or $\displaystyle\int_{\Heis^n} g=0$, the conclusion follows due to our
		convention concerning the $p$-mean $M_s^p$. Thus, we may assume that
		both integrals are positive. The proof is divided into three parts. 
		
		{\it Step 1.} {We first consider the particular case when the functions $f, g$ are compactly supported and normalized, i.e.,
			\begin{equation}\label{normalized}
			\displaystyle\int_{\Heis^n} f=  \displaystyle\int_{\Heis^n} g= 1.
			\end{equation}}
		Let us keep the notations from the proof of Theorem \ref{TJacobianDetIneq}, by identifying the density functions $\rho_0$ and $\rho_1$ of the measures $\mu_0$ and $\mu_1$ with $f$ and $g$, respectively. Since the Jacobian determinant inequality is equivalent to (\ref{Jacobi-inequality-elso-ekvivalens}), we have that 
		\begin{eqnarray}\label{IneqJacDet}
		\rho_s(\psi_s(x))^{-\frac{1}{2n+1}} \geq  
		\tau_{1-s}^n(\theta_x) (f(x))^{-\frac{1}{2n+1}}   + \tau_{s}^n(\theta_x) (g(\psi(x)))^{-\frac{1}{2n+1}} \ \ \mbox{ for }  \mbox{ a.e.} \ x\in D.
		\end{eqnarray}
		Choosing $y= \psi(x) $ in hypothesis (\ref{ConditionRescaledBBLWithWeights}) for points $x \in D$ we obtain:
		\begin{eqnarray} \label{IneqAssumptionBBL}
		h(\psi_s(x)) \geq M^{p}_s
		\left(\frac{f(x)}{\left(\tilde{\tau}_{1-s}^n(\theta_x)\right)^{2n+1}},\frac{g(\psi(x))}{\left(\tilde{\tau}_s^n(\theta_x)\right)^{2n+1}} \right). 
		\end{eqnarray}
		Since $p\geq -\frac{1}{2n+1}$, the monotonicity of the $p$-mean, relation $\tilde{\tau}_s^n = s^{-1}\tau_s^n$ and inequalities  (\ref{IneqJacDet}),  (\ref{IneqAssumptionBBL}) imply  
		that 
		$$h(\psi_s(x)) \geq  M^{-\frac{1}{2n+1}}_s
		\left(\frac{f(x)}{\left(\tilde{\tau}_{1-s}^n(\theta_x)\right)^{2n+1}},\frac{g(\psi(x))}{\left(\tilde{\tau}_s^n(\theta_x)\right)^{2n+1}} \right)\geq\rho_s(\psi_s(x)) \mbox{ for  a.e. } x \in D.$$
		Since $\mu_s=(\psi_s)_\#\mu_0$, an integration and change of variables give that
		\begin{eqnarray*}
			\int_{\mathbb H^n} h&\geq& \int_{\psi_s(D)} h(z)\dd\mathcal L^{2n+1}(z)=\int_{D} h(\psi_s(x)){\rm Jac}(\psi_s)(x)\dd\mathcal L^{2n+1}(x)\\ &\geq& \int_{D} \rho_s(\psi_s(x)){\rm Jac}(\psi_s)(x)\dd\mathcal L^{2n+1}(x)=\int_D \rho_0(y)\dd\mathcal L^{2n+1}(y)=\int_D f(y)\dd\mathcal L^{2n+1}(y)\\&=&1.
		\end{eqnarray*}

		{\it Step 2.} We assume that the functions $f, g$ are compactly supported and $\displaystyle 0<\int_{\mathbb H^n}f<\infty$ and $\displaystyle 0<\int_{\mathbb H^n}g<\infty.$
		To proceed further, we first recall the inequality for $p$- and $q$-means from Gardner \cite[Lemma
		10.1]{Gardner}, i.e., 
		\begin{eqnarray}\label{MspIneq}
		M_s^{p}(a,b)M_s^{q}(c,d) \geq M_s^{\eta}(ac, bd),
		\end{eqnarray} for every $a,b,c,d \geq 0, s \in (0,1)$ and
		$p, q \in \R$ such that $p+q\geq0$ with ${\eta}=\frac{pq}{p+q}$ when
		$p$ and $q$ are not both zero, and $\eta=0$ if $p=q=0$.

		Define $\tilde{f} =\displaystyle \frac{f}{\displaystyle\int_{\mathbb H^n}f}$, $\tilde{g} =
		\displaystyle\frac{g}{\displaystyle\int_{\mathbb H^n}g}$ and $\tilde{h} =
		\left(M_s^{\frac{p}{1+(2n+1)p}}\left(\displaystyle\int_{\mathbb H^n}f, \displaystyle\int_{\mathbb H^n}g\right)\right)^{-1} h$.
		Clearly,  we have the relations $\displaystyle \int_{\mathbb H^n}{\tilde f} =\int_{\mathbb H^n}  \tilde{g} = 1$.

		We apply inequality (\ref{MspIneq}) with the choice of $q= \frac{-p}{1 +(2n+1)p}$ and $p\geq \frac{-1}{2n+1}$. Notice that $p+q\geq 0$ is satisfied and we have 
		that $\eta= -\frac{1}{2n+1}$. By hypothesis (\ref{ConditionRescaledBBLWithWeights}) we have that for every $x,y\in \mathbb H^n$ and $z\in Z_s(x,y)$, 
		\begin{eqnarray*}
			\tilde{h}(z) &=& \left(M_s^{\frac{p}{1+(2n+1)p}}\left(\int_{\mathbb H^n}f, \int_{\mathbb H^n}g\right)\right)^{-1} h(z) = M_s^{-\frac{p}{1+(2n+1)p}}\left(\frac{1}{\displaystyle\int_{\mathbb H^n}f}, \frac{1}{\displaystyle\int_{\mathbb H^n}g}\right) h(z) \\
			&\geq& M_s^{-\frac{p}{1+(2n+1)p}}\left(\frac{1}{\displaystyle\int_{\mathbb H^n}f}, \frac{1}{\displaystyle\int_{\mathbb H^n}g}\right) M^{p}_s \left(\frac{f(x)}{\left(\tilde{\tau}_{1-s}^n(\theta(y,x))\right)^{2n+1}},\frac{g(y)}{\left(\tilde{\tau}_{s}^n(\theta(x,y))\right)^{2n+1}} \right)  \\
			&\geq&M_s^{-\frac{1}{2n+1}
			}\left(\frac{\tilde{f}(x)}{\left(\tilde{\tau}_{1-s}^n(\theta(y,x))\right)^{2n+1}},\frac{\tilde{g}(y)}{\left(\tilde{\tau}_{s}^n(\theta(x,y))\right)^{2n+1}}
			\right).
		\end{eqnarray*}
		Now we are in the position to apply Step 1
		for the functions $\tilde{f}, \tilde{g}$ and $\tilde{h}$, obtaining
		that
		$\displaystyle\int_{\mathbb H^n}\tilde h\geq 1,$ which is equivalent to
		\begin{eqnarray*}
			\int_{\mathbb H^n} h \geq
			M_s^{\frac{p}{1+(2n+1)p}}\left(\int_{\mathbb H^n}f, \int_{\mathbb H^n}g\right).
		\end{eqnarray*}

		{\it Step 3.} {We now consider the general case when $f$ and $g$ are not necessarily compactly supported. The integrable functions $f$ and $g$ can be approximated in $L^1(\mathbb H^n)$ from below by upper semicontinuous compactly supported functions; let $\{f_k\}_{k \in \mathbb N}$ and $\{g_k\}_{k \in \mathbb N}$ be these approximating function sequences. 
			We observe that  hypothesis (\ref{ConditionRescaledBBLWithWeights}) is inherited by the  triplet $\{h,f_k, g_k\}_{k \in \mathbb N}$  via the monotonicity of $M_s^p(\cdot,\cdot)$, i.e.,
			\begin{eqnarray*}
				h(z) \geq M^{p}_s
				\left(\frac{f_k(x)}{\left(\tilde{\tau}_{1-s}^n(\theta(y,x))\right)^{2n+1}},\frac{g_k(y)}{\left(\tilde{\tau}_s^n(\theta(x,y))\right)^{2n+1}} \right) 
				\  {\rm for\ all}\ (x,y)\in \Heis^n\times \Heis^n, z\in Z_s(x,y).
			\end{eqnarray*}
			By applying Step 2 for every $k\in \mathbb N$,  it yields that
			$$\int_{\mathbb H^n} h  \geq
			M_s^{\frac{p}{1+(2n+1)p}}\left(\int_{\mathbb H^n}f_k, \int_{\mathbb H^n}g_k\right).$$
			Letting $k \to \infty,$ we conclude the proof.
			\hfill $\square$
			
			\begin{remark}\rm  If $\displaystyle\int_{\Heis^n} f=+\infty$ or $\displaystyle\int_{\mathbb H^n} g=+\infty$
				we can apply a standard approximation argument, similar to Step 3 from the previous proof, obtaining that $\displaystyle\int_{\Heis^n} h=+\infty$.
			\end{remark}
			\medskip
			

			\begin{corollary}\label{CLighterWeightedBBL}  {\bf (Uniformly weighted Borell-Brascamp-Lieb inequality on $\mathbb H^n$)}
				Fix $s\in (0,1)$ and $p \geq -\frac{1}{2n+1}.$  Let $f,g,h:\mathbb
				H^n\to [0,\infty)$ be  Lebesgue integrable
				functions satisfying
				\begin{eqnarray}&\label{1-ConditionRescaledBBLWithoutWeights}
				h(z) \geq M^{p}_s \left(\frac{f(x)}{(1-s)^2},\frac{g(y)}{s^2}\right) \ \  for\ all\ (x,y)\in \Heis^n\times
				\Heis^n, z\in Z_s(x,y).
				\end{eqnarray}
				Then the following inequality holds:
				\begin{eqnarray*}
					\int_{\Heis^n} h \geq M^\frac{p}{1+(2n+1)p}_s \left(
					\int_{\Heis^n} f, \int_{\Heis^n}g\right).
				\end{eqnarray*}
			\end{corollary}
			
			{\it Proof.} Directly follows by Theorem \ref{TRescaledBBLWithWeights} and relation (\ref{kapcsolat-volume-es-CD}).  \hfill $\square$ 
			
			\begin{corollary}\label{CRescaledBBLWithoutWeights}  {\bf (Non-weighted Borell-Brascamp-Lieb inequality on $\mathbb H^n$)}
				Fix $s\in (0,1)$ and $p \geq -\frac{1}{2n+3}.$  Let $f,g,h:\mathbb
				H^n\to [0,\infty)$ be  Lebesgue integrable
				functions satisfying
				\begin{eqnarray}&\label{ConditionRescaledBBLWithoutWeights}
				h(z) \geq M^{p}_s (f(x),g(y)) \ \  for\ all\ (x,y)\in \Heis^n\times
				\Heis^n, z\in Z_s(x,y).
				\end{eqnarray}
				Then the following inequality holds:
				\begin{eqnarray}\label{correction}
				\int_{\Heis^n} h \geq \frac{1}{4}M^\frac{p}{1+(2n+3)p}_s \left(
				\int_{\Heis^n} f, \int_{\Heis^n}g\right).
				\end{eqnarray}
			\end{corollary}

			{\it Proof.}  Let us
			first assume that $\displaystyle\int_{\mathbb H^n}f= \int_{\mathbb H^n}g= 1.$ By using the 
			notations from Theorems \ref{TJacobianDetIneq} \& \ref{TRescaledBBLWithWeights},  we explore hypothesis
			(\ref{ConditionRescaledBBLWithoutWeights}) only for the pairs
			$(x,\psi(x))\in A_0\times A_1$ with $x\in D$. In particular,
			$x^{-1}\cdot \psi(x)\notin L^*$ for every $x\in D.$
			
			Fix  $s\in (0,1)$ and $p \geq -\frac{1}{2n+3}$. By  the
			$p$-mean inequality (\ref{MspIneq}), we have for every  $x\in D$ that
			$$ M_s^p(f(x),g(\psi(x))) \geq 
			M_s^{-\frac{1}{2}} (v^0_{1-s}(\psi(x),x), v^0_s(x,\psi(x)))
			M_s^{\frac{p}{2p+1}} \left( \frac{f(x)}{v^0_{1-s}(\psi(x),x)},
			\frac{g(\psi(x))}{v^0_s(x,\psi(x))} \right).
			$$
			According to 
			(\ref{kapcsolat-volume-es-CD}),   for every $ x\in D$ we have  
			$$M_s^{-\frac{1}{2}} (v^0_{1-s}(\psi(x),x), v^0_s(x,\psi(x))) \geq M_s^{-\frac{1}{2}} ((1-s)^2, s^2) = \frac{1}{4}.$$
			Now, by hypothesis (\ref{ConditionRescaledBBLWithoutWeights}) and relation (\ref{kapcsolat-volume-es-CD}) we
			have for every $x\in D$ and $z= Z_s(x,\psi(x))$ that
			$$h(z)\geq\frac{1}{4} M_s^{\frac{p}{2p+1}} \left( \frac{f(x)}{\left(\tilde{\tau}_{1-s}^n(\theta_x)\right)^{2n+1}}, \frac{g(\psi(x))}{\left(\tilde{\tau}_{s}^n(\theta_x)\right)^{2n+1}} \right). 
			$$ By the assumption $p\geq \frac{-1}{2n+3}$ we have $\frac{p}{2p+1}\geq-\frac{1}{2n+1}$. A similar argument as in the proof of 
			Theorem \ref{TRescaledBBLWithWeights} (see relation (\ref{IneqAssumptionBBL})) yields that $$\displaystyle\int_{\mathbb H^n}h\geq
			\frac{1}{4}.$$
			
			The general case follows again as in Theorem
			\ref{TRescaledBBLWithWeights}, replacing the power $p$ by
			$\frac{p}{2p+1}$; therefore, 
			$$\int_{\mathbb H^n}h \geq
			\frac{1}{4}M_s^{\frac{p}{1+(2n+3)p}}\left(\int_{\mathbb H^n}f, \int_{\mathbb H^n}g\right),$$
			which concludes the proof.
			\hfill $\square$
			
			\medskip
			
			\begin{remark}\rm 
				We notice that  we
				pay a price in Corollary \ref{CRescaledBBLWithoutWeights} for missing out  of Heisenberg volume distortion
				coefficients $\tilde{\tau}_{1-s}^n$ and $\tilde{\tau}_{s}^n$ from (\ref{ConditionRescaledBBLWithWeights}) or the weights $(1-s)^2$ and $s^2$  from (\ref{1-ConditionRescaledBBLWithoutWeights}), respectively. Indeed, unlike in the
				Euclidean case (where the volume distortions are identically 1), we
				obtain $\frac{1}{4}$ as a correction factor in the right hand side
				of (\ref{correction}). Note however that the constant $\frac{1}{4}$
				is {\it sharp} in (\ref{correction}); details are
				postponed to Remark \ref{remark-1/4}.
			\end{remark}
			
			\medskip
			
			All three versions of the Borell-Brascamp-Lieb inequality imply a corresponding Pr\'ekopa-Leindler-type inequality by simply setting $p = 0$ and using the convention $M_s^0(a,b) = a^{1-s}b^{s}$ for all $a,b \geq 0$ and $s \in (0,1)$; for sake of completeness we state them in the sequel.

			\begin{corollary}\label{C-weighted-prekopa-leindler} {\bf (Weighted Pr\'ekopa-Leindler inequality on $\mathbb H^n$)}
				Fix $s\in (0,1)$. 
				Let $f,g,h:\mathbb H^n\to [0,\infty)$ be Lebesgue integrable
				functions satisfying 
				$$
				h(z) \geq 
				\left(\frac{f(x)}{\left(\tilde{\tau}_{1-s}^n(\theta(y,x))\right)^{2n+1}}\right)^{1-s}\left(\frac{g(y)}{\left(\tilde{\tau}_s^n(\theta(x,y))\right)^{2n+1}} \right)^s \ \  for\ all\ (x,y)\in \Heis^n\times
				\Heis^n, z\in Z_s(x,y).
				$$
				Then the following inequality holds:
				\begin{eqnarray*}
					\int_{\Heis^n} h \geq \left(\int_{\Heis^n}
					f\right)^{1-s}\left( \int_{\Heis^n} g \right)^s.
				\end{eqnarray*}
			\end{corollary}
			
			\begin{corollary}\label{C-uniformlyweighted-prekopa-leindler} {\bf (Uniformly weighted Pr\'ekopa-Leindler inequality on $\mathbb H^n$)}
				Fix $s\in (0,1)$. 
				Let $f,g,h:\mathbb H^n\to [0,\infty)$ be Lebesgue integrable
				functions satisfying
				$$
				h(z) \geq 
				\left(\frac{f(x)}{(1-s)^2}\right)^{1-s}\left(\frac{g(y)}{s^2} \right)^s \ \  for\ all\ (x,y)\in \Heis^n\times
				\Heis^n, z\in Z_s(x,y).
				$$
				Then the following inequality holds:
				\begin{eqnarray*}
					\int_{\Heis^n} h \geq \left(\int_{\Heis^n}
					f\right)^{1-s}\left( \int_{\Heis^n} g \right)^s.
				\end{eqnarray*}
			\end{corollary}
			
			\begin{corollary}\label{C-nonweighted-prekopa-leindler} {\bf (Non-weighted Pr\'ekopa-Leindler inequality on $\mathbb H^n$)}
				Fix $s\in (0,1)$. 
				Let $f,g,h:\mathbb H^n\to [0,\infty)$ be Lebesgue integrable
				functions satisfying
				$$
				h(z) \geq 
				\left(f(x)\right)^{1-s}\left(g(y)\right)^s \ \  for\ all\ (x,y)\in \Heis^n\times
				\Heis^n, z\in Z_s(x,y).
				$$
				Then the following inequality holds:
				\begin{eqnarray*}
					\int_{\Heis^n} h \geq \frac{1}{4}\left(\int_{\Heis^n}
					f\right)^{1-s}\left( \int_{\Heis^n} g \right)^s.
				\end{eqnarray*}
			\end{corollary}
			
			\medskip
			
			Let us conclude this section with an observation. 
			\begin{remark} \rm
				It is possible to obtain a slightly improved version of the Borell-Brascamp-Lieb and Pr\'ekopa-Leindler-type inequalities by requiring that condition (\ref{ConditionRescaledBBLWithWeights}) holds only for $y=\psi(x)$ where $\psi$ is the optimal transport map between two appropriate absolutely continuous probability measures $\mu_0$ and $\mu_1$ given in terms of the densities $f$ and $g$. We leave the details to the interested reader.
			\end{remark}
			\medskip
			
			\section{Geometric aspects of Brunn-Minkowski inequalities on $\mathbb H^n$}\label{section-BM}  We first notice that different versions of the Brunn-Minkowski inequality have been studied  earlier in the setting of the Heisenberg group. In particular, Leonardi and Masnou \cite{LM} considered the  multiplicative Brunn-Minkowski inequality on $\mathbb H^n$, i.e., if $A,B\subset \mathbb H^n$ are compact sets, then 
			\begin{equation}\label{multiplicative-BM}
			\mathcal L^{2n+1}(A\cdot B)^\frac{1}{N}\geq \mathcal L^{2n+1}(A)^\frac{1}{N}+\mathcal L^{2n+1}( B)^\frac{1}{N}
			\end{equation}
			for some $N\geq 1,$ where $'\cdot'$ denotes the Heisenberg group law. It turned out that (\ref{multiplicative-BM}) fails for the homogeneous dimension $N=2n+2$, see Monti \cite{Monti};  moreover,  it fails even for all $N>2n+1$ as shown by Juillet \cite{Juillet-IMNR}. However,  inequality  (\ref{multiplicative-BM}) holds for the topological dimension $N=2n+1$, see \cite{LM}.
			
			In this subsection we shall present several geodesic Brunn-Minkowski inequalities on $\mathbb H^n$ and discuss their geometric features.
			
			\medskip
			
			\noindent {\bf Proof of Theorem \ref{Brunn-Minkowski}.} 
			We have nothing to prove when both sets have zero $\mathcal L^{2n+1}$-measure. 
			
			Let $A,B\subset \mathbb H^n$ be two nonempty measurable sets such that at least one of them has positive $\mathcal L^{2n+1}$-measure.  We first claim that $\Theta_{A,B}<2\pi$. To check this we recall that
			{$$\Theta_{A,B}=\sup_{A_0, B_0} \inf_{(x,y) \in A_0 \times B_0} \left\{|\theta|\in [0,2\pi]:(\chi,\theta)\in \Gamma_1^{-1}(x^{-1}\cdot
				y)\right\},$$
				where the sets $A_0$ and $B_0$ are nonempty, full measure subsets of $A$ and $B$, respectively.}
			Arguing by
			contradiction, if  $\Theta_{A,B}=2\pi$, it follows that up to a set of null $\mathcal L^{2n+1}$-measure, we have for
			every $(x,y)\in A\times B$ that
			$$x^{-1}\cdot y\in \Gamma_1(\chi,\pm 2\pi)\subset L=\{0_{\mathbb C^n}\}\times \mathbb R.$$
			In particular, up to a set of null $\mathcal L^{2n+1}$-measure, $A^{-1}\cdot B\subset \{0_{\mathbb C^n}\}\times
			\mathbb R,$ thus $\mathcal L^{2n+1}(A^{-1}\cdot B)=0.$ Therefore, the
			multiplicative Brunn-Minkowski inequality (\ref{multiplicative-BM}) for $N=2n+1$ gives that $$\mathcal
			L^{2n+1}(A^{-1}\cdot B)^\frac{1}{2n+1} \geq \mathcal
			L^{2n+1}(A^{-1})^\frac{1}{2n+1}+\mathcal
			L^{2n+1}(B)^\frac{1}{2n+1},$$ 
			which implies that $\mathcal L^{2n+1}(A)= \mathcal L^{2n+1}(B)=0$, a
			contradiction.
			
			{
				Fix $s\in (0,1)$ and let $$c_1^s=\sup_{A_0, B_0} \inf_{(x,y) \in A_0 \times B_0} \tilde \tau_{1-s}^n (\theta(y,x)) \ \ {\rm and}\ \ c_2^s=\sup_{A_0, B_0} \inf_{(x,y) \in A_0 \times B_0}\tilde \tau_{s}^n (\theta(x,y))
				,$$ where $A_0$ and $B_0$ are nonempty, full measure subsets of $A$ and $B$. 
				Since the function $\theta\mapsto
				\tilde \tau_s^n(\theta)$ is increasing on $[0,2\pi)$, cf. Lemma \ref{lemma-novekvo}, it turns out that 
				$$c_1^s=\tilde \tau_{1-s}^n(\Theta_{A,B}) \ \ {\rm and}\ \ c_2^s=\tilde \tau_s^n(\Theta_{A,B}).$$
				Due to the fact that $\Theta_{A,B}<2\pi$, we have  $0<c_1^s,c_2^s<+\infty.$} We now distinguish two cases. 
			
			{\bf Case 1:}  $\mathcal
			L^{2n+1}(A)\neq 0\neq \mathcal
			L^{2n+1}(B).$ Let $p=+\infty$, $f(x)=(c_1^s)^{2n+1} \mathbbm{1} _A(x),$
			$g(y)=(c_2^s)^{2n+1} \mathbbm{1} _B(y)$ and $h(z)=\mathbbm{1}
			_{Z_s(A,B)}(z).$ Since 
			(\ref{ConditionRescaledBBLWithWeights}) holds,
			we may apply  Theorem \ref{TRescaledBBLWithWeights} with the above choices, 
			obtaining
			\begin{eqnarray*}
				\mathcal L^{2n+1}(Z_s(A,B)) &\geq & M^\frac{1}{2n+1}_s \left((c_1^s)^{2n+1}\mathcal L^{2n+1}(A),
				(c_2^s)^{2n+1} \mathcal L^{2n+1}(B)\right) \\
				&=& \left(\tau_{1-s}^n(\Theta_{A,B})\mathcal
				L^{2n+1}(A)^\frac{1}{2n+1}+\tau_{s}^n(\Theta_{A,B})\mathcal
				L^{2n+1}(B)^\frac{1}{2n+1}\right)^{2n+1}.
			\end{eqnarray*}

			{\bf Case 2:} $\mathcal
			L^{2n+1}(A)\neq 0= \mathcal
			L^{2n+1}(B)$ or $\mathcal
			L^{2n+1}(A)= 0\neq \mathcal
			L^{2n+1}(B).$     We consider the first sub-case; the second one is treated in a similar way. 
			By the first part of the proof, we have that $\Theta_{A,B} < 2\pi$. By setting $\mu_0 = \frac{\mathcal L^{2n+1}|_{A}}{\mathcal L^{2n+1}(A)}$ and $\mu_1 = \delta_x$ the point-mass associated to a point $x \in B$, the Jacobian determinant inequality (\ref{Jacobi-degenerate}) can be explored in order to obtain 
			\begin{eqnarray*}
				\mathcal L^{2n+1}(Z_s(A,B)) &\geq& \mathcal L^{2n+1}(Z_s(A,\{x\})) = \mathcal L^{2n+1}\left(\cup_{y \in A}Z_s( y, x)\right) \geq \mathcal L^{2n+1}(\psi_s(A))  \\ 
				&=& \int\limits_A {\rm Jac} (\psi_s)(y) \dd \mathcal L^{2n+1}(y)
				\\ &\geq& \int\limits_A \left(\tau_{1-s}^n(\theta_y)\right)^{2n+1} \dd \mathcal L^{2n+1}(y) \geq \left(\tau_{1-s}^n(\Theta_{A,\{x\}})\right)^{2n+1} \mathcal L^{2n+1}(A) \\
				&\geq& \left(\tau_{1-s}^n(\Theta_{A,B})\right)^{2n+1} \mathcal L^{2n+1}(A),
			\end{eqnarray*}
			where we used that $\Theta_{A,\{x\}}\geq \Theta_{A,B}$. \hfill $\square$



			\begin{remark}\rm Let $\lambda>0$. Since $(\delta_\lambda (x))^{-1}\cdot \delta_\lambda (y)=\delta_\lambda (x^{-1}\cdot y)$ for every $x,y\in \mathbb H^n,$ it turns out that  $\Theta_{\delta_\lambda(A),\delta_\lambda(B)}=\Theta_{A,B}$ for every sets $A,B\subset \mathbb H^n$. As a consequence, the weighted Brunn-Minkowski inequality is  invariant under the dilation of the sets.
			\end{remark}
			
			The arguments in Theorem \ref{Brunn-Minkowski} put the {\it measure
				contraction property} $\textsf{MCP}(0,2n+3)$ of  Juillet \cite[Theorem 2.3]{Juillet-IMNR} into the right perspective. In particular, it explains the appearance of the somewhat mysterious value $2n+3$ of the exponent: 

			\begin{corollary}\label{MCP-1}  {\bf (Measure contraction property on $\mathbb H^n$)}
				The measure contraction property {\rm{\textsf{ MCP}}}$(0,2n+3)$ holds
				on $\mathbb H^n$, i.e., for every  $s\in [0,1]$, $x\in \mathbb H^n$
				and nonempty measurable set $E\subset \mathbb H^n$,
				$$\displaystyle \mathcal L^{2n+1}(Z_s(x,E)) \geq \left(\tau^n_s\left(\Theta_{\{x\},E}\right)\right)^{2n+1} \mathcal L^{2n+1}(E) \geq s^{2n+3}\mathcal L^{2n+1}(E).$$
			\end{corollary}
			{\it Proof.} The first inequality is nothing but the weighted Brunn-Minkowski inequality for $A = \{x\}$ and $B = E$ (see also Case 2 in the proof of Theorem \ref{Brunn-Minkowski}). 
			As $\tau^n_s \geq s^{\frac{2n+3}{2n+1}}$, the proof is complete.
			\hfill $\square$\\

			The geodesic Brunn-Minkowski inequality carries more information on the sub-Riemannian geometry of the Heisenberg group. {To illustrate this aspect, we give  the geometric interpretation of the expression $\Theta_{A,B}$ appearing in Theorem \ref{Brunn-Minkowski} for sets $A, B\subset \mathbb H^n$  with positive measure
				and of the Heisenberg distortion coefficients $\tau_{1-s}^n(\Theta_{A,B})$ and $\tau_s^n(\Theta_{A,B})$ that appear as weights in the 
				Brunn-Minkowski inequality.}
			
			We say that $A$ and $B$ are {\it essentially horizontal} if there exist full measure subsets $A_0 \subset A$ and $B_0 \subset B$ such that for every $x_0\in A_0$ there exists $y_0\in B_0\cap  H_{x_0},$ where $$H_{x_0}=\left\{y=(\zeta, t)\in \mathbb H^n: t =t_0+ 2
			{\rm Im} \langle\zeta_0 , {\zeta}\rangle\right\}$$ denotes the horizontal
			plane at $x_0=(\zeta_0, t_0).$ 
			In such a case, 
			for some $\chi_0\in \mathbb C^n$ we have $x_0^{-1}\cdot y_0 =(\chi_0,0)=\Gamma_1(\chi_0,0)$, i.e., $\Theta_{A,B} = 0$.

			We now turn our attention to the case when the sets $A$ and $B$ are not essentially horizontal to each other. Bellow we indicate an example showing that in such a case, the Heisenberg distortion coefficients
			$\tau_{1-s}^n(\Theta_{A,B})$ and $\tau_{s}^n(\Theta_{A,B})$ can even take arbitrarily large values. 
			
			To be more precise, let $s \in (0,1)$  and consider the CC-balls $A_r = B((0_{\mathbb C^n},t_1), r)$ and $B_r = B((0_{\mathbb C^n},t_2), r)$ in  $\mathbb H^n$ for sufficiently small values of $r>0$ and $t_1\neq t_2.$ Clearly, the sets $A_r$ and $B_r$ are horizontally far from each other, i.e., $B_r\cap H_{x_0}=\emptyset$ for every $x_0\in A_r.$  
			The geodesics joining the elements of $A_r$ and
			$B_r$ largely deviate from the $t$-axis and $Z_s(A_r,B_r)$ becomes a
			large set w.r.t. $A_r$ and $B_r$; see Figure \ref{figure-2} for $n=1$. More precisely, we have
			\begin{figure}
				\centering
				\includegraphics[scale=0.45]{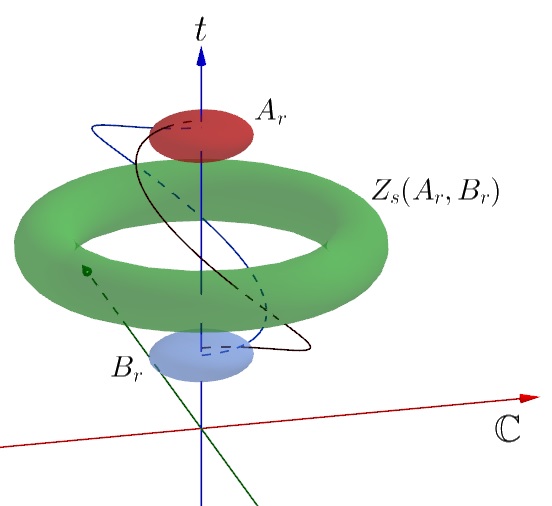}
				\includegraphics[scale=0.45]{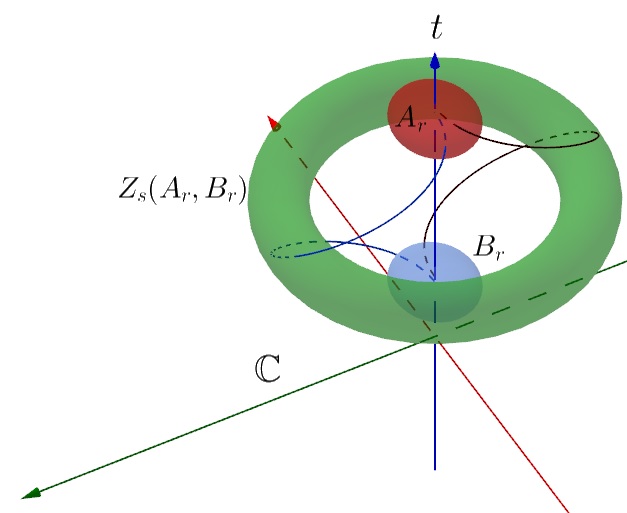}
				\caption{Heisenberg geodesics in $\mathbb H^1$ viewed from two different positions joining points in $A_r$ and $B_r$, and the  set $Z_s(A_r,B_r)$ of $s$-intermediate points.}\label{figure-2}
			\end{figure}

			\begin{proposition}\label{proposition-estimate} Let $A_r = B(0_{\Heis^n}, r)$, $B_r = B((0_{\mathbb C^n},1), r)$  and $s\in (0,1)$. Then
				\begin{itemize}
					\item[{\rm(i)}] $\Lmeas^{2n+1}(Z_s(A_r,B_r))=\omega(r^{3})  \mbox{ as } r\to  0;$\footnotemark
					\item[{\rm(ii)}] $2\pi-\Theta_{A_r,B_r}=\mathcal O(r) \mbox{ as }  r\to  0; $
					\item[{\rm(iii)}]  $\tau_{s}^n(\Theta_{A_r,B_r})=\omega\left(r^\frac{1-2n}{1+2n}\right) \mbox{ as }  r\to  0. $
				\end{itemize}	
				\addtocounter{footnote}{0}
				\footnotetext{$f(r)=\omega(g(r))$ as $r\to 0$ if there exist
					$c,\delta>0$ such that $|f(r)|\geq c|g(r)|$ for every $r\in
					(0,\delta)$.} 	
			\end{proposition}

			{\it Proof.} (i) Note first that for every $r>0$ one has
			$\mathcal
			L^{2n+1}(A_r)=\mathcal L^{2n+1}(B_r)=c_0r^{2n+2}$  for some  $c_0>0.$ By using the same notations as in (\ref{kesobb-jo-lesz}), it yields
			\begin{eqnarray*}
				\Lmeas^{2n+1}(Z_s(A_r,B_r)) &\geq&\Lmeas^{2n+1}(Z_s(0_{\mathbb H^n},B_r)) \geq\frac{c_2}{c_1^{2n-1}}\frac{\Lmeas^{2n+1}(B((0_{\mathbb C^n},1),r)\setminus L)}{r^{2n-1}}=\frac{c_2}{c_1^{2n-1}}\frac{c_0r^{2n+2}}{r^{2n-1}}\\
				&=& c_3 r^3,
			\end{eqnarray*}
			as $r\to 0$, where $c_3>0$ depends on $n\in \mathbb N$ and $s\in (0,1)$.

			(ii)\&(iii) By elementary behavior of the Heisenberg geodesics (\ref{explicit-geodet}) it follows that $\Theta_{A_r,B_r}\to 2\pi$ as $r\to 0.$ In fact, a similar estimate as in (\ref{kesobb-kell-peldaba}) shows that $\sin(\Theta_{A_r,B_r}/2)=\mathcal O(r)$ as $r\to 0$, which implies that  $2\pi-\Theta_{A_r,B_r}=\mathcal O(r)$ as $r\to 0$. 
			By the latter estimate and (\ref{concentration}) we have that 
			$\tau_{s}^n(\Theta_{A_r,B_r})\geq c_4 r^\frac{1-2n}{2n+1}$
			as $r\to 0,$ where $c_4>0$ depends on $n\in \mathbb N$ and $s\in (0,1)$.  $\hfill \square$
			\\


			\noindent In particular, Proposition \ref{proposition-estimate} implies that  $$\frac{\mathcal
				L^{2n+1}(Z_s(A_r,B_r))}{\mathcal L^{2n+1}(A_r)}\to +\infty\ \ {\rm as}\ \ r\to
			0;$$ this is the reason why the weights
			$\tau_{1-s}^n(\Theta_{A_r,B_r})$ and $\tau_{s}^n(\Theta_{A_r,B_r})$
			appear in the  Brunn-Minkowski inequality (\ref{BM-1}) in order to compensate the size of
			$Z_s(A_r,B_r)$ w.r.t. $A_r$ and $B_r$.  Quantitatively, the left hand side of  (\ref{BM-1}) is 
			$$\mathcal
			L^{2n+1}(Z_s(A_r,B_r))^\frac{1}{2n+1}=\omega\left(r^\frac{3}{2n+1}\right),$$ while the right hand side has the growth
			$$\tau_{s}^n(\Theta_{A_r,B_r})\mathcal
			L^{2n+1}(A_r)^\frac{1}{2n+1}=\omega\left(r^\frac{1-2n}{1+2n}\right)r^\frac{2n+2}{2n+1}=\omega\left(r^\frac{3}{1+2n}\right)$$ as $r\to 0$, which is in a perfect concordance with the
			competition of the two sides of (\ref{BM-1}). \\

			A direct consequence of  Theorem \ref{Brunn-Minkowski}, Corollary \ref{CRescaledBBLWithoutWeights} and
			estimate (\ref{tau-becsles}) reads as follows:
			
			\begin{corollary}\label{CBrunn-Minkowski-2}  {\bf (Non-weighted Brunn-Minkowski inequalities on $\mathbb H^n$)}
				Let $s\in (0,1)$ and $A$ and $B$ be two nonempty measurable sets of $\mathbb
				H^n$. Then the following inequalities hold:
				\begin{itemize}
					\item[{\rm (i)}] $\displaystyle \mathcal L^{2n+1}(Z_s(A,B))^\frac{1}{2n+1} \geq
					(1-s)^\frac{2n+3}{2n+1}\mathcal
					L^{2n+1}(A)^\frac{1}{2n+1}+s^\frac{2n+3}{2n+1}\mathcal
					L^{2n+1}(B)^\frac{1}{2n+1};$
					\item[{\rm (ii)}] $\displaystyle \mathcal L^{2n+1}(Z_s(A,B))^\frac{1}{2n+3} \geq
					\left(\frac{1}{4}\right)^{\frac{1}{2n+3}} \left( (1-s)\mathcal
					L^{2n+1}(A)^\frac{1}{2n+3}+s\mathcal
					L^{2n+1}(B)^\frac{1}{2n+3}\right).$
				\end{itemize}
			\end{corollary}

			The other main result of Juillet \cite[Lemma 3.1]{Juillet-IMNR} and Corollary \ref{MCP-1}  implicitly show that the non-weighted
			Brunn-Minkowski inequalities on $\mathbb H^n$ (see Corollary
			\ref{CBrunn-Minkowski-2})  are sharp.
			
			\begin{remark}\label{remark-1/4} \rm {\bf (Optimality of the constant $\frac{1}{4}$ in Corollaries \ref{CRescaledBBLWithoutWeights} and  \ref{CBrunn-Minkowski-2} (ii))} We deal just
				with Corollary \ref{CBrunn-Minkowski-2} (ii).  Let us assume that we can
				put a larger value instead of $\frac{1}{4}$ in our conclusion, i.e.,
				$\frac{1}{4}+\eta$ with $\eta>0.$ Let $A=A_r$ and $B=B_r$ be the
				Euclidean balls of radius $r$ and centers $a=(-1,0,...,0)\in \mathbb H^n$ and
				$b=(1,0,...,0)\in \mathbb H^n$, respectively. According to our hypothesis, 
				$$\mathcal L^{2n+1}(Z_{1/2}(A_r,B_r))\geq  \left(\frac{1}{4}+\eta\right)\mathcal
				L^{2n+1}(B_r).$$ On the other hand, by Juillet \cite[Theorem
				1]{Juillet-kicsi} and relation (\ref{Jacobian-Juillet}) one has that
				$$\limsup_{r\to 0}\frac{\mathcal L^{2n+1}(Z_{1/2}(A_r,B_r))}{\mathcal L^{2n+1}(B_r)}\leq 2^{2n+1}\frac{ \Jac (a\cdot\Gamma_{1/2})}{\Jac
					(a\cdot\Gamma_1)}(\Gamma_1^{-1}(a^{-1}\cdot
				b))=2^{2n+1}\cdot\frac{1}{2^{2n+3}}=\frac{1}{4},$$ a contradiction.
			\end{remark}
			
			\begin{remark} \rm
				We notice that instead of  $\Theta_{A,B}$ in the weighted Brunn-Minkowski inequality, we can use a better quantity depending on the optimal mass transport
				$$\hat{\Theta}_{A, \psi}= \sup_{A_0} \inf_{x \in A_0} \left\{|\theta|\in [0,2\pi]:(\chi,\theta)\in \Gamma_1^{-1}(x^{-1}\cdot \psi(x))\right\},$$
				where the set $A_0$ is a nonempty, full measure subset of $A$ and $\psi$ is the optimal transport map resulting from the context. Since $\hat{\Theta}_{A, \psi} \geq\Theta_{A,B}$ and  $\tau_s^n$ is increasing, one has  $\tau_s^n\left(\hat{\Theta}_{A,\psi}\right)\geq \tau_s^n(\Theta_{A,B})$. In this way, one can slightly improve the Brunn-Minkowski inequality (\ref{BM-1}). A further improvement can be obtained by replacing $\tau^n_s$ by $\hat{\tau}^n_s$ from (\ref{better-concentration}).
			\end{remark}

			\medskip
			
			\section{Concluding remarks and further questions}\label{SectionFinalRemarks}
			
			The purpose of this final section is to indicate open research  
			problems that  are closely related to our results and can be considered as starting points of further investigations.
			
			Let us mention first that there have been several different approaches to functional inequalities for sub-Riemannian geometries. One such possibility was initiated  by Baudoin, Bonnefont and Garofalo \cite{Bau-Bon-Gar} via the Bakry-\'Emery  {\it carr\'e du champ} operator by introducing an analytic  curvature-dimension inequality on sub-Riemannian manifolds. A challenging problem is to establish the relationship between their and our results, similarly as Erbar,  Kuwada and Sturm \cite{EKS} performed recently by proving the equivalence of the entropic curvature-dimension condition and Bochner's inequality formulated in terms of the Bakry-\'Emery   operator on metric measure spaces.

			One of the standard proofs of the isoperimetric inequality in $\mathbb R^n$ is  based on the  Brunn-Minkowski inequality. In 1982, Pansu \cite{Pansu} conjectured that the extremal set in the isoperimetric inequality  in $\mathbb H^1$ is the so-called bubble set. This is a topological ball whose boundary is foliated by geodesics. In our notation, the bubble sphere can be given as $\{ \Gamma^1_s(\chi, 2\pi): |\chi|=1, s\in [0,1] \}$.  Although there are several partial answers to this question supporting the conjecture (under $C^2$-smoothness or axially-symmetry of domains), the general case is still unsolved; see the monograph of Capogna, Danielli, Pauls and Tyson \cite{Capogna}.  
			We believe that our  Brunn-Minkowski inequality (e.g. Theorem \ref{Brunn-Minkowski}) could provide a new approach to Pansu's conjecture; a deeper understanding of the behavior of the optimal transport map is indispensable.

			Closely related to isoperimetric inequalities are sharp Sobolev inequalities. The method of optimal mass transportation is an efficient tool to prove such results, see Cordero-Erausquin,  Nazaret and Villani \cite{CE-N-Villani} and Villani \cite[Chapter 6]{Villani}. Moreover,   Bobkov and Ledoux \cite{BL1, BL2} established sharp Sobolev-type inequalities on $\mathbb R^n$ by using a version of the  Brunn-Minkowski inequality and properties of the solutions of the Hamilton-Jacobi equation given by the infimum-convolution operator. Since the latter is well understood on $(\mathbb H^n,d_{CC},\mathcal
			L^{2n+1})$, see Manfredi and Stroffolini \cite{MS}, it seems plausible to approach sharp Sobolev inequalities in $\mathbb H^n$ by the Brunn-Minkowski inequality (\ref{BM-1}).  We note that Frank and Lieb \cite{FL} obtained recently sharp Hardy-Littlewood-Sobolev-type inequalities on $\mathbb H^n$ by a careful analysis of a convolution operator. Sharp Sobolev inequalities with general exponents are still open in the Heisenberg group.

			We expect that our method will open the possibility to study geometric inequalities on generic Sasakian manifolds verifying a lower bound assumption for the Ricci curvature. If $(M,d_{SR},\mu)$ is a $2n+1$ dimensional Sasakian manifold equipped with a natural sub-Riemannian structure, where $d_{SR}$ is the sub-Riemannian distance and $\mu$ is the corresponding Riemannian volume form on $M$,  the Ricci curvature lower bound is  formulated by controlling from below the tangent vectors 
			from the canonical distribution in terms of the Tanaka-Webster connection. This notion requires two parameters, $k_1,k_2\in \mathbb R$, depending on the specific components of the vectors from the distribution, see Lee 	\cite{Lee-arxiv}, Lee, Li and Zelenko \cite{LLZ}, and Agrachev and Lee \cite{AL} for $n=1$. In \cite{LLZ} it is
			proved that a $2n+1$ dimensional Sasakian manifold with $(k_1,k_2)$ Ricci curvature lower bound satisfies the generalized measure contraction property $\mathcal{MCP}(k_1,k_2,2n,2n+1)$. If $(M,d_{SR},\mu)=(\mathbb H^n,d_{CC},\mathcal
			L^{2n+1})$, it turns out that $\mathcal{MCP}(0,0,2n,2n+1)=${\rm{\textsf{ MCP}}}$(0,2n+3)$. Note that the  Heisenberg group $\mathbb H^n$ is the  simplest Sasakian manifold with vanishing Tanaka-Webster curvature, in a similar way as the Euclidean space $\mathbb R^n$ is the standard flat space among $n$-dimensional Riemannian manifolds.  
			It would be interesting to extend the results from our paper to this more general setting. We expect that by direct computations one can determine explicit forms of the Sasakian distortion coefficient $\tau_s^{k_1,k_2,n}$ which should reduce to the Heisenberg distortion coefficient $\tau_s^n$ whenever $M=\mathbb H^n$ (and $k_1=k_2=0$).  
			
		In order to avoid further technical difficulties, in the present paper we focused to the Heisenberg groups $\mathbb H^n$. Note that our method also works on Carnot groups of step two, on the $3$-sphere or on more general sub-Riemannian manifolds which have well-behaving cut-locus, see e.g. Boscain and Rossi \cite{BR}, Rifford \cite{R1, R2} and Rizzi \cite{Rizzi}.  Contrary to groups of step two, the structure of sub-Riemannian cut-locus in generic sub-Riemannian manifolds may have a pathological behavior, see e.g. Figalli and Rifford \cite[\S 5.8, p. 145]{FR}, and the geodesics in the Riemannian approximants may converge to singular geodesics. 

			After posting the first version of the present work to the mathematical community, follow-up works have been obtained by establishing {\it intrinsic} geometric inequalities on corank 1 Carnot groups (by Balogh, Krist\'aly and Sipos \cite{BKS-Jacobian}) and on ideal sub-Riemannian manifolds (by Barilari and Rizzi \cite{Barilari-Rizzi}) by different methods than ours. Naturally, the Heisenberg distortion coefficient $\tau_s^n$ introduced in the present paper and those from  the latter works \textit{coincide} on $\mathbb H^n$. This confirms the efficiency of the approximation arguments in suitable sub-Riemannian geometric contexts. In addition, as C. Villani suggested in \cite[p. 43]{Villani-cikk}, the results in the present paper (together with those from \cite{BKS-Jacobian} and \cite{Barilari-Rizzi}) motivate  the so-called "grande unification" of  geometric inequalities appearing in  Riemannian, Finslerian and sub-Riemannian geometries.

\nocite{*}
\bibliographystyle{plain}
\bibliography{references}

\begin{thebibliography}{10}

\bibitem{AL}
A.~Agrachev and P.~W.~Y. Lee.
\newblock Generalized {R}icci curvature bounds for three dimensional contact
  subriemannian manifolds.
\newblock {\em Math. Ann.}, 360(1-2):209--253, 2014.

\bibitem{AR}
L.~Ambrosio and S.~Rigot.
\newblock Optimal mass transportation in the {H}eisenberg group.
\newblock {\em J. Funct. Anal.}, 208(2):261--301, 2004.

\bibitem{Bacher}
K.~Bacher.
\newblock On {B}orell-{B}rascamp-{L}ieb inequalities on metric measure spaces.
\newblock {\em Potential Anal.}, 33(1):1--15, 2010.

\bibitem{BKS-CR}
Z.~M. Balogh, A.~Krist\'aly, and K.~Sipos.
\newblock Geodesic interpolation inequalities on {H}eisenberg groups.
\newblock {\em C. R. Math. Acad. Sci. Paris}, 354(9):916--919, 2016.

\bibitem{BKS-Jacobian}
Z.~M. {Balogh}, A.~{Krist{\'a}ly}, and K.~{Sipos}.
\newblock Jacobian determinant inequality on corank 1 {C}arnot groups with
  applications.
\newblock {\em ArXiv e-print: 1701.08831}, pages 1--28, January 2017.

\bibitem{Barilari-Rizzi}
D.~Barilari and L.~Rizzi.
\newblock Sub-riemannian interpolation inequalities: ideal structures.
\newblock {\em ArXiv e-print: 1705.05380}, pages 1--45, May 2017.

\bibitem{Bau-Bon-Gar}
F.~Baudoin, M.~Bonnefont, and N.~Garofalo.
\newblock A sub-{R}iemannian curvature-dimension inequality, volume doubling
  property and the {P}oincar\'e inequality.
\newblock {\em Math. Ann.}, 358(3-4):833--860, 2014.

\bibitem{BL1}
S.~G. Bobkov and M.~Ledoux.
\newblock From {B}runn-{M}inkowski to {B}rascamp-{L}ieb and to logarithmic
  {S}obolev inequalities.
\newblock {\em Geom. Funct. Anal.}, 10(5):1028--1052, 2000.

\bibitem{BL2}
S.~G. Bobkov and M.~Ledoux.
\newblock From {B}runn-{M}inkowski to sharp {S}obolev inequalities.
\newblock {\em Ann. Mat. Pura Appl. (4)}, 187(3):369--384, 2008.

\bibitem{BR}
U.~Boscain and F.~Rossi.
\newblock Invariant {C}arnot-{C}aratheodory metrics on {$S^3,\ {\rm SO}(3),\
  {\rm SL}(2)$}, and lens spaces.
\newblock {\em SIAM J. Control Optim.}, 47(4):1851--1878, 2008.

\bibitem{Capogna}
L.~Capogna, D.~Danielli, S.~D. Pauls, and J.~T. Tyson.
\newblock {\em An introduction to the {H}eisenberg group and the
  sub-{R}iemannian isoperimetric problem}, volume 259 of {\em Progress in
  Mathematics}.
\newblock Birkh\"auser Verlag, Basel, 2007.

\bibitem{McCann}
D.~Cordero-Erausquin, R.~J. McCann, and M.~Schmuckenschl{\"a}ger.
\newblock A {R}iemannian interpolation inequality \`a la {B}orell, {B}rascamp
  and {L}ieb.
\newblock {\em Invent. Math.}, 146(2):219--257, 2001.

\bibitem{CE-N-Villani}
D.~Cordero-Erausquin, B.~Nazaret, and C.~Villani.
\newblock A mass-transportation approach to sharp {S}obolev and
  {G}agliardo-{N}irenberg inequalities.
\newblock {\em Adv. Math.}, 182(2):307--332, 2004.

\bibitem{EKS}
M.~Erbar, K.~Kuwada, and K.-T. Sturm.
\newblock On the equivalence of the entropic curvature-dimension condition and
  {B}ochner's inequality on metric measure spaces.
\newblock {\em Invent. Math.}, 201(3):993--1071, 2015.

\bibitem{FJ}
A.~Figalli and N.~Juillet.
\newblock Absolute continuity of {W}asserstein geodesics in the {H}eisenberg
  group.
\newblock {\em J. Funct. Anal.}, 255(1):133--141, 2008.

\bibitem{FR}
A.~Figalli and L.~Rifford.
\newblock Mass transportation on sub-{R}iemannian manifolds.
\newblock {\em Geom. Funct. Anal.}, 20(1):124--159, 2010.

\bibitem{FL}
R.~L. Frank and E.~H. Lieb.
\newblock Sharp constants in several inequalities on the {H}eisenberg group.
\newblock {\em Ann. of Math. (2)}, 176(1):349--381, 2012.

\bibitem{GHL}
S.~Gallot, D.~Hulin, and J.~Lafontaine.
\newblock {\em Riemannian geometry}.
\newblock Universitext. Springer-Verlag, Berlin, 1987.

\bibitem{Gardner}
R.~J. Gardner.
\newblock The {B}runn-{M}inkowski inequality.
\newblock {\em Bull. Amer. Math. Soc. (N.S.)}, 39(3):355--405, 2002.

\bibitem{Gromov}
M.~Gromov.
\newblock Carnot-{C}arath\'eodory spaces seen from within.
\newblock In {\em Sub-{R}iemannian geometry}, volume 144 of {\em Progr. Math.},
  pages 79--323. Birkh\"auser, Basel, 1996.

\bibitem{Juillet-IMNR}
N.~Juillet.
\newblock Geometric inequalities and generalized {R}icci bounds in the
  {H}eisenberg group.
\newblock {\em Int. Math. Res. Not. IMRN}, (13):2347--2373, 2009.

\bibitem{Juillet-kicsi}
N.~Juillet.
\newblock On a method to disprove generalized {B}runn-{M}inkowski inequalities.
\newblock In {\em Probabilistic approach to geometry}, volume~57 of {\em Adv.
  Stud. Pure Math.}, pages 189--198. Math. Soc. Japan, Tokyo, 2010.

\bibitem{Juillet-calculus}
N.~Juillet.
\newblock Diffusion by optimal transport in {H}eisenberg groups.
\newblock {\em Calc. Var. Partial Differential Equations}, 50(3-4):693--721,
  2014.

\bibitem{Lee-arxiv}
P.~W.~Y. {Lee}.
\newblock {Ricci curvature lower bounds on Sasakian manifolds}.
\newblock {\em ArXiv e-print: 1511.09381}, November 2015.

\bibitem{LLZ}
P.~W.~Y. Lee, C.~Li, and I.~Zelenko.
\newblock Ricci curvature type lower bounds for sub-{R}iemannian structures on
  {S}asakian manifolds.
\newblock {\em Discrete Contin. Dyn. Syst.}, 36(1):303--321, 2016.

\bibitem{LM}
G.~P. Leonardi and S.~Masnou.
\newblock On the isoperimetric problem in the {H}eisenberg group {${\Bbb
  H}^n$}.
\newblock {\em Ann. Mat. Pura Appl. (4)}, 184(4):533--553, 2005.

\bibitem{LV}
J.~Lott and C.~Villani.
\newblock Ricci curvature for metric-measure spaces via optimal transport.
\newblock {\em Ann. of Math. (2)}, 169(3):903--991, 2009.

\bibitem{MS}
J.~Manfredi and B.~Stroffolini.
\newblock A version of the {H}opf-{L}ax formula in the {H}eisenberg group.
\newblock {\em Comm. Partial Differential Equations}, 27(5-6):1139--1159, 2002.

\bibitem{McCann-PhD}
R.~J. McCann.
\newblock {\em A convexity theory for interacting gases and equilibrium
  crystals}.
\newblock ProQuest LLC, Ann Arbor, MI, 1994.
\newblock Thesis (Ph.D.)--Princeton University.

\bibitem{McCann_Adv_Math}
R.~J. McCann.
\newblock A convexity principle for interacting gases.
\newblock {\em Adv. Math.}, 128(1):153--179, 1997.

\bibitem{McCann-GAFA}
R.~J. McCann.
\newblock Polar factorization of maps on {R}iemannian manifolds.
\newblock {\em Geom. Funct. Anal.}, 11(3):589--608, 2001.

\bibitem{Monti}
R.~Monti.
\newblock Brunn-{M}inkowski and isoperimetric inequality in the {H}eisenberg
  group.
\newblock {\em Ann. Acad. Sci. Fenn. Math.}, 28(1):99--109, 2003.

\bibitem{Ohta}
S.~Ohta.
\newblock Finsler interpolation inequalities.
\newblock {\em Calc. Var. Partial Differential Equations}, 36(2):211--249,
  2009.

\bibitem{Pansu}
P.~Pansu.
\newblock Une in\'egalit\'e isop\'erim\'etrique sur le groupe de {H}eisenberg.
\newblock {\em C. R. Acad. Sci. Paris S\'er. I Math.}, 295(2):127--130, 1982.

\bibitem{R2}
L.~Rifford.
\newblock Ricci curvatures in {C}arnot groups.
\newblock {\em Math. Control Relat. Fields}, 3(4):467--487, 2013.

\bibitem{R1}
L.~Rifford.
\newblock {\em Sub-{R}iemannian geometry and optimal transport}.
\newblock Springer Briefs in Mathematics. Springer, Cham, 2014.

\bibitem{Rizzi}
L.~Rizzi.
\newblock Measure contraction properties in {C}arnot groups.
\newblock {\em Calc. Var.}, 55(3, DOI: 10.1007/s00526-016-1002-y), 2016.

\bibitem{Sturm1}
K.-T. Sturm.
\newblock On the geometry of metric measure spaces. {I}.
\newblock {\em Acta Math.}, 196(1):65--131, 2006.

\bibitem{Sturm2}
K.-T. Sturm.
\newblock On the geometry of metric measure spaces. {II}.
\newblock {\em Acta Math.}, 196(1):133--177, 2006.

\bibitem{Villani}
C.~Villani.
\newblock {\em Topics in optimal transportation}, volume~58 of {\em Graduate
  Studies in Mathematics}.
\newblock American Mathematical Society, Providence, RI, 2003.

\bibitem{Villani1}
C.~Villani.
\newblock {\em Optimal transport}, volume 338 of {\em Grundlehren der
  Mathematischen Wissenschaften [Fundamental Principles of Mathematical
  Sciences]}.
\newblock Springer-Verlag, Berlin, 2009.
\newblock Old and new.

\bibitem{Villani-cikk}
C.~Villani.
\newblock In\'egalit\'es isop\'erim\'etriques dans les espaces metriques
  mesur\'es.
\newblock {\em S\'eminaire Bourbaki, Ast\'erisque}, 69\`eme
  ann\'ee(1127):1--50, 2017.

\end{thebibliography}

\vspace{0.5cm} \noindent {\footnotesize{\sc  Mathematisches Institute,
Universit\"at Bern,
             Sidlerstrasse 5,
 3012 Bern, Switzerland.}\\
 Email: \textsf{zoltan.balogh@math.unibe.ch}\\

 \noindent {\footnotesize {\sc Department of Economics, Babe\c s-Bolyai University, Str. Teodor Mihali 58-60, 400591
 		Cluj-Napoca, Romania \& Institute of Applied Mathematics, \'Obuda University,
  B\'ecsi \'ut 96, 1034 Budapest, Hungary.}\\ Email:
{\textsf{alex.kristaly@econ.ubbcluj.ro}}\\

 \noindent {\footnotesize{\sc   Mathematisches Institute,
Universit\"at Bern,
             Sidlerstrasse 5,
 3012 Bern, Switzerland.\\
} Email: {\textsf{kinga.sipos@math.unibe.ch}\\

\end{document}